\documentclass[a4paper]{article}
\usepackage[utf8]{inputenc}
\usepackage{mathtools}
\usepackage{amsthm}
\usepackage{amssymb}
\usepackage{parskip}
\usepackage{natbib}
\usepackage{hyperref}
\usepackage{xcolor}
\usepackage{graphicx, wrapfig, subcaption}
\usepackage{epstopdf}
\usepackage[a4paper, margin=2.5cm]{geometry}
\usepackage{placeins}
\usepackage{authblk}

\newcommand{\adj}[0]{\text{adj}}
\newcommand{\diag}[0]{\text{diag}}
\newcommand{\br}[0]{\mathbf{r}} % def. bold r
\newcommand{\re}[0]{\text{Re}}
\newcommand{\im}[0]{\text{Im}}

\newtheorem{thm}{Theorem}[section]

\title{Bifurcations of Neural Fields on the Sphere}
\author[1,*]{Len Spek}
\author[1]{Stephan A. van Gils}
\author[1,2]{Yuri A. Kuznetsov}
\author[3]{M\'onika Polner}
\affil[1]{Department of Applied Mathematics, University of Twente, Enschede, The Netherlands}
\affil[2]{Department of Mathematics, Utrecht University, Utrecht, The Netherlands}
\affil[3]{Bolyai Institute, University of Szeged, Szeged, Hungary}
\affil[*]{Corresponding author: l.spek@utwente.nl}
\date{December 2022}

\begin{document}
\maketitle

\begin{abstract}
    \noindent A common model to study pattern formation in large groups of neurons is the neural field. We investigate a neural field with excitatory and inhibitory neurons, like \cite{wilson_excitatory_1972}, with transmission delays and gap junctions. We build on the work of \cite{visser_standing_2017} by investigating pattern formation in these models on the sphere. Specifically, we investigate how different amounts of gap-junctions, modelled by a diffusion term, influence the behaviour of the neural field. We look in detail at the periodic orbits which are generated by Hopf bifurcation in the presence of spherical symmetry. For this end, we derive general formulas to compute the normal form coefficients of these bifurcations up to third order and predict the stability of the resulting branches and formulate a novel numerical method to solve delay equations on the sphere. 
    \end{abstract}

\section{Introduction}
Neural fields are a useful model to understand the macroscopic spatiotemporal behaviour of the brain, which can be measured by an electroencephalogram (EEG) \cite{jirsa_field_1996}. Synchronised waves of neural activity observed in patients with Parkinson's disease have been linked with an increase of gap junctions, electrical connections in the brain \cite{schwab_pallidal_2014}. Neural fields offer a natural way to study this link \cite{spek_neural_2020, spek_dynamics_2022}. The aforementioned works considered neural fields on an interval and rectangle, respectively. Boundary conditions however, play a critical role in the pattern formation in neural fields \cite{ingber_neocortical_2011}, hence a spherical domain is a more natural topology when considering activity in the brain \cite{jirsa_spatiotemporal_2002}. 

This work is a generalisation of the paper by \cite{visser_standing_2017} on neural fields with transmission delays on spherical domains. We expand these results in several directions: we add a diffusion term to model gap junctions. Moreover, we consider two distinct populations of excitatory and inhibitory neurons, similar to the models of Wilson and Cowan \cite{wilson_mathematical_1973}, in contrast to the effective single population models of Amari \cite{amari_dynamics_1977}. We also expand on the analysis of the Hopf-bifurcations with spherical symmetry, particularly, investigate the equivariant normal forms for $l=2,3$, \cite{iooss_hopf_1989,sigrist_hopf_2010}. 

We find a non-trivial interaction between the delay induced oscillations and the diffusion term. Specifically, when the diffusion coefficient between the excitatory neurons is significantly higher than between inhibitory neurons, periodic orbits with higher order harmonics become stable, compared to the cases where the diffusion coefficients are equal or zero. This is in contrast to the findings in previous less detailed models \cite{spek_neural_2020}, where a large diffusion coefficient destabilises all spatially non-uniform patterns. 

Our main contributions are the derivation of analytical conditions which determine when the trivial equilibrium loses stability due to a bifurcation, a systematic method to obtain analytical formulas for the normal form coefficients in terms of the parameters of the model, and a numerical discretisation which is suitable to solve reaction diffusion equations on a sphere where the reaction terms contain distributed delays. These contributions are specifically developed for neural fields, but can be applied to many models with spherical symmetry. 

\subsection{Background}
Wilson and Cowan introduced, in the 1970s, neural fields to study the interactions of large collections of neurons in a mean-field setting \cite{wilson_excitatory_1972,wilson_mathematical_1973}. They considered two interacting populations of neurons, excitatory and inhibitory, located in some abstract space. Amari reduced this model by combining the excitatory and inhibitory neurons into a single field \cite{amari_dynamics_1977}. Nunez and collaborators used these models to show how the long-range connections are essential in generating the $\alpha$-rhythm of an EEG \cite{nunez_brain_1974, nunez_electric_2006}. They also recognised the importance of choosing the proper boundary conditions for the model \cite{ingber_neocortical_2011}. The cortex is topologically close to a sphere, and neural fields naturally produce standing waves in such topologies without a boundary \cite{ingber_neocortical_2011}. Visser et al. investigated a neural field model with transmission delays on a spherical domain \cite{visser_standing_2017}. They used bifurcation theory in the presence of symmetry to study the formation of spatio-temporal patterns. They observed that the standing and rotating waves found with this analysis have counterparts in the EEG patterns of epileptic and schizophrenic brain conditions. 

Gap junctions are electrical connections between neurons, which directly exchange ions through a connexin protein, in contrast to synaptic connections where neurotransmitters are involved. In terms of neural fields, they can be modelled as a simple diffusion process \cite{coombes_tutorial_2014} and have been studied in \cite{spek_neural_2020}. They can have an impact on the activity of large collections of neurons, for example they are thought to be related to certain anomalous brain activity in patients with Parkinson's disease \cite{schwab_pallidal_2014,schwab_synchronization_2014}.

When analysing spatio-temporal patterns in dynamical systems, the underlying symmetries of the model and the domain on which the model is defined, play an important role. Bressloff, Cowan and collaborators have for example used this in their study of visual hallucinations, which arise from a Turing instability in a neural field model with the Euclidean symmetry group  \cite{ermentrout_mathematical_1979,bressloff_geometric_2001,bressloff_what_2002}. This symmetry group contains all the translations and rotations of free space. 

Equivariant bifurcation theory is the systematic study of bifurcations in symmetric dynamical systems. The book by Golubitsky et al. offers a detailed overview of the theory and its applications \cite{golubitsky_singularities_1988}. A key concept here is spontaneous symmetry breaking, where a state that is invariant with respect to a symmetry group loses stability. In turn, possibly multiple, branches appear which have less symmetry than the original state. For Hopf bifurcations with spherical symmetry, \cite{golubitsky_hopf_1985} gives an overview of the different possible branches. The representation of the Lie group of spherical symmetries, the vector space of orthogonal $3\times 3$ matrices, can be decomposed as the infinite sum of irreducible subspaces $V_l$. Each of these $V_l$ are spanned by the spherical harmonics of order $l$. The equivariant Hopf theorem guarantees that for every maximal isotropy subgroup $\Sigma$, there exists a unique branch of periodic solutions bifurcating from the origin with $\Sigma$ as their group of symmetries. For $l=2$, Iooss and Rossi investigated the stability of these branches \cite{iooss_hopf_1989} and for $l=3$ this analysis was recently performed by Sigrist \cite{sigrist_hopf_2010,sigrist_symmetric_2011}.

In this paper, we study neural fields with transmission delays. For bifurcation analysis, these delays require the study of our model in the context of abstract delay differential equations. A framework which allows a natural and rigorous foundation of the bifurcation theory of delay equations is the sun-star calculus \cite{diekmann_delay_1995}. As our model combines diffusion with delays, a recent generalisation of the sun-star calculus in \cite{janssens_class_2019,janssens_class_2020} is required. This generalisation has been applied to Hopf bifurcations in \cite{spek_neural_2020}, hence we use the results and convenient notations of this paper to derive formulas for the normal form coefficients of Hopf bifurcations with spherical symmetry.

Detail descriptions of numerical methods for time integration of neural fields with transmission delays are rare in the literature, e.g. see  \cite{faye_theoretical_2010, lima_numerical_2015,hutt_numerical_2014,polner_space-time_2017}. The main difficulty is the efficient treatment of the space dependent delay term. The method presented here is built upon the work in \cite{visser_efficient_2017}. We propose to use cubic Hermite splines for the interpolation of the history, finite differences for the discretisation of the diffusion term and for the actual time integration an implicit-explicit method is employed. The main benefit is that interpolation, numerical integration and diffusion become linear operations, hence leading to an efficient numerical scheme.   

\subsection{Overview}
In Section~\ref{sec:neural_field}, we construct the neural field model that consists of two populations of neurons, one excitatory and one inhibitory, with synaptic connections between all neurons and gap junctions between neurons of the same population. For the synaptic connections, we also take into account a spatial transmission delay due to the finite propagation speed of the action potentials and a sigmoidal activation function. This model is then reformulated into an abstract delay differential equation, allowing us to use the corresponding relevant theoretical tools. Finally, we derive an analytical condition for the eigenvalues of the linearisation around the trivial equilibrium, which determines its stability. The key idea which makes this possible is to switch to the basis of spherical harmonics, which leads to explicit formulas for the eigenfunctions. 

In Section~\ref{sec:Hopf}, we first review some aspects of the theory for Hopf bifurcations in the presence of spherical symmetry. We derive the normal forms for Hopf bifurcations on $V_l \oplus V_l$ with $l=1, \ 2, \,3$. Here $l$ corresponds to the order of the spherical harmonics, which span the eigenspace of the bifurcation. We summarise the conditions for the stability of certain periodic orbits, or branches, in terms of the normal form coefficients, as given in \cite{iooss_hopf_1989} and \cite{sigrist_hopf_2010}. Using the homological equation for the Hopf bifurcation, we derive analytic expressions for the normal form coefficients of these bifurcations. Here the decomposition into spherical harmonics again plays a key role.

In Section~\ref{sec:numerics}, the numerical method used to simulate the neural field model is presented. We first generate a mesh based on a geodesic projection on an icosahedron, which is then refined using triangular subdivisions. This allows us to discretize the nonlinear term over the elements. Diffusion is discretized using finite differences. For time integration, an implicit-explicit (IMEX) scheme is used, where the linear diffusion term is evaluated implicitly, and the nonlinear reaction term with delays is evaluated explicitly. 

In Section~\ref{sec:results}, some numerical examples are presented. We compute bifurcation diagrams for different values of diffusion parameters, investigate Hopf bifurcations and compute the normal form coefficients using formulas developed in Section~\ref{sec:Hopf}. Our theoretical predictions are verified by simulating the full neural field model using the method described in Section~\ref{sec:numerics}.

\section{The neural field model}
\label{sec:neural_field}
In this section, we construct the neural field model which is investigated in this paper. Next, we show how this model fits within the larger framework of abstract delay equations, so that existing theoretical tools to study the dynamical behaviour are applicable. In particular, we analyse the stability and bifurcations of the trivial equilibrium of the system. 

\subsection{Constructing the model}
Following \cite[Section 1.3]{spek_neural_2020}, consider two populations of excitatory and inhibitory neurons, where the voltages are given by $u_e$ and $u_i$, respectively. We assume that the neurons are placed on the unit sphere $\Omega:= S^2(\mathbb{R})$ and take a mean field approach, so $u_e(t,\br)$ ($u_i(t,\br)$) represents the averaged voltage of the excitatory (inhibitory) neurons at the place $\br\in \Omega$. Here $\br=\br(\theta,\phi)$ is a point on the sphere with coordinates $\theta\in [0,\pi]$ the polar angle and $\phi\in[0,2\pi)$ the azimuthal angle.  

When neurons receive no input, the voltage $u$ decays to the resting state with rate $\alpha$. Neurons are connected in two ways, via electrical connections, also called gap junctions, and synaptic connections. The electrical connections are very short range and only connect neurons of the same type, i.e., inhibitory to inhibitory (or excitatory to excitatory), hence, can be modelled by a diffusion term $d (\nabla \cdot \nabla)$, with diffusion coefficient $d$. Here $(\nabla \cdot \nabla)$ is the Laplace-Beltrami operator on the sphere.

For the synaptic connections, the amount of action potentials a cell sends out is scaled by an activation function $S$, which is modelled as a sigmoid. The incoming action potentials are weighted by the strength of the connectivity $J$. We assume that the connectivity decreases exponentially based on the great circle distance between neurons. Furthermore, the signal delay $\tau$ consists of a fixed delay $\tau_0$, representing dendritic integration and the processes at the synapse, and a delay depending on the great circle distance and speed of signal propagation $c$.

If we combine these modelling assumptions, we arrive at the following delay integro-differential equation
%\small
%\begin{equation}\label{eq:model_sphere}
%\begin{split}
%\dot u_e(t,\br)=& d_e (\nabla \cdot \nabla) u_e(t,\br) -\alpha_e u_e(t,\br) + \int_\Omega J_{ee}(\br,\br')S(u_e(t-\tau(\br,\br'),\br')) + J_{ei}(\br,\br')S(u_i(t-\tau(\br,\br'),\br'))d\br'\\    
%\dot u_i(t,\br)=& d_i (\nabla \cdot \nabla) u_i(t,\br) -\alpha_i u_i(t,\br) + \int_\Omega J_{ie}(\br,\br')S(u_e(t-\tau(\br,\br'),\br')) + J_{ii}(\br,\br')S(u_i(t-\tau(\br,\br'),\br'))d\br'
%\end{split}
%\end{equation}
%\normalsize
%for $\br\in \Omega$ and $t\geq 0$. 
\begin{equation}\label{eq:model_sphere}
    \dot{\mathbf{u}}(t)(\br) = \begin{pmatrix}
    d_e&0\\0&d_i
    \end{pmatrix}(\nabla \cdot \nabla)\mathbf{u}(t)(\br) -\begin{pmatrix}
    \alpha_e&0\\0&\alpha_i
    \end{pmatrix}\mathbf{u}(t)(\br)+\int_\Omega J(\br,\br')S(\mathbf{u}(t-\tau(\br,\br'),\br'))d\br'
\end{equation}
for $\br\in \Omega$ and $t\geq 0$, where we combined $u_e$ and $u_i$ in a vector and $J$ is the connectivity matrix defined as
\begin{equation}\label{Jmatrix}
    \mathbf{u}(t)(\br) := \left(\begin{matrix} u_e(t,\br)\\ u_i(t,\br) \end{matrix} \right),\quad
    J(\br,\br') = \begin{pmatrix}
    J_{ee}(\br,\br') & J_{ei}(\br,\br')\\
    J_{ie}(\br,\br') & J_{ii}(\br,\br')
    \end{pmatrix}. 
\end{equation}
Here $S$ acts element-wise on the components of $\mathbf{u}$. Moreover, $\alpha_e>0$, $\alpha_i>0$ and the diffusion coefficients $d_e>0, d_i>0$.

We make some natural choices for the connectivity kernel $J$, the delay $\tau$ and firing rate function $S$ 
\begin{equation}\label{eq:JtauS}
\begin{split}
J_{x}(\br,\br')&:= \eta_{x} \exp \left(-\frac{\arccos(\br\cdot \br')}{\sigma_x}\right)\\
\tau(\br,\br')&:= \tau_0 + \frac{\arccos(\br\cdot \br')}{c}\\
S(u) &:= \frac{1}{1+e^{-\gamma(u-\delta)}}-\frac{1}{1+e^{\gamma \delta}},
\end{split}
\end{equation}
with $\sigma_x>0$, $\eta_{ee}\eta_{ei}<0$, $\eta_{ie}\eta_{ii}<0$, $\tau_0>0$, $c>0$, $\delta\geq 0$, $\gamma>0$. Here $\arccos({\br\cdot \br'})$ gives the arclength between $\br$ and $\br'$. To simplify the further analysis, we have shifted the non-linearity $S$ such that $S(0)=0$. Hence $\mathbf{u}(t)$ should be understood as the voltage difference from the resting state.

\subsection{Abstract delay differential equations}
To study the dynamical behaviour of the neural field model, we reformulate \eqref{eq:model_sphere} in terms of an abstract delay differential equation. In this form, we can leverage existing theorems and methods for stability and bifurcations of equilibria in \cite{spek_neural_2020}.

Set $\mathbf{u}(t)\in Y:=C^{0,\alpha}(\Omega; \mathbb{R}^2)$, the Banach space of Hölder continuous functions with exponent $\alpha$, $\tfrac{1}{2}<\alpha\leq 1$. This space is common for functions defined on a spherical domain as the expansion into spherical harmonics converges uniformly on $Y$, see Theorem~\ref{thm:spherical_harmonics}. 

As \eqref{eq:model_sphere} is a system of delay equations, the proper state needs to also contain the history $\mathbf{u}_t(s):=\mathbf{u}(t+s)$ for $s \in [-h,0]$, where $h:=\max\{\tau(\br,\mathbf{r'}),\br,\mathbf{r'}\in\Omega\}>0$. The state space is then given by $X:=C([-h,0],Y)$.

Define the linear unbounded operator $B:D(B)\rightarrow Y$ and the compact non-linear operator $G:X\rightarrow Y$ as
\begin{equation}
\begin{split}
B\mathbf{u}(t)&:= \diag(d_e,d_i) (\nabla \cdot \nabla)\mathbf{u}(t) -\diag(\alpha_e,\alpha_i)  \mathbf{u}(t)\\
G(\mathbf{u}_t)&:= \int_\Omega J(\cdot,\br')S\left(\mathbf{u}_t(-\tau(\cdot,\br'))(\br')\right) d\br' ,
\end{split}
\end{equation}
where $D(B):=\{y\in Y| (\nabla \cdot \nabla)y\in Y \}$ and $S$ acts element-wise on the components of $\mathbf{u}_t$.

The Fréchet derivative $DG$ at $\mathbf{u}_t$ given in \cite{gils_local_2013} as
\begin{equation}\label{eq:F_derivative} 
DG(\mathbf{u}_t)\psi = \int_\Omega J(\cdot,\br') S' \left(\mathbf{u}_t(-\tau(\cdot,\br'))(\br')\right) \psi(-\tau(\cdot,\br'),\mathbf{r'}) d\br',
\end{equation}
which is similarly compact. Using these operators \eqref{eq:model_sphere} can be written as an abstract delay differential equation with an initial state (history) $\varphi\in X$ as
\begin{equation}\label{eq:abstractDDE}
\begin{cases}\dot{\mathbf{u}}(t)&=B\mathbf{u}(t) + G(\mathbf{u}_t)\\
\mathbf{u}_0&=\varphi.
\end{cases}
\end{equation}
In this form and with the above-mentioned properties of $B$ and $G$, we can apply all the relevant theorems of \cite{spek_neural_2020}.

\subsection{Stability of the trivial equilibrium}
As we have shifted the nonlinearity $S$ in \eqref{eq:JtauS} such that $S(0)=0$, the constant function $\mathbf{u}\equiv 0$ is an equilibrium of the system, called the trivial equilibrium. In this section, we study its stability. 

Define the linear operator $A:D(A)\rightarrow X$ as the linearisation of \eqref{eq:abstractDDE} at the trivial equilibrium, where $D(A):=\{\varphi \in C^1([-h,0];Y)|\,\varphi(0)\in D(B), \dot{\varphi}(0)=B\varphi(0)+DG\varphi\}$. To study the spectrum, we consider the complexified version of $X$. By \cite[Theorem 27]{spek_neural_2020}, the essential spectrum of $A$ is equal to the essential spectrum of $B$. As the Laplace-Beltrami operator is self-adjoint, it has an empty essential spectrum and all eigenvalues are real. Hence, $\sigma_{ess}=\emptyset$ when $d_e,d_i>0$.

Due to \cite[Proposition VI.6.7]{engel_one-parameter_1999}, we have that $\lambda \in \mathbb{C}$ is an eigenvalue of $A$ with corresponding eigenvector $\psi\in X$ if and only if $\psi(s)(\br)= e^{\lambda s}\mathbf{q}(\br)$, $s\in [-h,0]$, where $\mathbf{q}:=(q_e,q_i)^T \in Y$ is a non-trivial solution of
\begin{equation}\label{eq:char}
\Delta(\lambda) \mathbf{q} := \lambda \mathbf{q} - B\mathbf{q} - DG(\psi) = 0.
\end{equation}
Here $DG(\psi)$ is the first Fréchet derivative of $G$ at the trivial equilibrium
\begin{equation}
\begin{split}
DG(\psi)(\br) &:= \int_\Omega S'(0)J(\br,\br')\psi(-\tau(\br,\br'))(\br')d\br'\\
&= \int_\Omega S'(0)g(\br \cdot \br',\lambda){\bf q}(\br')d\br',
\end{split}
\end{equation}
where the matrix $g(\br\cdot \br',z):= \exp(-z\,\tau(\br,\mathbf{r'}))J(\br,\br')$ is well-defined, as the kernel $J$ and delay function $\tau$ in \eqref{eq:JtauS} only depend on $\br \cdot \br'$, that is
\[g(s,z) = \begin{pmatrix}
g_{ee}(s,z) & g_{ei}(s,z)\\
g_{ie}(s,z) & g_{ii}(s,z)
\end{pmatrix}, \quad g_x(s,z)=\exp(-z\cdot\tau(s))J_x(s),\]
with $s\in[-1,1]$ and $z\in \mathbb{C}$.

The trick in solving the characteristic equation \eqref{eq:char} is to use the expansion into spherical harmonics. The complex-valued spherical harmonics $Y_l^m$ of degree $l$ and order $|m|\leq l$, are the eigenfunctions of the Laplace-Beltrami operator $(\nabla \cdot \nabla)$ with a fixed spherical radius, i.e., 
\begin{equation}\label{eigenfunction_Laplacian}
    (\nabla \cdot \nabla) Y_l^m = -l(l+1) Y_l^m.
\end{equation}

The spherical harmonics form an orthonormal basis and converge uniformly in the space of Hölder continuous functions.
\begin{thm}[\cite{ragozin_polynomial_1970, ragozin_uniform_1971}]
\label{thm:spherical_harmonics}
Let $v \in C^{0,\alpha}(S^2,\mathbb{C})$, then the spherical harmonics $Y_l^m$ form an orthonormal basis in $C^{0,\alpha}(S^2,\mathbb{C})$ and 
\[\sum_{l=0}^\infty \sum_{m=-l}^l v_l^m Y_l^m = v \qquad \text{uniformly on }\Omega\]
where the coefficients $v_l^m$ are given by
\[v_l^m = \int_\Omega v(\br) \overline{Y_l^m(\br)}d\br.\]
\end{thm}

Then, by \cite[Theorem 2]{visser_standing_2017} it follows that
\begin{equation}\label{eq:diagonality}
\int_\Omega g_{ee}(\br \cdot \br',z)Y_l^m(\br')d\br' = G_{ee,l}(z) Y_l^m(\br),
\end{equation}
with coefficients $G_{ee,l}$ given as
\begin{equation}\label{eq:defintion_G}
G_{ee,l}(z) := \int_\Omega \int_\Omega g_{ee}(\br \cdot \br',z)Y_l^m(\br')d\br' \overline{Y_l^m(\br)} d\br = 2\pi \int_{-1}^1 g_{ee}(s,z) P_l(s) ds,
\end{equation}
where $P_l$ is the Legendre polynomial of degree $l$. Note that these coefficients are independent of $m$. Similarly, we can define $G_{ei,l},G_{ie,l}$ and $G_{ii,l}$, which can be combined into the matrix as 
\[G_l(z) = \begin{pmatrix}
G_{ee,l}(z) & G_{ei,l}(z)\\
G_{ie,l}(z) & G_{ii,l}(z)
\end{pmatrix}.\]

By Theorem~\ref{thm:spherical_harmonics}, any solution $\mathbf{q}$ of \eqref{eq:char} can be decomposed into spherical harmonics $Y_l^m$ and due to \eqref{eq:diagonality} it is sufficient to look at each $Y_l^m$ individually.

Let us take $\mathbf{q}(\br)=Y_l^m(\br) \left(\begin{smallmatrix} v_e \\ v_i \end{smallmatrix} \right)$ and substitute this into the characteristic equation \eqref{eq:char}. By \eqref{eq:diagonality} we get that $\Delta(\lambda)\mathbf{q} = E_l(\lambda) \left(\begin{smallmatrix} v_e \\ v_i \end{smallmatrix} \right)$, where the matrix $E_l(z) \in \mathbb{C}^{2\times 2}$ is given by
\[E_l(z):= \begin{pmatrix}
    z + \alpha_e + l(l+1)d_e - S'(0)G_{ee,l}(z) & -S'(0) G_{ei,l}(z) \\[7pt]
-S'(0) G_{ie,l}(z) & z + \alpha_i +l(l+1) d_i  - S'(0) G_{ii,l}(z)
\end{pmatrix}.\]
where $z\in \mathbb{C}$. Note that the diffusion parameters only effect $E_l$ when $l\geq 1$.

Therefore, $\lambda\in \mathbb{C}$ is an eigenvalue with $\mathbf{q}(\br)=Y_l^m(\br)\mathbf{v}$ when 
\begin{equation}\label{eq:E_l_determinant}
  \mathcal{E}_l(\lambda):= \det E_l(\lambda) = 0 \text{ for some } l\geq 0,
\end{equation}
where $\mathbf{v}:=\left(\begin{smallmatrix} v_e \\ v_i \end{smallmatrix} \right)$ are given by the null-space of $E_l(\lambda)$. In this case, the linearised system has $2l+1$ linearly independent solutions that are eigenfunctions of the form $\psi_m(\theta)(\br)=e^{\lambda \theta}Y_l^m(\br)\mathbf{v}$, $m=-l,\cdots,l$. 

Furthermore, we can solve the resolvent system $\Delta(z)\mathbf{q}=\mathbf{y}$, when $z\in \rho(A)$, the resolvent set, by expanding $\mathbf{q}$ and $\mathbf{y}$ into spherical harmonics with coefficients $\mathbf{q}_{\,l}^m$ and $\mathbf{y}_l^m$, respectively. Then they can be expressed as
\begin{equation}\label{eq:resolvent}
\mathbf{q}_{\,l}^m = E_l^{-1}(z) \mathbf{y}_l^m 
\end{equation}
and will play an important role in the computation of the normal form coefficients in Section~\ref{sec:computing_nfc}.

\section{Hopf bifurcations with spherical symmetry}\label{sec:Hopf with symmetry}
\label{sec:Hopf}
In this section, we investigate what happens if the trivial equilibrium loses its stability while it undergoes a Hopf bifurcation. The spherical symmetry of the neural field \eqref{eq:model_sphere} plays a key role in determining what patterns, standing and rotating waves will emerge. In this mathematical context, spherical symmetry means that the vector field in our neural field model \eqref{eq:model_sphere} commutes with the representation of the orthogonal group $O(3)$. Moreover, Hopf bifurcations add a symmetry with respect to phase shifts, given by the symmetry group $S^1$. First, we construct the proper normal forms and investigate the branching equations. Since our model is spherically symmetric, this implies that the normal form of the bifurcation should be symmetric too. By considering the relevant symmetries, we can construct the normal form, which can predict the dynamical behaviour close to the bifurcation. The normal form predicts the existence of certain limit cycles which branches from the bifurcation. These cycles are invariant with respect to different subgroups of $O(3)\times S^1 $. Finally, we obtain formulas to compute the normal form coefficients for a Hopf bifurcation. For a more in-depth overview of bifurcations with symmetries, we recommend \cite{golubitsky_singularities_1988}.

\subsection{General theory of Hopf bifurcations with symmetry}
Let $\Gamma$ be a symmetry group. A vector field $f: \mathbb{R}^{p}\to \mathbb{R}^{p}$ is called {\em equivariant} with respect to the symmetry group $\Gamma$ if 
\begin{equation}
f(\gamma \mathbf{z}) = \gamma f(\mathbf{z}) \text{ for all } \gamma\in\Gamma \text{ and } \mathbf{z}\in \mathbb{C}^{p}.
\label{eq:equivariance}
\end{equation}
A Hopf bifurcation occurs when the linearisation of $f$ has purely imaginary eigenvalues and the corresponding eigenspace is of the form $V \bigoplus V$, where $V$ is an absolutely irreducible representation of $\Gamma$. 

Near the Hopf bifurcation, we expect to find branches of periodic solutions. A spatiotemporal symmetry of a periodic solution $\mathbf{z}(t)$ is an element $(\gamma,\psi) \in \Gamma \times S^1$ such that 
\begin{equation}
    (\gamma, \psi)\cdot \mathbf{z}(t) := \gamma \mathbf{z}(t + \psi) = \mathbf{z}(t).
\end{equation}
Here $S^1$ is the group of phase shifts, where we assume that the period has been rescaled to $2\pi$. The isotropy subgroup $\Sigma_{\mathbf{z}(t)}$ of $\mathbf{z}(t)$ is defined as
\begin{equation}
    \Sigma_{\mathbf{z}(t)} := \{(\gamma, \psi) \in \Gamma \times S^1 \mid \gamma \mathbf{z}(t+\psi) = \mathbf{z}(t) \}.
\end{equation}
The fixed point subspace of an isotropy subgroup $\Sigma \subseteq \Gamma \times S^1$ is given by 
\begin{equation}
    \text{Fix}(\Sigma) = \{\mathbf{z}\in \mathbb{C}^p | \sigma \mathbf{z} = \mathbf{z}, \, \forall \sigma \in \Sigma\}.
\end{equation}
An isotropy subgroup $\Sigma$ is called maximal if there does not exist an isotropy subgroup $\Sigma'$ of $\Gamma \times S^1$ such that $\Sigma \subset \Sigma' \subset \Gamma \times S^1$, where the inclusions are strict. When $\text{dim Fix}(\Sigma)=2$, then $\Sigma$ is maximal and we say that $\Sigma$ is $\mathbb{C}$-axial.

The equivariant Hopf theorem guarantees the existence of periodic orbits with a maximal isotropy subgroup $\Sigma$, under some non-degeneracy conditions \cite{fiedler_global_2006}. We will call such periodic orbits maximal branches. Note that also periodic orbits with non-maximal isotropy subgroups can exist, but their existence is not guaranteed. %We will show in Section~\ref{sec:results} that this is indeed the case for our model.

\subsection{Constructing $O(3)\times S^1$ equivariant normal forms}
When analysing a Hopf bifurcation in the presence of a symmetry group $\Gamma$, we use the Birkhoff normal form \cite{golubitsky_singularities_1988}, which is equivariant with respect to the symmetry group $\Gamma \times S^1$, and so are all its truncations at finite order.

The natural irreducible representation of $O(3)$ of dimension $2l+1$ is given by $V_l$, the space of spherical harmonics $Y^m_l$ of degree $l$ and order $m=-l,\cdots, l$. The natural action of $O(3)$ on $V_l$ is generated by three symmetries: inversion of the sphere, an infinitesimal azimuthal rotation $\theta'$, an infinitesimal polar rotation $\varphi'$. The natural action of $S^1$ on $V_l$ is multiplication by a complex factor $e^{i \psi}$, where $\psi\in [0,2\pi)$. The matrices describing the actions of these symmetries on $V_l$ are 
\begin{equation}\label{eq:generating_matrices}
\begin{split}
M_{inv} &= (-1)^l I_{2l+1},\\
M_{\theta'} &= \left(w_{-l} | w_{-l+1} | \cdots | w_l\right),\\
M_{\varphi'} &= \diag (e^{-i l \varphi'},e^{i(-l+1) \varphi'},\cdots ,e^{i(l-1)\varphi'},e^{i l \varphi'}),\\
M_{\psi} &= e^{i\psi} I_{2l+1},
\end{split}    
\end{equation}
where $I_{2l+1}$ is the identity matrix and the columns of the matrix $M_{\theta'}$ are given by 
\[w_m = \left(0,\cdots,0,-\frac{1}{2}\sqrt{(l+m)(l-m+1)}\theta',1,\frac{1}{2}\sqrt{(l-m)(l+m+1)}\theta',0,\cdots,0\right)^T,\]
with $1$ in the $m$th row, for $m=-l,\cdots, l$, hence $M_{\theta'}$ is tridiagonal. For more details, we refer to \cite{sigrist_hopf_2010}.

A vector $x$ in the space $V_l\oplus V_l$ can be written as
\begin{equation}
    x(\theta, \phi) = \sum_{m=-l}^l z_m Y_l^m(\theta, \phi) + \overline{z}_m \overline{Y_l^m(\theta, \phi)}, \label{eq:eigenspace}
\end{equation}
so the action of $O(3)$ on $V_l\oplus V_l$ is determined by its action on 
\[\mathbf{z} = (z_{-l}, z_{-l+1},\cdots,z_l)^T \in \mathbb{C}^{2l+1}.\]
The action of $O(3)$ on $\mathbf{z}$ is given by the same matrices as the action of $O(3)$ on $V_l$, i.e., \eqref{eq:generating_matrices}.

Let $F_k:=(F_{k,-l},\cdots, F_{k,l})^T$ denote the Taylor expansion of a function $f: \mathbb{C}^{2l+1}\to \mathbb{C}^{2l+1}$ up to order $k$, where each component $F_{k,r}, r=-l,\cdots,l$ is a polynomial of $z_{-l},\cdots,z_{l},\overline{z}_{-l},\cdots,\overline{z}_l$ of order up to $k$. To find the Birkhoff normal form, we use that $F_k$ is equivariant with respect to $O(3)\times S^1$ for all orders $k$. To prove equivariance of the whole symmetry group, it is sufficient to prove equivariance with respect to the generating symmetries \eqref{eq:generating_matrices}.

The matrices $M_{inv}, M_{\varphi'}$ and $M_{\psi}$ are diagonal, so for the equivariance condition \eqref{eq:equivariance} it is sufficient to consider individual monomials of the form 
\begin{equation}\label{eq:monomials} 
z_{-l}^{i_{-l}}\cdots z_{l}^{i_{l}}\overline{z}_{-l}^{j_{-l}}\cdots \overline{z}_{l}^{j_{l}}, \quad i_{-l}+\cdots+i_l+j_{-l}+\cdots+j_l\leq k.
\end{equation}
First, observe that $F_k$ is equivariant with respect to $M_{\psi}$ if and only if in \eqref{eq:monomials}
\[\sum_{m=-l}^l i_m - j_m = 1,\]
that is, $F_{k,r}$ consists of monomials of odd order. This immediately implies equivariance with respect to $M_{inv}$, when using $\psi = \pi$. Moreover, $F_k$ is equivariant with respect to $M_{\varphi'}$ if and only if $F_{k,r}$ consists of monomials for which  
\[\sum_{m=-l}^l m(i_m-j_m) = r.\]

To find $F_k$ that are equivariant with respect to $M_{\theta'}$, we take $F_{k,r}$ as a linear combination of all the monomials of order up to $k$, which satisfy the two conditions above. We then use Mathematica to solve the systems of equations given by
\[\lim_{\theta'\to 0} \frac{F_k(M_{\theta'}\mathbf{z})-M_{\theta'}F_k(\mathbf{z})}{\theta'}=0\]
for all $\mathbf{z}=(z_{-l},\cdots,z_{l},\overline{z}_{-l},\cdots,\overline{z}_l)^T\in \mathbb{C}^{2l+1}$.

\subsection{The $O(3)\times S^1$ equivariant normal forms for $l=0,1,2,3$}
\label{sec:nf}
This section gives an overview of the Birkhoff normal forms with $O(3)\times S^1$ symmetry, where we compute the terms up to third order. The normal forms for $l=2$ and $l=3$ contain many terms, hence we give them in Appendix~\ref{sec:app_normal_form}.

Introduce the following notations
\begin{align*}
    \mathbf{z} :=& (z_{-l}, z_{-l+1},\cdots,z_l)^T\\
    \hat{\mathbf{z}} :=& \left((-1)^l \bar{z}_l, (-1)^{l-1}\bar{z}_{l-1},\cdots,(-1)^{-l}\bar{z}_{-l}\right)^T\\
    |\mathbf{z}|^2 :=& \sum_{m=-l}^l |z_m|^2,\\
    P^l(\mathbf{z}) :=& z_0^2 + 2\sum_{m=1}^l (-1)^m z_m z_{-m}.
\end{align*}
For $l=0$, we find the standard Poincaré normal form for the Andronov-Hopf bifurcation:
\begin{equation}
\dot{z} = \mu z + g_{0,1} z|z|^2 + \mathcal{O}(|z|^5).    
\end{equation}
The normal forms allow us to predict the stability of limit cycles in the neighbourhood of the bifurcation. The first Lyapunov coefficient $l_1 = \tfrac{1}{\omega} \re(g_{0,1})$ determines the stability of the resulting limit cycle, where the eigenvalues at the Hopf bifurcation are given by $\lambda = \pm I \omega$. From standard theory follows that there exists a unique limit cycle which is stable if $l_1$ is negative.

For $l=1$, we find the same normal form as in \cite{visser_standing_2017}
\begin{equation}\label{eq:nf1}
\dot{\mathbf{z}} = \mu \mathbf{z} + g_{1,1}\mathbf{z}|\mathbf{z}|^2 + g_{1,2}\hat{\mathbf{z}} P^1(\mathbf{z}) + \mathcal{O}(|\mathbf{z}|^5).    
\end{equation}

When $l=2$, the normal form is as in \cite{iooss_hopf_1989}
\begin{equation}\label{eq:nf2}
\dot{\mathbf{z}} = \mu \mathbf{z} + g_{2,1}\mathbf{z}|\mathbf{z}|^2 + g_{2,2}\hat{\mathbf{z}} P^2(\mathbf{z}) + g_{2,3} \mathbf{C}(\mathbf{z}) + \mathcal{O}(|\mathbf{z}|^5),   
\end{equation}
where $\mathbf{C(z)}$ is given in Appendix~\ref{sec:app_normal_form_l2}.

Finally, when $l=3$, the normal form is as given in \cite{sigrist_hopf_2010} 
\begin{equation}\label{eq:nf3}
\dot{\mathbf{z}} = \mu \mathbf{z} + g_{3,1}\mathbf{z}|\mathbf{z}|^2 + g_{3,2}\hat{\mathbf{z}} P^3(\mathbf{z}) + g_{3,3} \mathbf{Q}(\mathbf{z})+ g_{3,4} \mathbf{R}(\mathbf{z}) + \mathcal{O}(|\mathbf{z}|^5)    
\end{equation}
with $\mathbf{Q(z)}$ and $\mathbf{R(z)}$ given in Appendix~\ref{sec:app_normal_form_l3}.

\subsection{Amplitude equations for $l=1$}
\label{sec:amp_eq}
We can rewrite the normal forms in terms of the amplitudes and phases by using the substitution $z_m=r_m e^{i \psi_m}$ for $m=-l,\dots,l$. In this section, we will truncate the normal forms at third order. Due to the symmetry of the normal form, the amplitudes can be written in terms of a single phase difference $\psi = 2\psi_0- \psi_{-1}-\psi_{1}$ as
\begin{equation}
\begin{split}
    \dot{r}_{-1} &= \re(\mu) r_{-1} + \re(g_{1,1})r_{-1}\left(r_{-1}^2+r_0^2+r_1^2\right) + 2\re(g_{1,2})r_{-1}r_1^2 - \re(g_{1,2}) r_0^2 r_1 \cos(\psi)\\
    &+ \im(g_{1,2}) r_0^2 r_1 \sin(\psi)\\
    \dot{r}_{0} &= \re(\mu) r_{0} + \re(g_{1,1})r_{0}\left(r_{-1}^2+r_0^2+r_1^2\right) + \re(g_{1,2})r_{0}^3 - 2 \re(g_{1,2}) r_{-1} r_0 r_1 \cos(\psi)\\
    &- 2\im(g_{1,2}) r_{-1} r_0 r_1 \sin(\psi)\\
    \dot{r}_{1} &= \re(\mu) r_{1} + \re(g_{1,1})r_{1}\left(r_{-1}^2+r_0^2+r_1^2\right) + 2\re(g_{1,2})r_{-1}^2r_1 - \re(g_{1,2}) r_{-1} r_0^2  \cos(\psi)\\
    &+ \im(g_{1,2}) r_{-1} r_0^2 \sin(\psi).
\end{split}\label{eq:amp_eq}
\end{equation}
For the phase equations, the same observation holds and lead to 
\begin{equation}
\begin{split}
    \dot{\psi}_{-1} &= \im(\mu) + \im(g_{1,1})\left(r_{-1}^2+r_0^2+r_1^2\right) + 2 \im(g_{1,2}) r_1^2  - r_{-1}^{-1}r_0^2 r_1(\im(g_{1,2}) \cos(\psi) \\
    &+ \re(g_{1,2}) \sin(\psi))\\
    \dot{\psi}_{0} &= \im(\mu) + \im(g_{1,1})\left(r_{-1}^2+r_0^2+r_1^2\right) + \im(g_{1,2})r_{0}^2 - 2 \im(g_{1,2}) r_{-1} r_1 \cos(\psi) \\
    &+ 2\re(g_{1,2}) r_{-1} r_1 \sin(\psi)\\
    \dot{\psi}_{1} &= \im(\mu) + \im(g_{1,1})\left(r_{-1}^2+r_0^2+r_1^2\right) + 2\im(g_{1,2})r_{-1}^2 - r_1^{-1} r_{-1} r_0^2 (\im(g_{1,2}) \cos(\psi) \\
    &+ \re(g_{1,2}) \sin(\psi)).
\end{split}\label{eq:phase_eq}
\end{equation}
The differential equation for the phase difference $\psi$ is given by
\begin{equation}
\begin{split}
\dot{\psi} &= -2 \im(g_{1,2}) \left(r_{-1}^2-r_0^2+r_{-1}^2\right) + \im(g_{1,2})\left(-4 r_{-1}r_1 + r_0^2\left(\frac{r_1}{r_{-1}}+\frac{r_{-1}}{r_1}\right)\right)\cos(\psi)\\
&+ \re(g_{1,2})\left(4 r_{-1}r_1 + r_0^2\left(\frac{r_1}{r_{-1}}+\frac{r_{-1}}{r_1}\right)\right)\sin(\psi).
\end{split}\label{eq:phase_diff}
\end{equation}

If we set the phase difference $\psi\equiv 0$, then the amplitude equations \eqref{eq:amp_eq} decouple from the phase equations 
\begin{equation}\label{eq:amplitude_eq_l1}
    \dot{r}_m = \re(\mu) r_m + \re(g_{1,1})r_m \left(r_{-1}^2+r_0^2+r_1^2\right) + \re(g_{1,2}) \hat{r}_m\left(r_0^2-2r_1r_{-1}\right) 
\end{equation}
for $m=-1,0,1$, where $\hat{r}=(-r_1,r_0,-r_{-1})^T$. By introducing the new coordinates $\rho_1 := r_{-1}^2+r_0^2+r_1^2$ and $\rho_2:= r_0^2-2r_{-1}r_1$, system \eqref{eq:amplitude_eq_l1} reduces to 
\begin{align*}
\dot{\rho_1} =& 2\rho_1(\re(\mu)+ \re(g_{1,1})\rho_1) + 2\re(g_{1,2}) \rho_2^2\\
\dot{\rho_2} =& 2\rho_2\left(\re(\mu)+(\re(g_{1,1})+\re(g_{1,2}))\rho_1 \right).
\end{align*}
Furthermore, \eqref{eq:phase_diff} for $\psi=0$ is given by
\[\dot{\psi} = \im(g_{1,2})\rho_2 \frac{(r_{-1}+r_1)^2}{r_{-1}r_1}.\]

%Make an assumption that r_i(t) e^{i \mu t + \phi(0)}
We are interested in fixed points $(\rho_1,\rho_2,\psi)$ with $\psi\equiv 0$ when $\re(\mu)>0$. This would correspond to (quasi)-periodic solutions of the normal form. This implies that $\rho_2=0$, which in turn demands that either $\rho_1=0$ or $\rho_1=\tfrac{\re(\mu)}{-\re(g_{1,2})}$. The case of $\rho_1=0$ implies that $r_{-1}=r_0=r_1=0$, which is the trivial solution, where $\dot{\psi}$ is undefined. However, $\rho_1=\tfrac{\re(\mu)}{-\re(g_{1,2})}>0$ when $\re(g_{1,2})<0$ which gives a valid solution to the amplitude equations. If we substitute $(\rho_1,\rho_2,\psi)=\left(\tfrac{\re(\mu)}{-\re(g_{1,2})},0,0\right)$ into the phase equations \eqref{eq:phase_eq}, we find that
\[\dot{\psi}_{-1}=\dot{\psi}_{0}=\dot{\psi}_{1}=\im(\mu)-\frac{\im(g_{1,1})}{\re(g_{1,1})}\re(\mu)\]
Hence, the corresponding solution of the normal form is periodic, as the derivatives of the phases are constant and equal. Projecting this family of limit cycles to the eigenspace $V_1 \bigoplus V_1$, we find rotating waves, where the remaining $3$ free dimensions correspond to the axis of rotation and the initial phase of the rotating wave. 

\begin{figure}
    \centering
    \includegraphics[scale=0.8]{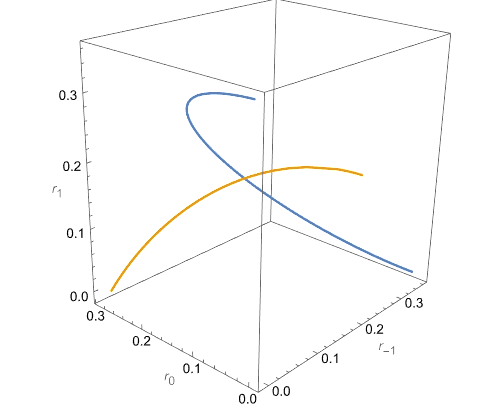}
    \caption{Curves corresponding to rotating waves \eqref{eq:rotating_wave} (blue) and standing waves \eqref{eq:standing_wave} (orange) for the values $\re(\mu)=0.1, \re(g_{1,1})=-0.810, \re(g_{1,2})=-0.405$.}
    \label{fig:ameq_l1}
\end{figure}

To investigate the stability of such limit cycles, we compute the eigenvalues of the linearisation of the coupled systems \eqref{eq:amp_eq} and \eqref{eq:phase_diff}. We can parametrise this family of limit cycles by $r_1$ as
\begin{equation}\label{eq:rotating_wave}
(r_{-1},r_0,r_1,\psi)= \left(\sqrt{\frac{\re(\mu)}{-\re(g_{1,2})}} - r_1, \sqrt{2r_1\left(\sqrt{\tfrac{\re(\mu)}{-\re(g_{1,2})}} - r_1\right)}, r_1,0 \right)    
\end{equation}
see also Figure~\ref{fig:ameq_l1}. We find that the eigenvalues are independent of the choice of $0\leq r_1 \leq \sqrt{\tfrac{\re(\mu)}{-\re(g_{1,2})}}$ and are given as $0,-2\re(\mu),-2\re(\mu) \tfrac{\re(g_{1,2})}{\re(g_{1,1})} \pm 2\re(\mu)\tfrac{\im(g_{1,2})}{\re(g_{1,1})}i$. The zero eigenvalue is due to the fact that we have a curve of fixed points, but as the other eigenvalues have negative real part when $\re(\mu),\tfrac{\re(g_{1,2})}{\re(g_{1,1})}>0$, we conclude that these rotating waves are stable for these normal form coefficients. 

If we set $\psi\equiv\pi$, then by similar reasoning we find a curve of fixed points
\begin{equation}\label{eq:standing_wave}
(r_{-1},r_0,r_1,\psi) = \left(r_1,\sqrt{\frac{\re(\mu)}{-\re(g_{1,1})-\re(g_{1,2})}-2r_1^2},r_1,\pi \right)    
\end{equation}
for $0\leq r_1\leq \sqrt{\tfrac{\re(\mu)}{-2(\re(g_{1,1})+\re(g_{1,2}))}}$. Substituting these into the amplitude \eqref{eq:amp_eq} and phase equations \eqref{eq:phase_eq}, we get that
\[\dot{r}_{-1}=\dot{r_0}=\dot{r_1}=0\qquad \dot{\psi}_{-1}=\dot{\psi}_{0}=\dot{\psi}_{1}=\im(\mu)-\frac{\im(g_{1,1})+\im(g_{1,2})}{\re(g_{1,1})+\re(g_{1,2})}\re(\mu)\]
Projecting this family of limit cycles into the eigenspace $V_1 \bigoplus V_1$, we find standing waves, where the remaining $3$ free dimensions correspond to the orientation and the initial phase of the standing wave. 

To investigate the stability of this limit cycle, we compute the eigenvalues of the linearisation of the coupled systems \eqref{eq:amp_eq} and \eqref{eq:phase_diff}, which are $0, -2\re(\mu)$ and a double eigenvalue $2\re(\mu) \tfrac{\re(g_{1,2})}{\re(g_{1,1})+\re(g_{1,2})}$. Hence,the standing waves are stable when $\re(\mu)>0$ and $\tfrac{\re(g_{1,2})}{\re(g_{1,1})+\re(g_{1,2})}<0$. 

\subsection{Branching equations for $l=2$ and $l=3$}\label{sec:branch}
The complete analysis of the amplitude equations corresponding to the normal forms for $l=2$ and $l=3$ proves to be intractable. The equivariant Hopf theorem guarantees however that at least the maximal branches exist. To determine their stability, conditions on the normal forms have been derived by Iooss \cite{iooss_hopf_1989} for the $l=2$ and by Sigrist \cite{sigrist_hopf_2010} for the $l=3$. In this section, we summarise their main results.

For $l=2$, there exist five branches: one axisymmetric solution, which is a limit cycle in the $Y_2^0$ direction and is symmetric along rotations in the azimuthal direction, one standing wave with $D_4$ symmetry, the symmetry group of the square, one standing wave with tetrahedral symmetry and two rotating waves with $D_4$ symmetry, the first rotating in the plane of the square and the second rotating around the central axis of the square. Iooss and Rossi gives conditions when these branches are stable \cite{iooss_hopf_1989}. Note that here, the authors use a different notation for the third order normal form coefficients, which can be identified with our notations as: $a=g_{2,1}, b=-\tfrac{1}{2}g_{2,2}, c=g_{2,3}$. For the axisymmetric branch and the standing wave with $D_4$ symmetry, one condition is dependent on the fifth order normal form coefficients, which we do not need to consider for the purposes of this paper. Furthermore, the second rotating wave is generically unstable, if $\re(g_{2,3})\neq 0$. 

\begin{figure}%[H]
    \centering
    \includegraphics[scale=0.8]{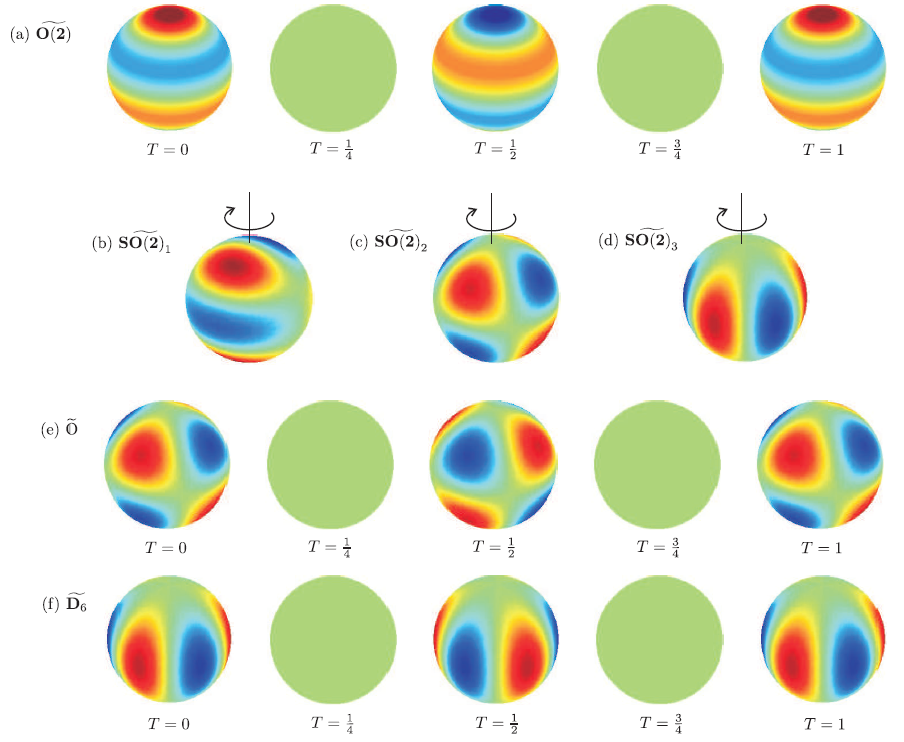}
    \caption{The six periodic solution branches with $\mathbb{C}$-axial symmetry. (a), (e) and (f) illustrate the evolution of the three standing waves over one period and (b), (c) and (d) illustrate the rotating wave solutions showing the axis and direction of rotation \cite{sigrist_hopf_2010}.}
    \label{fig:branches}
\end{figure}

For $l=3$, the equivariant Hopf theorem guarantees the existence of six branches of periodic solutions, see Figure~\ref{fig:branches}. Each of the branches is characterised by a subgroup of $O(3)\times S^1$, which represents the isotropy subgroup $\Sigma$ of the periodic orbit. These are: $\widetilde{O(2)}$, the symmetry group of the circle, $\widetilde{SO(2)}$, the group of rotations, which is equivalent to $S^1$, $\widetilde{\mathcal{O}}$ and the symmetry group of the octahedron and $\widetilde{D_6}$, the symmetry group of the hexagon. The three branches corresponding to $\widetilde{SO(2)}$ are all rotating waves, while the other three are standing waves. The equivariant Hopf theorem also gives conditions when each of these branches are stable, see Table 8 in \cite{sigrist_hopf_2010}. Note that the authors use a different notation for the normal form coefficients, which can be again identified with our as: $A=g_{3,1}, B=g_{3,2}, C=g_{3,3}, D=g_{3,4}$.

In addition to the maximal branches, \cite{iooss_hopf_1989} also considered quasi-periodic solutions and \cite{sigrist_hopf_2010} sub-maximal branches and stationary spiral patterns. Their existence cannot, however, be guaranteed by the equivariant Hopf theorem.

\subsection{Homological equation}
In order to obtain formulas for the normal form coefficients in terms of the parameters of our original problem, we use the so-called homological equation. It relates on one side the dynamics of our model, written in terms of the coordinates on the center manifold, $\mathbf{z}$ to the dynamics of the normal form on the other side. 

The construction of the homological equation uses the variation of constants formula. For delay equations, this requires a detour to the space $X^{\odot *}$. For a comprehensive overview of the sun-star calculus, see \cite{diekmann_delay_1995} and \cite{janssens_class_2019,janssens_class_2020,spek_neural_2020}. 

In short, the state space for delay equations is the space of continuous functions $X=C([-h,0];Y)$, as we need to include the history in the state space. In this case, the generator $A$ of the solution operator of a linear delay equation is given by translation in time
\[A \varphi = \dot{\varphi}, \qquad \varphi \in D(A)\]
and the extension of the state is part of the domain of $A$
\[D(A):= \{\varphi \in X| \dot{\varphi}(0) = B\varphi(0) + DG\varphi \}.\]
Hence, to consider (non-linear) perturbations, we need to perturb the domain of $A$. This means that a perturbation effects $\varphi(0)$, which is not a bounded perturbation in the space of continuous functions. 

This problem is circumvented by moving to a larger space. The sun-star calculus gives a natural procedure to construct this larger space $X^{\odot *}$. This space contains the space $Y\times X$, so the extension and translation can be decoupled and both included in the action of $A^{\odot *}$. The space $X$ can be embedded continuously into $X^{\odot *}$ with embedding $j$, where $j\varphi=(\varphi(0),\varphi)$. Similarly, $Y$ also has a natural continuous embedding $\ell$, $\ell y = (y,0)$. It can be shown that a perturbation of $A^{\odot *}$ with range $\ell Y$ is bounded in $X^{\odot *}$ and solutions of the perturbed equations map back into $X$. Hence, we can define a proper variation of constants formula.

Let $\mathcal{H}:\mathbb{C}^{2l+1}\times \mathbb{C}^{2l+1}\rightarrow X$ be the coordinates to the center manifold, 
\begin{equation}
\mathcal{H}(\mathbf{z},\overline{\mathbf{z}})= \sum_{m=-l}^l z_m \psi_m + \bar{z}_m\bar{\psi}_m + \sum_{2\leq |\mathbf{j}|+|\mathbf{k}| \leq 3} \frac{1}{\mathbf{j}! \mathbf{k}!} h_{\mathbf{j}\,\mathbf{k}}\mathbf{z}^\mathbf{j}\overline{\mathbf{z}}^\mathbf{k} + O(\|(\mathbf{z},\overline{\mathbf{z}})\|^4),
\end{equation}
where $\psi_m, \bar\psi_m$ are critical eigenfunctions, $h_{\mathbf{j}\,\mathbf{k}}\in X$, $\mathbf{z}=(z_{-l},\cdots,z_{l})$ and where $\mathbf{j},\mathbf{k}$ are multi-indices of length $2l+1$ and
\[|\mathbf{j}|:= \sum_{m=-l}^l |j_m|, \qquad \mathbf{j}! := \prod_{m=-l}^l j_m!, \qquad \mathbf{z}^\mathbf{j}:= \prod_{m=-l}^l z_m^{j_m}. \]
Furthermore, let $R:=G-DG$ represent the nonlinear terms and $F:\mathbb{C}^{2l+1} \rightarrow \mathbb{C}^{2l+1}$ the relevant normal form. 

\cite{spek_neural_2020} give the following form for the homological equation
\begin{equation}\label{eq:homological_equation}
A^{\odot *} j \mathcal{H}(\mathbf{z},\overline{\mathbf{z}}) + \ell R(\mathcal{H}(\mathbf{z},\overline{\mathbf{z}})) = j \mathcal{H}_{\mathbf{z}}(\mathbf{z},\overline{\mathbf{z}})F(\mathbf{z}) + j \mathcal{H}_{\overline{\mathbf{z}}}(\mathbf{z},\overline{\mathbf{z}})\overline{F(\mathbf{z})}
\end{equation}
for all $\mathbf{z} \in \mathbb{C}^{2l+1}$, where $\mathcal{H}_{\mathbf{z}}$ denotes the partial derivative of $\mathcal{H}$ with respect to $\mathbf{z}$. Writing both sides as a Taylor series, we can equate the coefficients of different powers of $z_m$ and $\bar{z}_m$ to construct equations for the center manifold coefficients $h_{\mathbf{j}\,\mathbf{k}}$ and normal form coefficients. Due to symmetries used in the construction of the normal forms, there are multiple equations which contain the normal form coefficients. So by choosing the right powers, the problem can be simplified drastically. 

\subsection{Computing normal form coefficients}
\label{sec:computing_nfc}
Let us start with a simple $l=0$ Hopf bifurcation. By this we mean that $A$ has a pair of eigenvalues  $\lambda= i \omega, \overline{\lambda} = -i\omega$ for which  $\mathcal{E}_0(i \omega)=0$ and $\mathcal{E}_l(i \hat{\omega})\neq 0$ for all $l\geq 1$ or $\omega \neq \hat{\omega}\in \mathbb{R}$, with corresponding eigenfunctions $\psi(\theta)(\br) := e^{i\omega \theta}Y^0_0(\br) v$ and its complex conjugate $\overline{\psi}(\theta)(\br) = e^{-i\omega \theta}Y^0_0(\br) \overline{v}$, where $v=(v_e,v_i)$, normalised such that $\|v\|_2:=\overline{v}^Tv=1$.

As we discussed in Section~\ref{sec:nf}, the normal form for this bifurcation is given by
\begin{equation}
    \dot z = \mu z + g_{0,1} z |z|^2 +\mathcal{O}(|z|^5), \ z\in\mathbb{C},
\end{equation}
with the normal form coefficient $g_{0,1}\in\mathbb{C}$. We aim to compute the first Lyapunov coefficient $l_1 = \tfrac{1}{\omega} \re(g_{0,1})$, as its sign determines the stability of the resulting limit cycle.

This bifurcation is identical to the generic Hopf bifurcation for which in \cite{spek_neural_2020} the formulas for computing the center manifold coefficients $h_{2\,0},h_{1\,1}$ and normal form coefficient $g_{0,1}$ were derived as
\begin{equation}\label{eq:normal_form_computation} 
\begin{split}
h_{2\,0}(\theta)=& e^{2i\omega\theta}\Delta^{-1}(2 i \omega) D^2G(\psi,\psi)\\
h_{1\,1}(\theta)=& \Delta^{-1}(0) D^2G(\psi,\bar{\psi})\\ 
2g_{0,1} \psi(\theta) =& \frac{1}{2\pi i} \oint_{\partial C_\lambda} e^{z\theta}\Delta^{-1}(z)(D^3G(\psi,\psi,\bar{\psi}) + D^2G(h_{2\,0},\bar{\psi}) + 2 D^2G(\psi,h_{1\,1}))\,dz.
\end{split}
\end{equation}

Here, the $n$th Fréchet derivative $D^nG$ around the trivial equilibrium, given in \cite{gils_local_2013} is
\begin{equation}\label{eq:nth_derivative} 
D^nG(\psi_1,\cdots,\psi_n)= S^{(n)}(0) \int_\Omega J(\cdot,\br') \left(\prod_{j=1}^n \psi_j(-\tau(\cdot,\br'),\mathbf{r'})\right) d\br',
\end{equation}
where the product between brackets is an element-wise product of vectors, resulting from the fact that $G$ contains no cross-terms of $u_e$ and $u_i$. Denote this component-wise product by $*$, also known as the of Schur or Hadamard product.

For the computation of $h_{2\,0}$, and $h_{1\,1}$ we use the resolvent formula \eqref{eq:resolvent}. To do this, we first need to compute the coefficients of the expansion of $D^2 G(\psi,\psi)$ and other functions involved in \eqref{eq:normal_form_computation} into spherical harmonics. In general, the definition of $G_l$ in  \eqref{eq:defintion_G} and \cite[Theorem 2]{visser_standing_2017} imply that for any $f\in Y$
\begin{equation}\label{eq:expansion1}
 \int_\Omega \int_\Omega g(\br\cdot \br',z)f(\br') d\br' \; \overline{Y_l^m(\br)} d\br = G_l(z) \int_\Omega f(\br) \overline{Y_l^m(\br)}d\br.   
\end{equation}
Note that in this case, both $g$ and $G_l$ are $2\times 2$ matrices. 

Apply first \eqref{eq:nth_derivative} and \eqref{eq:expansion1} to obtain 
\begin{align*}
D^2G(\psi,\psi)(\br)=&S''(0) \int_\Omega g(\br\cdot \br',2 i \omega)(v*v) Y_0^0(\br') Y_0^0(\br') d\br' ,\\
\int_\Omega D^2G(\psi,\psi)(\br) \overline{Y_l^m(\br)}d\br =& S''(0) G_l(2 i \omega) (v*v)\int_\Omega Y_0^0(\br) Y_0^0(\br) \overline{Y_l^m(\br)}d\br\\
=& \frac{1}{2\sqrt{\pi}} S''(0) G_l(2 i \omega) (v*v)\int_\Omega Y_0^0(\br) \overline{Y_l^m(\br)}d\br\\[5pt]
=& \begin{cases}\frac{1}{2\sqrt{\pi}} S''(0) G_l(2 i \omega) (v*v) & l=m=0\\
0 & \text{otherwise.}
\end{cases}
\end{align*}

We can now use the resolvent formula to find $h_{2\,0}$ and by a similar procedure $h_{1\,1}$
\begin{equation}
\begin{split}
h_{2\,0}(\theta,\br) =&\frac{1}{2\sqrt{\pi}} e^{2 i\omega \theta} S''(0)  Y_0^0(\br) Q_0(2 i\omega)(v*v),\\
h_{1\,1}(\theta,\br) =&\frac{1}{2\sqrt{\pi}}S''(0) Y_0^0(\br)Q_0(0) (v*\overline{v}),
\end{split}
\end{equation}
where, for ease of notation we introduce the matrix $Q_l(z):= E_l^{-1}(z)G_l(z)$, $l\geq 0$. Repeating these steps for the last equation of \eqref{eq:normal_form_computation}
\begin{equation*}
\begin{split}
e^{z\theta}&\Delta^{-1}(z)\left(D^3G(\psi,\psi,\bar{\psi}) + D^2G(h_{2\,0},\bar{\psi}) + 2 D^2G(\psi,h_{1\,1})\right) \\
=& \frac{1}{4 \pi}e^{z\theta}Y^0_0(\br) E_0^{-1}(z) G_0(i\omega)\Bigl[S'''(0)(v*v*\overline{v})+ (S''(0))^2 \Bigl((Q_0(2 i\omega)(v*v))*\overline{v}\\
& \hspace{7.2cm} +2 (Q_0(0) (v*\overline{v}))*v\Bigr)\Bigr].
\end{split}
\end{equation*}
By the definition of $E_0$ and Cauchy's integral formula, we find that 
\[\frac{1}{2\pi i} \oint_{\partial C_\lambda} E_0^{-1} (z) \, dz = \frac{1}{2\pi i} \oint_{\partial C_\lambda} \frac{\adj(E_0(z))}{\mathcal{E}_0(z) } \, dz = \frac{\adj(E_0( i\omega))}{\mathcal{E}'_0( i\omega)},\]
where $\adj$ denotes the adjugate matrix. Note that, all other terms in the integral in \eqref{eq:normal_form_computation} are analytic at $z=i \omega$, as $\lambda= i\omega$ is a simple zero of $\mathcal{E}_0$. Hence, the integral in \eqref{eq:normal_form_computation} can be evaluated as 
\begin{equation}
    \begin{split}
    2 g_{0,1} e^{ i\omega \theta}Y^0_0(\br) v
    =& \frac{1}{4 \pi}e^{ i\omega \theta}Y^0_0(\br) \frac{\adj(E_0( i\omega))}{\mathcal{E}'_0( i\omega)}G_0( i\omega) (S'''(0)(v*v*\overline{v}) \\
    &+(S''(0))^2 ((Q_0(2 i\omega)(v*v))*\overline{v}+2 (Q_0(0) (v*\overline{v}))*v)),
    \end{split}
\end{equation}
from which we obtain
\begin{equation}
    \begin{split}
    g_{0,1} = \frac{1}{8 \pi} \overline{v}^T \frac{\adj(E_0( i\omega))}{\mathcal{E}'_0( i\omega)}G_0( i\omega) (S'''(0)(v*v*\overline{v}) &+(S''(0))^2 (Q_0(2 i\omega)(v*v))*\overline{v}\\
    &+2(S''(0))^2 (Q_0(0) (v*\overline{v}))*v).
    \end{split}
\end{equation}

With some effort, along the same line, the normal form coefficients for $l=1,2,3$ can also be calculated. A detailed overview of this calculation can be found in Appendix~\ref{sec:computing_nfc_app}. One interesting note is that for $l=2$, $g_{2,3}$ is only dependent on $S''(0)$ and not on $S'''(0)$.

\section{Numerical discretisation}\label{sec:numerics}
As a starting point, we extend the numerical method developed for delayed neural field equations in \cite{visser_standing_2017} to systems. The details of this method were not included in the above paper, hence we give an overview here, based on unpublished personal communication. Next, we include the discretisation of the diffusion term and then describe some time integration methods that are suitable and efficient for both the delay and diffusion part. The numerical method described in this paper is implemented in Python.

\begin{figure}%[H]
	\centering
		\includegraphics[scale=.4]{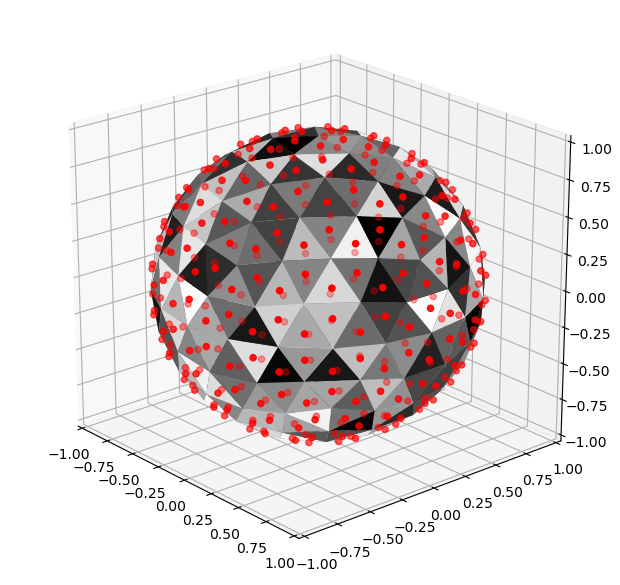}
		\caption{Mesh and centroids (red dots) after three refinement steps. }
		\label{fig:mesh}
\end{figure}

We construct first a regular triangulation of the surface of the sphere. The mesh generation procedure, as described also in \cite{baumgardner_icosahedral_1985}, starts with projecting the twelve vertices of the regular icosahedron onto the sphere. The mesh consists of 20 spherical triangles that are bounded by 30 geodesic arcs. A mesh refinement means that the edges are bisected, generating from each triangle four and then the new vertices are again projected onto the sphere, see Figure~\ref{fig:mesh} for a typical mesh with three refinement steps. 
The computational space domain $\Omega$ consists of a tessellation of $m$ non-overlapping spherical triangles $\Omega_j$ as
\begin{equation}
    \Bigl\{ \Omega_j : \bigcup_{j=1}^m \bar\Omega_j = \Omega,\ \Omega_i\cap\Omega_j=\emptyset,\ i\not=j \bigr\}.
\end{equation}
The mesh points of our discretisation are the centers of mass of all triangles, called \emph{centroids}, and are denoted by ${\mathbf r}_j$, $j=1,\cdots,m$.  

\subsection{Discretisation of the nonlinear term}\label{sec:discretize_nonlinearity}
Consider the neural field system in the form \eqref{eq:abstractDDE} and discretize the integral in the nonlinear part $G$ using the centroids as quadrature points as 
\begin{equation}\label{eq:nonlinearity_approx}
    \begin{split}
        G(u_t)(\br_j)=&\int_\Omega J(\br_j,\br')S(u(t-\tau(\br_j,\br'),\br')) d\br'=\sum_{\nu=1}^m\int_{\Omega_\nu} J(\br_j,\br')S(u(t-\tau(\br_j,\br'),\br')) d\br'\\    
        \approx & \sum_{\nu=1}^m J(\br_j,\br_\nu)S(u(t-\tau(\br_j,\br_\nu),\br_\nu)) |{\Omega_\nu}|,\quad \br_\nu\in\Omega_\nu,
        \end{split}
\end{equation}
for any $\br_j\in\Omega_j$, where $u$ is one of the components of $\mathbf u$ and $J$ one of the four connectivities in \eqref{eq:model_sphere}.

On any fixed mesh there are $m^2$ delays in the system and from \eqref{eq:nonlinearity_approx} it is clear that for every mesh element $\Omega_j$ the history needs to be interpolated for $m$ time delays, i.e., for $\tau(\br_j,\br_\nu)$, $\nu=1,\cdots,m$. In this section we focus on the efficient evaluation of the nonlinear term containing the delays. This is based on the work proposed in \cite{visser_efficient_2017}, which we adapted to our system. The main idea of this method is as follows. To determine the state of the system at given time lags we use cubic Hermite splines to interpolate over the history. The chosen interpolation method requires that the time step is constant. It turns out that both interpolation and numerical integration are linear operations and hence they lead to an efficient numerical scheme.  

Introduce the following notations for the delays and the state of the system (segment of the solution) at time $t$ 
\begin{equation}
    \begin{split}
        & \tau_{j,\nu}:=\tau(\br_j,\br_\nu), \quad  u_{t,\nu}(s):=u(t+s,\br_\nu),\, s\in[-h,0].
    \end{split}
\end{equation}
For a fixed time step $\Delta t$, let $k$ be the smallest integer multiple of the time step that covers the whole history interval $[-h,0]$, i.e.,
\begin{equation}
    k:=\min_{\kappa\in\mathbb{N}} \{\kappa : \kappa\Delta t\geq h\}.
\end{equation}
Set $t_i=i\Delta t$ and $\br_\nu\in\Omega_\nu$, $\nu=1,\cdots m$, and represent the approximations of the solution and its derivative at a space-time mesh point as   
\begin{equation}
    u^h_{i,\nu}\approx u(t_i,\br_\nu)\quad \mbox{and} \quad \dot u^h_{i,\nu} \approx \dot u(t_i,\br_\nu).
\end{equation}

For the approximation of a segment of the solution and its derivative we have
\begin{equation}\label{eq:approximate_u_uderivative}
\begin{split}
    u^h_{i-\ell,\nu} &\approx u(t_{i-\ell},\br_\nu) = u(t_i-\ell\Delta t,\br_\nu) = u_{t_i}(-\ell\Delta t,\br_\nu) = u_{t_i,\nu}(-\ell\Delta t),\quad \ell=0,\cdots, k\\
    \dot{u}^h_{i-\ell,\nu} &\approx \dot{u}(t_{i-\ell},\br_\nu) = \dot{u}(t_i-\ell\Delta t,\br_\nu) = \dot{u}_{t_i}(-\ell\Delta t,\br_\nu)= \dot{u}_{t_i,\nu}(-\ell\Delta t),\quad \ell=1,\cdots,k.
\end{split}
\end{equation}
Note that, in addition to the approximation of the solution at each time step, also all derivatives in \eqref{eq:approximate_u_uderivative} up to, but not including, the current time step are known, since they have been evaluated during the forward time stepping. This allows us to use cubic Hermite splines, with one exception that will be clarified. Hence, to represent the approximation of the entire history $u_{t_i,\nu}$ we use the vectors
\begin{equation}
     %{\bf u}^h_{i,\nu} := \left(u^h_{i,\nu},\cdots,u^h_{i-k,\nu}\right)^T, \quad 
     %\dot{\bf u}^h_{i,\nu} := \left(\dot{u}^h_{i-1,\nu},\cdots,\dot{u}^h_{i-k,\nu}\right)^T,
     {\bf u}^h_{i,\nu} := \begin{pmatrix}
        u^h_{i,\nu}\\ \vdots \\ u^h_{i-k,\nu}
     \end{pmatrix},\quad
     \dot{\bf u}^h_{i,\nu} :=\begin{pmatrix}
        \dot{u}^h_{i-1,\nu}\\ \vdots\\ \dot{u}^h_{i-k,\nu}
     \end{pmatrix},
\end{equation}
that is, ${\bf u}^h_{i,\nu}(\ell) = u^h_{i-\ell,\nu}$ for $\ell=0,\cdots,k$ and $\dot{\bf u}^h_{i,\nu}(\ell) = \dot{u}^h_{i-\ell,\nu}$ with $\ell=1,\cdots,k$. 

The main question is: for any fixed $\br_j\in \Omega_j$, \emph{how to approximate $u_{t_i,\nu}(-\tau_{j,\nu})$ in \eqref{eq:nonlinearity_approx} if we know ${\bf u}^h_{i,\nu}$ and $\dot{\bf u}^h_{i,\nu}$?} 

Let $\ell$ be such that $\tau_{j,\nu}\in [\ell\Delta t, (\ell+1)\Delta t)$. We distinguish two cases:  
\begin{itemize}
    \item[(a)] When $\ell\geq 1$, then ${\bf u}^h_{i,\nu}(\ell)$, ${\bf u}^h_{i,\nu}(\ell+1)$ and the derivatives $\dot{\bf u}^h_{i,\nu}(\ell)$, ${\dot{\bf u}^h_{i,\nu}(\ell+1)}$ are known, which motivates to use cubic Hermite splines for the interpolation of $u_{t_i,\nu}(-\tau_{j,\nu})$, that is 
    \begin{equation}
    \begin{split}\label{eq:interpolation_l}
        u_{t_i,\nu}(-\tau_{j,\nu}) &\approx p_3(\tau_{j,\nu}-\ell\Delta t) {\bf u}^h_{i,\nu}(\ell)+ p_3(\tau_{j,\nu}-(\ell+1)\Delta t) {\bf u}^h_{i,\nu}(\ell+1)\\[5pt]
        & \quad - q_3(\tau_{j,\nu}-\ell\Delta t)\dot{\bf u}^h_{i,\nu}(\ell)-q_3(\tau_{j,\nu}-(\ell+1)\Delta t)\dot{\bf u}^h_{i,\nu}(\ell+1)\\
        &=\sum_{\ell=1}^k p_3(\tau_{j,\nu}-\ell\Delta t) {\bf u}^h_{i,\nu}(\ell)- \sum_{\ell=1}^k q_3(\tau_{j,\nu}-\ell\Delta t)\dot{\bf u}^h_{i,\nu}(\ell),
    \end{split}
    \end{equation}
    where the third order spline basis functions are
    \begin{equation}
        \begin{split}
            & p_3(s) = \begin{cases}
            2\left|\frac{s}{\Delta t}\right|^3 - 3\left|\frac{s}{\Delta t}\right|^2+1, & s\in [-\Delta t,\Delta t) \\
            0 & \text{otherwise}
            \end{cases} \\[7pt]
            & q_3(s) = \begin{cases}
            s\left(1-\left|\frac{s}{\Delta t}\right|\right)^2, & s\in [-\Delta t,\Delta t) \\
            0 & \text{otherwise.}
            \end{cases}
        \end{split}
    \end{equation}
    Note that the last equality in \eqref{eq:interpolation_l} holds as only two terms are nonzero in both sums.
    \item[(b)] In the special case when $\ell=0$, i.e., $\tau_{j,\nu}\in [0, \Delta t)$, only the three values ${\bf u}^h_{i,\nu}(0), {\bf u}^h_{i,\nu}(1)$ and $\dot{\bf u}^h_{i,\nu}(1)$ are known, hence we can only use a quadratic polynomial approximation 
    \begin{equation}
        u_{t_i,\nu}(-\tau_{j,\nu})\approx p_{30}(\tau_{j,\nu}) {\bf u}^h_{i,\nu}(0) + p_{31}(\tau_{j,\nu}-\Delta t) {\bf u}^h_{i,\nu}(1) - q_{31}(\tau_{j,\nu}-\Delta t)\dot{\bf u}^h_{i,\nu}(1),
    \end{equation}
    with quadratic basis functions
    \begin{equation}
    \begin{split}
        & p_{30}(s) = \begin{cases}
        \left(1-\frac{s}{\Delta t}\right)^2, & s\in [0,\Delta t) \\
        0 & \text{otherwise}
        \end{cases} \\[7pt]
        & p_{31}(s) = \begin{cases}
        1-\left|\frac{s}{\Delta t}\right|^2, & s\in [-\Delta t,0) \\
         p_3(s) & s\in[0,\Delta t) \\
        0 & \text{otherwise}
        \end{cases} \\[7pt]
        & q_{31}(s) = \begin{cases}
        s\left(1+\left|\frac{s}{\Delta t}\right|\right), & s\in [-\Delta t,0) \\
         q_3(s) & s\in[0,\Delta t) \\
        0 & \text{otherwise}.
        \end{cases}
    \end{split}
    \end{equation}
\end{itemize}

Combining (a) and (b) leads to the interpolation on the entire interval $[0,k\Delta t]$
\begin{equation}\label{eq:interpolation_delay}
    \begin{split}
        u_{t_i,\nu}(-\tau_{j,\nu}) & \approx  p_{30}(\tau_{j,\nu}) {\bf u}^h_{i,\nu}(0) + p_{31}(\tau_{j,\nu}-\Delta t) {\bf u}^h_{i,\nu}(1)+\sum_{\ell=2}^k p_3(\tau_{j,\nu}-\ell\Delta t) {\bf u}^h_{i,\nu}(\ell) \\
        & \quad -q_{31}(\tau_{j,\nu}-\Delta t)\dot{\bf u}^h_{i,\nu}(1) - \sum_{\ell=2}^k q_3(\tau_{j,\nu}-\ell\Delta t)\dot{\bf u}^h_{i,\nu}(\ell).
    \end{split}
\end{equation}
  
We would like to write this into a matrix form, hence define the spline basis vector functions ${\bf p}_3(\cdot)$ and ${\bf q}_3(\cdot)$ of the delays, such that each returns a row vector, i.e., 
\begin{equation}
    {\bf p}_3 (\tau_{j,\nu}) := \left(
    p_{30}(\tau_{j,\nu}),  p_{31}(\tau_{j,\nu}-\Delta t),  p_3(\tau_{j,\nu}-2\Delta t),  \cdots,p_3(\tau_{j,\nu}-k\Delta t) 
    \right)
\end{equation}
and 
\begin{equation}
    {\bf q}_3 (\tau_{j,\nu}) := \left(
     q_{31}(\tau_{j,\nu}-\Delta t),  q_3(\tau_{j,\nu}-2\Delta t), \cdots,  q_3(\tau_{j,\nu}-k\Delta t) \right)
\end{equation}
for $\nu=1,\cdots,m$. Then, \eqref{eq:interpolation_delay} reads 
\begin{equation}
    u_{t_i,\nu}(-\tau_{j,\nu}) \approx {\bf p}_3(\tau_{j,\nu}) {\bf u}^h_{i,\nu} - {\bf q}_3(\tau_{j,\nu}) \dot {\bf u}^h_{i,\nu},\quad \nu=1,\cdots,m.
\end{equation}
The main idea is now to use the above interpolation strategy to interpolate $S(u_{t_i,\nu}(-\tau_{j,\nu}))$, given $S({\bf u}^h_{i,\nu})$ as follows
\begin{equation}
    S\left(u_{t_i,\nu}(-\tau_{j,\nu})\right) \approx {\bf p}_3(\tau_{j,\nu}) {\bf v}_{i,\nu} - {\bf q}_3(\tau_{j,\nu}) \dot{\bf v}_{i,\nu},\quad \nu=1,\cdots,m,
\end{equation}
where
\begin{align*}
    {\bf v}_{i,\nu} &= S({\bf u}^h_{i,\nu})\text{ ($S$ acting element-wise)}\\[5pt]  
    \dot{\bf v}_{i,\nu} &= S'({\bf u}^h_{i,\nu}) * \dot{\bf u}^h_{i,\nu} \text{ (by the chain-rule and element-wise multiplication).}
\end{align*}

Let us turn our attention to the interpolation of the nonlinearity \eqref{eq:nonlinearity_approx}, i.e., 
\begin{equation}
    \begin{split}
        G(u_{t_i})(\br_j)\approx \sum_{\nu=1}^m J(\br_j,\br_\nu)S(u_{t_i,\nu}(-\tau_{j,\nu}))\, |{\Omega_\nu}|
        = \sum_{\nu=1}^m J_{j,\nu}\left({\bf p}_3(\tau_{j,\nu}) {\bf v}_{i,\nu} - {\bf q}_3(\tau_{j,\nu}) \dot{\bf v}_{i,\nu}\right)\, |{\Omega_\nu}|,
    \end{split}
\end{equation}
where $J_{j,\nu}=J(\br_j,\br_\nu)$. Next, we extend this to all centroids $\br_j$, $j=1,\cdots, m$ by constructing the matrices
\begin{equation}
    P:=\begin{pmatrix}
     J_{1,1}\, {\bf p}_3(\tau_{1,1})|\Omega_1| & J_{1,2} \,{\bf p}_3(\tau_{1,2})|{\Omega_2}| & \cdots & J_{1,m} \,{\bf p}_3(\tau_{1,m})|\Omega_m|\\
     J_{2,1} \,{\bf p}_3(\tau_{2,1})|\Omega_1| & J_{2,2} \,{\bf p}_3(\tau_{2,2})|{\Omega_2}| & \cdots & J_{2,m} \,{\bf p}_3(\tau_{2,m})|\Omega_m|\\
     \vdots & & & \\
     J_{m,1} \,{\bf p}_3(\tau_{m,1})|\Omega_1| & J_{m,2} \,{\bf p}_3(\tau_{m,2})|{\Omega_2}| & \cdots & J_{m,m} \,{\bf p}_3(\tau_{m,m})|\Omega_m|
    \end{pmatrix}\in\mathbb{R}^{m\times m(k+1)}
\end{equation}
and similarly 
\begin{equation}
    Q:=\begin{pmatrix}
     J_{1,1} \,{\bf q}_3(\tau_{1,1})|\Omega_1| & J_{1,2} \,{\bf q}_3(\tau_{1,2})|{\Omega_2}| & \cdots & J_{1,m} \,{\bf q}_3(\tau_{1,m})|\Omega_m|\\
     J_{2,1} \,{\bf q}_3(\tau_{2,1})|\Omega_1| & J_{2,2} \,{\bf q}_3(\tau_{2,2})|{\Omega_2}| & \cdots & J_{2,m} \,{\bf q}_3(\tau_{2,m})|\Omega_m|\\
     \vdots & & & \\
     J_{m,1} \,{\bf q}_3(\tau_{m,1})|\Omega_1| & J_{m,2} \,{\bf q}_3(\tau_{m,2})|{\Omega_2}| & \cdots & J_{m,m} \,{\bf q}_3(\tau_{m,m})|\Omega_m|
    \end{pmatrix}\in\mathbb{R}^{m\times m k}.
\end{equation}

With these notations, the discretized neural field equation (without diffusion) has the form 
\begin{equation}\label{eq:semidiscreteNF}
    \begin{pmatrix} \dot u_{i,1}^h \\ \vdots \\ \dot u_{i,m}^h \end{pmatrix} = 
    -\alpha \begin{pmatrix} u^h_{i,1} \\ \vdots \\ u^h_{i,m} \end{pmatrix} + 
    P \begin{pmatrix}  {\bf v}_{i,1} \\ \vdots \\  {\bf v}_{i,m} \end{pmatrix} -
    Q \begin{pmatrix} \dot{\bf v}_{i,1} \\ \vdots \\ \dot{\bf v}_{i,m} \end{pmatrix},
\end{equation}
which shows that both integration and interpolation can be expressed as a linear operation, which can be evaluated efficiently. Moreover, $P$ and $Q$ do not depend on the time step $i$ for which $t_i=i\Delta t$, hence they need to be computed only once at the beginning of the time simulation. 

This formulation holds for both $u_e$ and $u_i$ components of $\bf u$. Although, we need to compute four spline matrices corresponding to the four connectivities in \eqref{Jmatrix}, we only need to compute them once during the whole simulation. The forward time integration scheme to solve the resulting system will be discussed later when we have discretized also the diffusion term.

\subsection{Discretisation of the diffusion term}\label{sec:discretize_diffusion}
The numerical discretisation of the diffusion term is based on \cite{oostendorp_interpolation_1989, huiskamp_difference_1991}, where finite difference formulas were derived for the approximation of the surface Laplacian on a triangulated surface. 

Let $u$ be again one of the two components of $\bf u$. Consider a mesh point $\br_j$ surrounded by $N$ direct neighbours $\br_i$ at distance $h_{ij}$, $i=1,\cdots,N$, respectively. In \cite{oostendorp_interpolation_1989}, the mesh points are the vertices of the triangulation and direct neighbours are considered the neighbouring vertices. Since our computational points are the centroids, we set as direct neighbouring points of a mesh point the centroids of the adjacent triangles, hence in our case $N=3$. The general form of the approximation is
\begin{equation}\label{eq:approx_diffusion}
    (\nabla \cdot \nabla) u (\br_j) \approx \Delta_h u(\br_j) := \sum_{i=1}^N w_{ij}\left(u(\br_i)-u(\br_j)\right),
\end{equation}
with weights
\begin{equation}
    w_{ij} = \frac{4}{\overline{h}_j}\frac{1}{N}\frac{1}{h_{ij}}
\end{equation}
and $\bar{h}_j$ the mean distance of the neighbours to $\br_j$. The discrete Laplacian $\Delta_h$ can be expressed as a matrix $D$ consisting the weights, with elements
\[
    D_{ij} = \begin{cases}
            w_{ij} & i\not=j,\ \br_i\text{ is a direct neighbour of }\br_j\\
            -\frac{4}{\overline{h}_j}\left(\overline{\frac{1}{h_{j}}}\right) & i=j\\
            0 & i\not=j,\ \br_i\text{ is not a direct neighbour of }\br_j,
    \end{cases}
\]
where $\left(\overline{\frac{1}{h_{j}}}\right)$ is the average of $1/h_{ij}$ over the neighbours of $\br_j$. Note that, distance is always measured along the surface of the sphere. When $D$ is multiplied with the vector $u=\left(u(\br_1),\cdots,u(\br_m)\right)^T$, it leads to the approximation of the Laplacian in all mesh points, i.e.,
\begin{equation}
    (\nabla \cdot \nabla) u\approx \Delta_h u = D u.
\end{equation}
This can easily be incorporated in our framework described in the previous section. 

As the spherical harmonics of degree $l$ solve the eigenvalue problem \eqref{eigenfunction_Laplacian}, a straightforward analytical test for our numerical solution is to solve the differential equation
\begin{equation}\label{eq:test_diffusion}
    (\nabla \cdot \nabla) u = -Y_l^m, 
\end{equation}
for any $l$ and $|m|\leq l$. 

Remark that, although theoretical analysis of convergence and accuracy of the numerical method is not the scope of this paper, we verified the numerical solution using two error measures, also used in \cite{baumgardner_icosahedral_1985}. The first,
\begin{equation}
    \frac{\langle Y_l,Y_l \rangle}{\langle U,Y_l\rangle}\approx l(l+1)
\end{equation}
is a quotient of two inner products that approximates the eigenvalue of the spherical harmonic of degree $l$. Here $U$ is the numerical solution of \eqref{eq:test_diffusion} and $Y_l$ is one of the spherical harmonics $Y_l^m$ with order $|m|\leq l$. Consequently, if the computations were exact, then this quotient would equal the eigenvalue. The second measure, i.e., $\|U-u\|$ is the error between the exact and numerical solution. The results of the error computations are summarised in Table~\ref{table:errors} for three harmonic degrees and four mesh refinements.   

\begin{table}
\hspace{-.5cm}
\small
%\begin{center}
\begin{tabular}{c|c|c|c}
 &  $l=1$ & $l=2$ & $l=3$ \\ 
    \hline
\begin{tabular}{c}
\\
 $m$\\ \\
 \hline
320\\
1280\\
5120\\
20480
\end{tabular}
&
\begin{tabular}{ c|c } 
\\
$\langle Y_l,Y_l \rangle/\langle U,Y_l\rangle$ & $\|U-u\|$  \\ \\
\hline
 1.9856 & 0.0514 \\ 
 1.9873 & 0.0382\\ 
 1.9871 & 0.0155 \\ 
 1.9869 & 0.0207 \\ 
\end{tabular}
&
\begin{tabular}{ c|c } 
\\
 $\langle Y_l,Y_l \rangle /\langle U,Y_l\rangle$ & $\|U-u\|$ \\ \\
\hline
 5.9106 & 0.3048 \\ 
 5.9495 & 0.0968 \\ 
 5.9556 & 0.1628 \\ 
 5.9560 & 0.0923 \\ 
%\hline
\end{tabular}
&
\begin{tabular}{ c|c } 
\\
$\langle Y_l,Y_l \rangle/\langle U,Y_l\rangle$ & $\|U-u\|$ \\ \\
\hline
 10.9843 & 0.0641 \\ 
 11.2032 & 0.1075 \\ 
 11.2489 & 0.0848 \\ 
 11.2560 & 0.0779 \\ 
%\hline
\end{tabular}
\end{tabular}
%\end{center}
\normalsize
\caption{Error values of the Laplacian approximation for three harmonic degrees and four mesh refinements. $m$ is the number of mesh points (centroids) corresponding to the $n$th refinement, with $n=2,3,4,5$.}
\label{table:errors}
\end{table}

\subsection{Time integration}
In \cite{visser_standing_2017} an explicit Euler method was used for the time integration of the neural field equation in \eqref{eq:semidiscreteNF}. A fully explicit scheme places, however, severe restriction on the time step when the diffusion term is included. A fully implicit method requires the implicit treatment of the nonlinear term in each time step, which may be expensive. Moreover, relaxation of the time step restriction may not be worth for the extra work spent to solve the resulting implicit system at each time step. Furthermore, most popular implicit schemes are first order. All these motivate us to use an implicit-explicit (IMEX) scheme for time integration developed for the accurate approximation of time-dependent PDEs, see for example \cite{ascher_implicit-explicit_1995}. 

Consider the neural field system which has been discretized in space in Section~\ref{sec:discretize_nonlinearity} and Section~\ref{sec:discretize_diffusion}, respectively. We present the method again for one component of $\bf u$ but it can be extended to systems in a straightforward way. The system of ODEs has the form
\begin{equation}\label{eq:space_discrete}
        \dot{{u}}(t) = d\, \Delta_h {u} +F({u}),
\end{equation}
where $d$ is the diffusion coefficient and $F({u})$ contains the term $(-\alpha u)$ and the nonlinear term that we do not want to integrate implicitly. Thus, we integrate it explicitly and the diffusion term implicitly, yielding an IMEX scheme.

A comparison of several IMEX schemes, suitable for reaction-diffusion problems for pattern formation, is presented in \cite{ruuth_implicit-explicit_1995}. Our investigation shows that a well suited method for the numerical integration of \eqref{eq:space_discrete} is the modified version of Crank-Nicolson, Adams-Bashforth method (MCNAB) that has the form
\begin{align}
    \frac{u^{n+1}-u^n}{\Delta t}= \frac{3}{2} F(u^n)- \frac{1}{2} F(u^{n-1})+d\left( \frac{9}{16}\Delta_h u^{n+1} + \frac{3}{8}\Delta_h u^{n} + \frac{1}{16}\Delta_h u^{n-1}\right),
\end{align}
where $u^n$ denotes the solution at time step $n$. Note that, this method requires the storage of $F(u^{n-1})$. Its evaluation in the first time step is done using a one-step implicit Euler method. We can then express $u^{n+1}$ as
\begin{equation*}
    u^{n+1} = M_{MCNAB}^{-1}\left[ u^n+\Delta t \left( \frac{3}{2} F(u^n)- \frac{1}{2} F(u^{n-1})\right) + 
    d\, \Delta t\,  D \left( \frac{3}{8} u^{n} + \frac{1}{16} u^{n-1}\right)\right],
\end{equation*}
where $M_{MCNAB}= I- \tfrac{9}{16}d\, \Delta t\,  D$ and $I$ the identity matrix. 

Remark that, it is possible to use different time steps for the forward time integration, $\Delta t$, and spline interpolation of the history, $\Delta t_{\text{spline} }$. We only need to ensure that $\Delta t_{\text{spline }}$ is an integer multiple of $\Delta t$. 

\section{Exploring bifurcations using numerical examples}\label{sec:results}
In this section, we investigate numerically how the diffusion influences the formation of patterns in the neural field model. Using synthetic parameters we find the Hopf bifurcations described in Section~\ref{sec:Hopf with symmetry}. We construct bifurcation diagrams, compute eigenvalues and normal form coefficients, and make simulations of the full nonlinear system \eqref{eq:model_sphere} for different values of $d_i,d_e$ to analyse their effects on the dynamics. The numerical computations in all examples use the mesh created with four refinement steps, which results in 5120 triangles (centroids). 

\subsection{Constructing bifurcation diagrams}
An equilibrium has two different co-dimension one bifurcations, those of fold-type with a zero eigenvalue, $\lambda=0$ and those of Hopf-type with a pair of purely imaginary eigenvalues $\lambda = \pm i \omega$. Constructing a curve of fold-type bifurcations just requires solving the implicit equation $\mathcal{E}_l(0)=0$ for the chosen bifurcation parameters. Finding a curve of Hopf bifurcations is however, more difficult as one also needs to find $\omega$. 

We make some simplifications to our model \eqref{eq:model_sphere}, by assuming that the connectivity parameters only depend on the pre-synaptic neuron, i.e., $\eta_{e}:=\eta_{ie}=\eta_{ee}$, $\eta_i:=\eta_{ei}=\eta_{ii}$, $\sigma_e:=\sigma_{ie}=\sigma_{ee}$, $\sigma_i:=\sigma_{ei}=\sigma_{ii}$. This assumption, often implicitly made like in \cite{visser_standing_2017}, implies that a neuron acts the same regardless of the postsynaptic neuron. 

The bifurcation parameters are $\eta_{e}$ and $\eta_{i}$, and for ease of notation define $\hat{G}_{e,l}$ and $\hat{G}_{i,l}$ such that
\begin{equation*}
\eta_e \hat{G}_{e,l}(z) = S'(0) G_{ee,l}(z),\quad
\eta_i \hat{G}_{i,l}(z) = S'(0) G_{ii,l}(z).
\end{equation*}
To find a curve of Hopf bifurcations, we need to solve the equation $\mathcal{E}_l(i \omega)=0$ for some non-zero $\omega \in \mathbb{R}$. Using the notations above, the determinant in \eqref{eq:E_l_determinant} can be written as
\begin{equation}
\begin{split}
\mathcal{E}_l(i \omega) &= \left(i \omega +\alpha_e+l(l+1)d_e- \eta_e \hat{G}_{e,l}(i \omega)\right) \left(i\omega +\alpha_i+l(l+1)d_i-\eta_i \hat{G}_{i,l}(i \omega)\right) \\
&\quad-\eta_e \eta_i \hat{G}_{e,l}(i \omega) \hat{G}_{i,l}(i \omega)\\
&= (i \omega +\alpha_e+l(l+1)d_e)(i\omega +\alpha_i+l(l+1)d_i) -  \eta_e (i\omega +\alpha_i+l(l+1)d_i) \hat{G}_{e,l}(i \omega)\\
&\quad -\eta_i (i \omega +\alpha_e+l(l+1)d_e) \hat{G}_{i,l}(i \omega). 
\end{split}
\end{equation}
Setting the real and imaginary parts to zero, gives the following real linear system
\begin{equation}
\begin{split}
&\begin{pmatrix}
\re( (i\omega +\alpha_i+l(l+1)d_i) \hat{G}_{e,l}(i \omega)) & \re( (i \omega +\alpha_e+l(l+1)d_e) \hat{G}_{i,l}(i \omega) )\\
\im( (i\omega +\alpha_i+l(l+1)d_i) \hat{G}_{e,l}(i \omega)) & \im( (i \omega +\alpha_e+l(l+1)d_e) \hat{G}_{i,l}(i \omega) )\\
\end{pmatrix}
\begin{pmatrix}
\eta_e \\ \eta_i
\end{pmatrix}\\
&= \begin{pmatrix}
\re((i \omega +\alpha_e+l(l+1)d_e)(i\omega +\alpha_i+l(l+1)d_i))\\
\im((i \omega +\alpha_e+l(l+1)d_e)(i\omega +\alpha_i+l(l+1)d_i))
\end{pmatrix}.
\end{split}
\end{equation}
Given $\omega$, the solution of this linear system for $\eta_e,\eta_i$ produces a Hopf bifurcation. In the bifurcation diagrams that follow we set $\eta_{e}>0,\eta_{i}<0$ and fix the parameters $\sigma_{e}=2/9$, $\sigma_{i}=1/6$, $\tau_0=3$, $c=0.8$, $\alpha_e=\alpha_i=1$, $\gamma=8$, $\delta = 0$, like in \cite{visser_standing_2017}. 

\subsection{No diffusion and equal diffusion}
\begin{figure}%[H]
    \centering
    \includegraphics[scale=0.8,trim={0.8cm 0 0 0},clip]{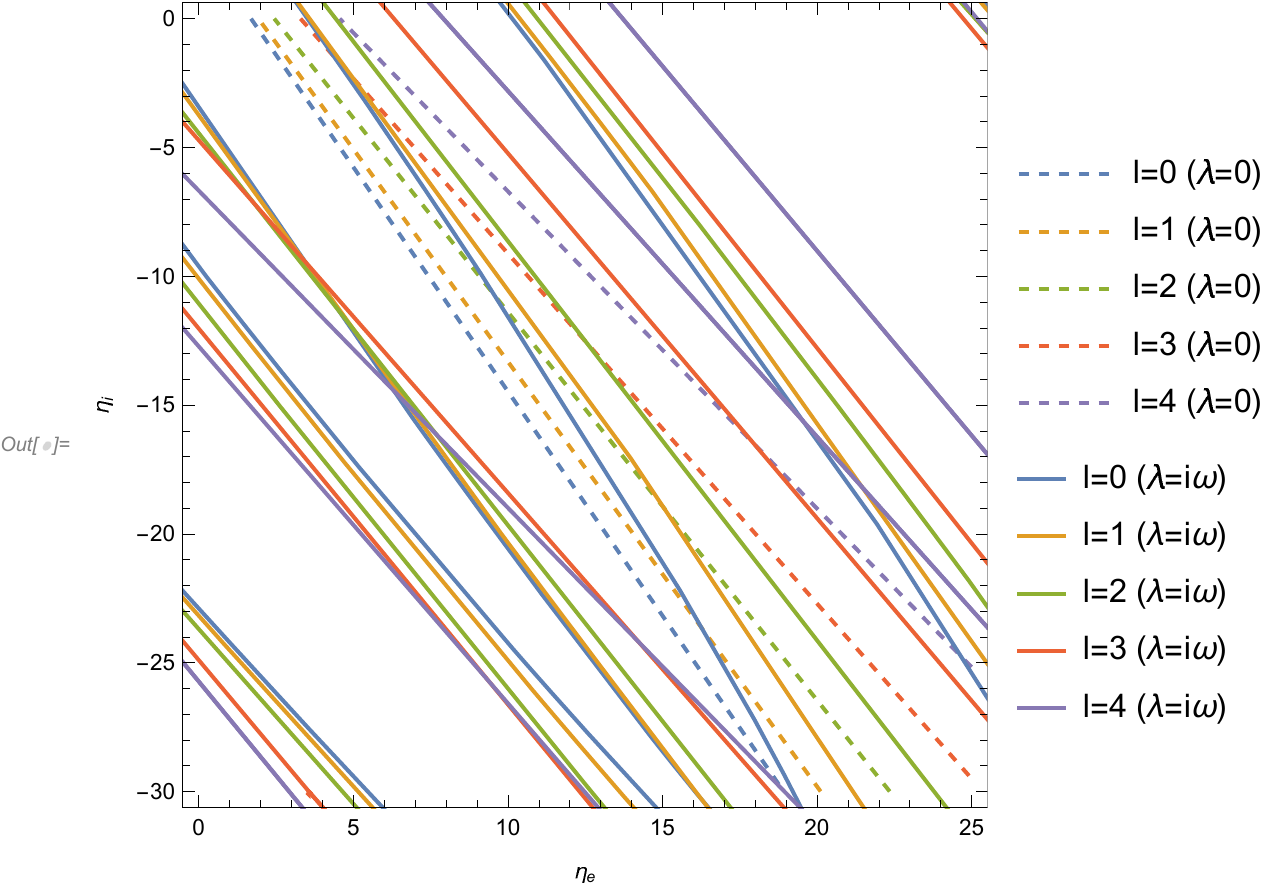}
    \caption{Bifurcation diagram in $\eta_{e}$ and $\eta_{i}$, up to $l=4$, with $d_e=d_i=0$. Bifurcations of Hopf-type are given by an unbroken curve, and those of fold-type are given by a dashed line. In the white region containing $(0,0)$ the trivial equilibrium is stable.}
    \label{fig:bif_diag_no_diffusion}
\end{figure}

First, we investigate the case when there is no diffusion, $d_e=d_i=0$. See Figure~\ref{fig:bif_diag_no_diffusion} for a bifurcation diagram with fold and Hopf-type bifurcation curves up to $l=4$.

When $\eta_e=\eta_i=0$, the trivial equilibrium is globally asymptotically stable as \eqref{eq:model_sphere} reduces to
\begin{align*}
\dot{u}_e(t) &= -\alpha_e u_e(t)\\
\dot{u}_i(t) &= -\alpha_i u_i(t).
\end{align*}
In the region containing $(\eta_e,\eta_i)=(0,0)$ and is bounded by fold and Hopf-type bifurcation curves, the trivial equilibrium is stable. Taking a point in this region and decreasing $\eta_i$ leads to the trivial equilibrium losing stability due to a Hopf-type bifurcation of many different orders $l$ close together. 

\begin{figure}
    \centering
    \includegraphics[scale=0.8,trim={0.8cm 0 0 0},clip]{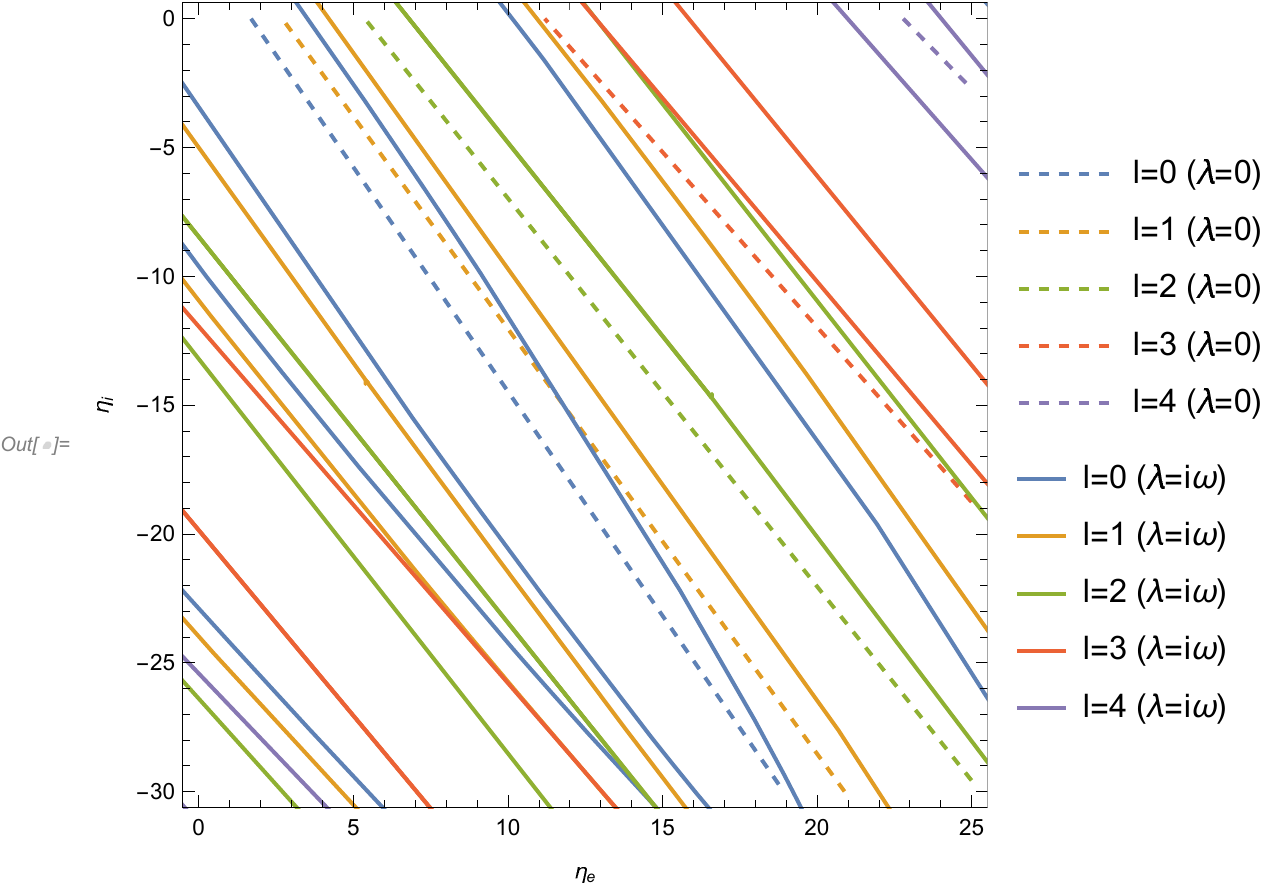}
    \caption{Bifurcation diagram in $\eta_{e}$ and $\eta_{i}$, up to $l=4$, with $d_e=d_i=0.2$. Bifurcations of Hopf-type are given by an unbroken curve, and those of fold-type are given by a dashed line. In the white region containing $(0,0)$ the trivial equilibrium is stable.}
    \label{fig:bif_diag_equal_diffusion}
\end{figure}

Set now the diffusion parameters non-zero but equal, $d_e=d_i=0.2$. Compared to the bifurcation diagram without diffusion, bifurcations of higher degree $l$ have moved away from the region of stability, Figure~\ref{fig:bif_diag_equal_diffusion}.

\subsection{Hopf bifurcation $l=0$}
\begin{figure}%[H]
    \centering
    \includegraphics[scale=0.8,trim={0.8cm 0 0 0},clip]{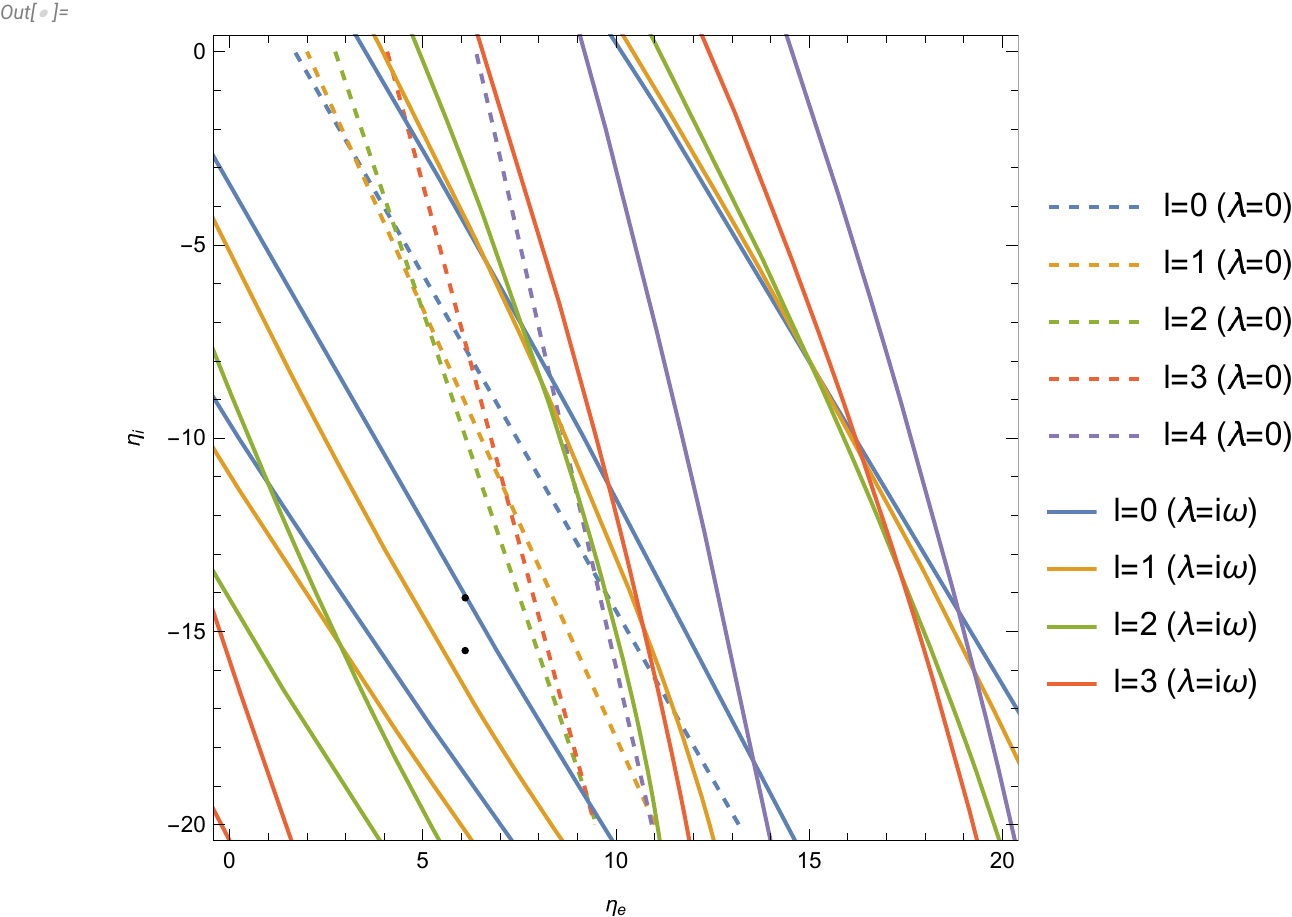}
    \caption{Bifurcation diagram in $\eta_{e}$ and $\eta_{i}$, up to $l=4$, with $d_e=0.02,d_i=0.2$. Bifurcations of Hopf-type are given by an unbroken curve, and those of fold-type by a dashed line. The points at $\eta_{e}=6.1$, $\eta_{i}=-14.134$ and $\eta_{i}=-15.5$ correspond to the parameter values of the Hopf bifurcation in Figure~\ref{fig:eigenvalues_l0} and of the simulation in Figure~\ref{fig:hopf_l0_time}, respectively.}
    \label{fig:bif_diag_l0}
\end{figure}

First, we set the diffusion in the inhibitory cells larger than in the excitatory cells. Our investigations show that in this case, an $l=0$ type Hopf bifurcation can be captured. In addition to the fixed parameters above, set $d_e = 0.02$, $d_i=0.2$, see Figure~\ref{fig:bif_diag_l0} for the bifurcation diagram. We observe that the fold-type bifurcations have moved closer to the $\eta_i$ axis, similarly for Hopf-type bifurcations with larger values of $l$. This makes the $l=0$ Hopf bifurcation curve clearly separated from the other Hopf bifurcations. Along this curve, at parameter values $\eta_{e}=6.1$, $\eta_{i}=-14.134$, we find that the eigenvalues and eigenvectors are $\lambda = \pm 0.802 i$ and $v_e=v_i=-0.652+0.275i$, respectively, see Figure~\ref{fig:eigenvalues_l0}. The normal form coefficient is $g_{0,1}=-0.336-0.030i$, hence the first Lyapunov coefficient can be computed $l_1=-0.419$. Since it is negative, this means that a stable limit cycle should exist locally, which is homogeneous in space, similar to the harmonic $Y_0^0$.

\begin{figure}%[H]
    \centering
    \includegraphics[scale=.9,trim={1cm 0 0 0},clip]{figures/eigenvalues_l0}
    \caption{Eigenvalues when an $l=0$ Hopf bifurcation occurs, with parameters $d_e=0.02$, $d_i=0.2$, $\eta_{e}=6.1$, $\eta_{i}=-14.134$. The critical eigenvalues and corresponding eigenvectors are $\lambda = \pm 0.802 i$, $v_e=v_i=-0.652+0.275i$, respectively. The normal form coefficient and first Lyapunov coefficient are $g_{0,1}=-0.336-0.030i$, $l_1=-0.419$.}
    \label{fig:eigenvalues_l0}
\end{figure}

\begin{figure}%[H]
     \centering
     \includegraphics[scale=0.7]{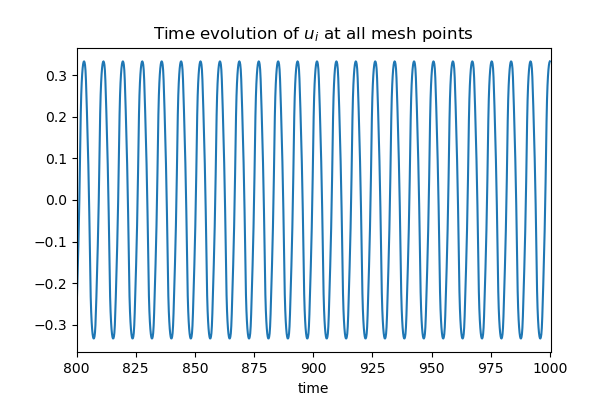}
     \caption{Time evolution of the $u_i$ state beyond a Hopf bifurcation with $l=0$, that oscillates homogeneously at all mesh points. Diffusion coefficients are $d_e = 0.02$, $d_i=0.2$ and the bifurcation parameters $\eta_{e}=6.1, \eta_{i}=-15.5$.}
     \label{fig:hopf_l0_time}
 \end{figure}

We simulate the full non-linear model  \eqref{eq:model_sphere} beyond the Hopf bifurcation, at $\eta_{e}=6.1$, $\eta_{i}=-15.5$ and observe a spatially homogeneous oscillation in both $u_e$ and $u_i$ states that oscillate with the same period and amplitude, see Figure~\ref{fig:hopf_l0_time}. The initial functions for this simulation are 
\[\varphi_e(t,r) = 0.1 \sin(\im(\lambda) t) Y_0^0(\br), \quad
  \varphi_i(t,r) = 0.1 \cos(\im(\lambda) t) Y_0^0(\br). \]

\subsection{Hopf bifurcation $l=1$}
\begin{figure}%[H]
    \centering
    \includegraphics[scale=.8,trim={.8cm 0 0 0},clip]{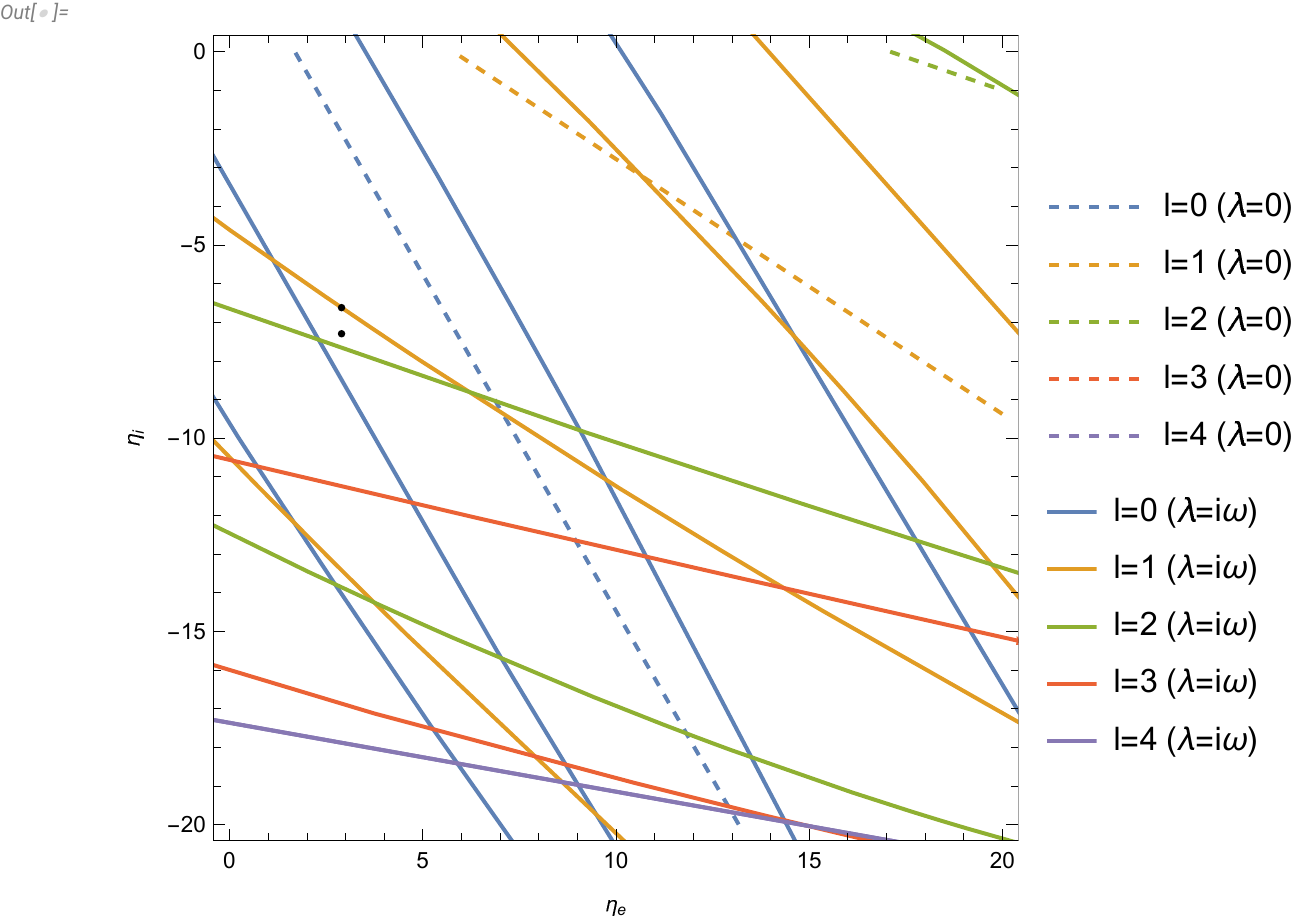}
    \caption{Bifurcation diagram in $\eta_{e}$ and $\eta_{i}$,  up to $l=4$, with $d_e=1$, $d_i=0.1$. Bifurcations of Hopf-type are given by an unbroken curve, and those of fold-type are given by a dashed curve. The points at $\eta_{e}=2.9$, $\eta_{i}=-6.624$ and $\eta_{i}=-7.3$ correspond to the parameter values of the Hopf bifurcation in Figure~\ref{fig:eigenvalues_l1} and of the simulation in Figures~\ref{fig:hopf-l1-Y11Y1m1},\ref{fig:hopf-l1-Y10}, respectively.}
    \label{fig:bif_diag_l1}
\end{figure}

Next, we consider the case where the diffusion in the excitatory cells is larger than in the inhibitory cells, $d_e=1,d_i=0.1$. In the bifurcation diagram in Figure~\ref{fig:bif_diag_l1}, we see the `opposite' pattern compared to the previous section. For a fixed value of $\eta_i$, we find that the fold-type bifurcations happen at a larger value of $\eta_e$, and the $l=3$ and $l=4$ fold-type bifurcations are not even present for this parameter range. Instead, we see a parameter range for which the $l=1$ Hopf bifurcation curve bounds the area where the trivial equilibrium is stable.

\begin{figure}%[H]
    \centering
    \includegraphics[scale=1,trim={0.8cm 0 0 0},clip]{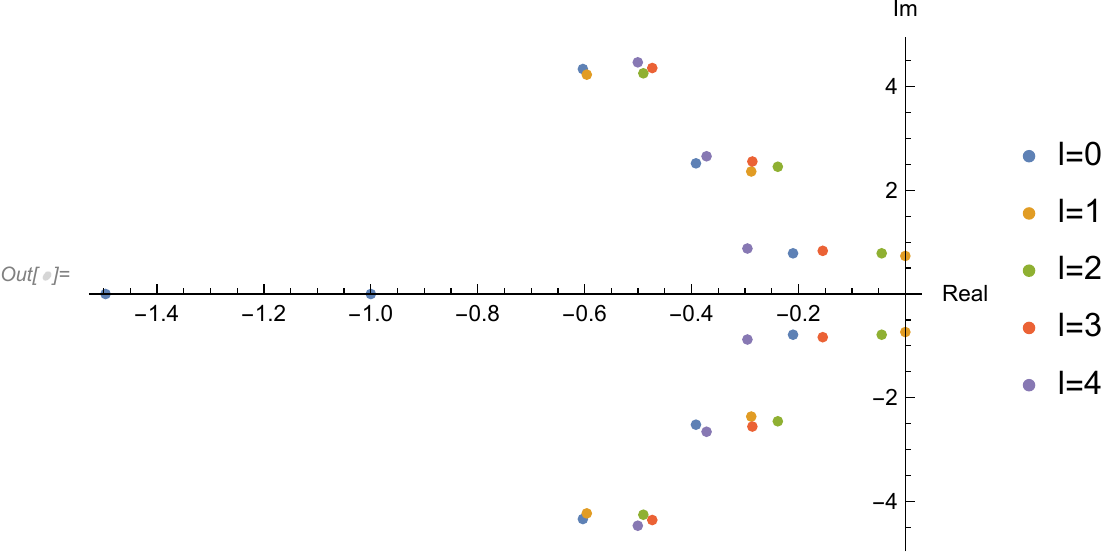}
    \caption{Eigenvalues at the $l=1$ Hopf bifurcation with parameters $d_e=1$, $d_i=0.1$,  $\eta_{e}=2.9$, $\eta_{i}=-6.624$. The eigenvalues and eigenvectors are $\lambda = \pm 0.734 i$, $v_e=-0.235-0.342i$, $v_i=-0.719-0.557i$, and normal form coefficients $g_{1,1}=-0.523+0.299i$, $g_{1,2}=-0.262+0.150i$.}
    \label{fig:eigenvalues_l1}
\end{figure}

We find an $l=1$ Hopf bifurcation with eigenvalues and eigenvectors $\lambda =\pm \omega i =\pm 0.734 i$, $v_e=-0.235-0.342i, v_i=-0.719-0.557i$, respectively, and normal form coefficients $g_{1,1}=-0.523+0.299i, g_{1,2}=-0.262+0.150i$ for parameters $\eta_{e}=2.9$, $\eta_{i}=-6.624$, see Figure~\ref{fig:eigenvalues_l1}. 

Our analysis of the amplitude equations \eqref{eq:amp_eq} in Section~\ref{sec:amp_eq} indicates that there should be a family of stable rotating waves when $\re(\mu)>0$, but the family of standing waves is unstable. 

%We can solve the normal form \eqref{eq:nf1} using a standard ODE-solver with $\mu=0.1+\lambda$, $g_{1,1},g_{1,2}$ as above and arbitrary initial conditions. We find that all these trajectories converge to a limit cycle where $\rho_1= 0.123$ and $\rho_2=\varphi=0$, as predicted.

We simulate the full nonlinear neural field system with $\eta_{e}=2.89$, $\eta_{i}=-7.3$.

\begin{figure}%[H]
\centering % <-- added
    \begin{subfigure}{0.45\textwidth}
    \includegraphics[width=\linewidth]{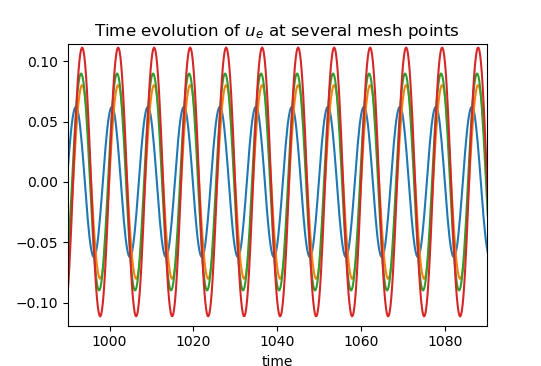}
    \caption{}
    %\label{fig:1}
    \end{subfigure}%\hfil % <-- added
    \begin{subfigure}{0.45\textwidth}
    \includegraphics[width=\linewidth]{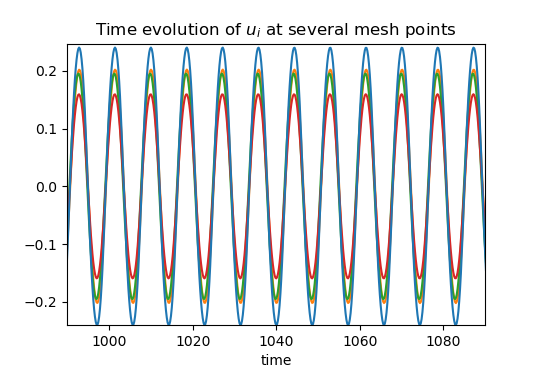}
    \caption{}
    %\label{fig:2}
    \end{subfigure}
    \begin{subfigure}{0.36\textwidth}
    \vspace{-1.5cm}
    \includegraphics[width=\linewidth]{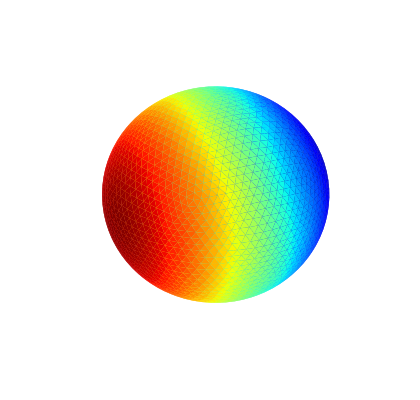}
    \vspace{-1.5cm}
    \caption{}
    %\label{fig:1}
    \end{subfigure}%\hfil % <-- added
\caption{Simulation of a rotating wave solution beyond an $l=1$ Hopf bifurcation with initial states as in \eqref{eq:initial_functions_l1_b}. (a) and (b) are the time evolution of the two states at some mesh points and (c) is the emerging solution pattern that rotates around the $z$-axis.}\label{fig:hopf-l1-Y11Y1m1}
\end{figure}

\paragraph{Rotating waves} Let us start with the initial functions 
     \begin{equation}\label{eq:initial_functions_l1_b}
     \varphi_e(t,r) = \varphi_i(t,r) = 0.1 \sqrt{2}\, \re(e^{i\omega t} Y_1^{-1}(\br)). 
     %=0.1 \left[ \cos(\omega t) Y_{1,1}+\sin(\omega t) Y_{1,-1}\right],
     \end{equation}
%where $Y_{1,1}$ and $Y_{1,-1}$ are the corresponding real spherical harmonics. 
In Figure~\ref{fig:hopf-l1-Y11Y1m1}, we can observe the time evolution of both states at some mesh points. They oscillate with the same frequency but larger amplitude in the $u_i$ component. This is a rotating wave solution, where Figure~\ref{fig:hopf-l1-Y11Y1m1} shows the pattern that travels around the $z$-axis.

When we slightly perturb the initial functions as
\begin{equation}\label{eq:initial_functions_l1_b_perturb}
     \varphi_e(t,r) = \varphi_i(t,r) = 0.1 \sqrt{2}\, \re(e^{i\omega t} Y_1^{-1}(\br)) + 10^{-4}\re{(e^{i\omega t}) Y_1^{0}(\br)}
     \end{equation}
the solution converges to the same rotating wave as in the unperturbed case, as expected from the theory.

\begin{figure}%[H]
\centering % <-- added
    \begin{subfigure}{0.32\textwidth}
    \includegraphics[width=\linewidth]{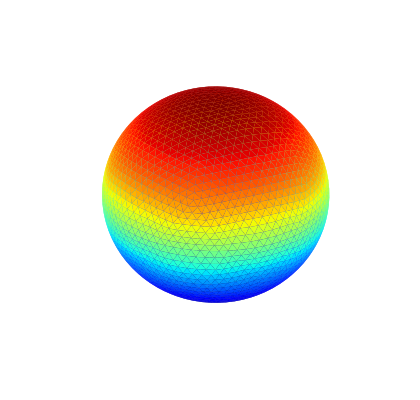}
    %\caption{$\xi=5\times 10^{-4}$}
    %\label{fig:1}
    \end{subfigure}%\hfil % <-- added
    \begin{subfigure}{0.32\textwidth}
    \includegraphics[width=\linewidth]{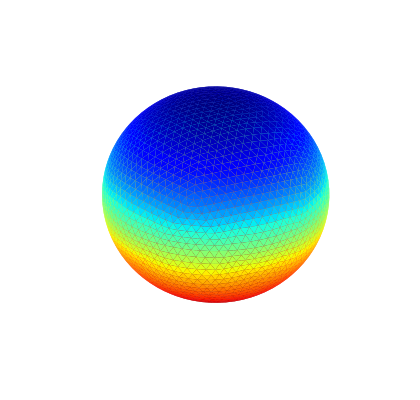}
    %\caption{$\xi=1\times 10^{-3}$}
    %\label{fig:2}
    \end{subfigure}%\hfil
\caption{Simulation of a standing wave solution beyond an $l=1$ Hopf bifurcation, with initial functions as in \eqref{eq:inital_functions_l1_a}. Snapshots at the beginning and at half of a time period.}\label{fig:hopf-l1-Y10}
\end{figure}

\paragraph{Standing waves} It is also possible to find a standing wave, when we choose the initial functions to be in the corresponding fixed point subspace 
     \begin{equation}\label{eq:inital_functions_l1_a}
         \varphi_e(t,r) = \varphi_i(t,r) = 0.1\, \re(e^{i\omega t}) Y_1^{0}(\br). 
     \end{equation}
The emerging standing wave solution is drawn in Figure~\ref{fig:hopf-l1-Y10}.

Using a small perturbation of the initial functions as 
\begin{equation}\label{eq:inital_functions_l1_a_perturbed}
         \varphi_e(t,r) = \varphi_i(t,r) = 0.1\, \re(e^{i\omega t}) Y_1^{0}(\br) + 10^{-4} \re(e^{i\omega t} Y_1^{-1}(\br)),
\end{equation}
we find that after a long transient the solution converges to a rotating wave solution similarly to what was obtained with \eqref{eq:initial_functions_l1_b}. We conjecture that the unstable manifold of the standing waves intersects with the stable manifold of the rotating waves. 
%Expectation: find a rotating wave

\subsection{Hopf bifurcation $l=2$} 
\begin{figure}%[H]
    \centering
    \includegraphics[scale=.8,trim={0.8cm 0 0 0},clip]{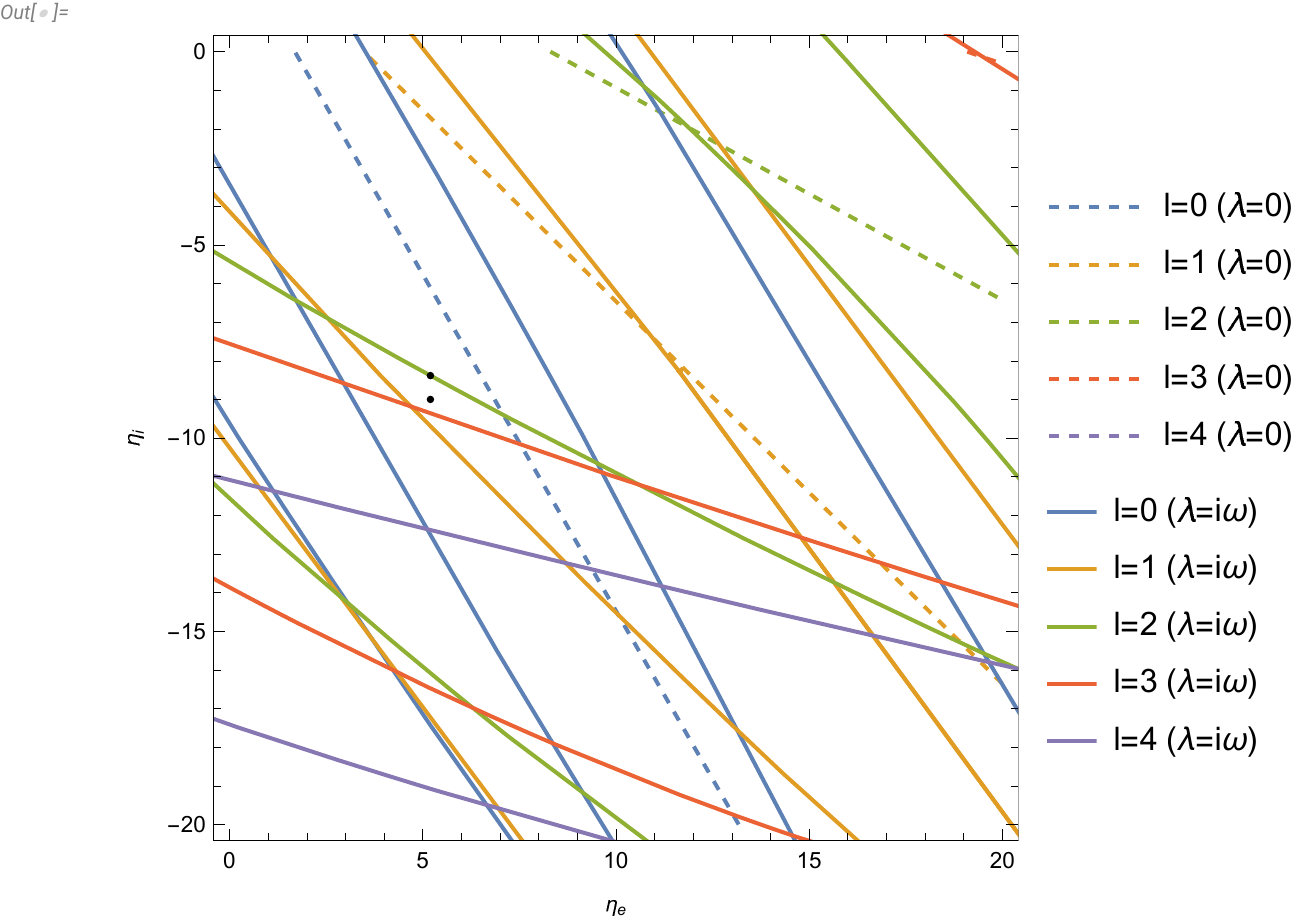}
    \caption{Bifurcation diagram in $\eta_{e}$ and $\eta_{i}$, up to $l=4$, with $d_e=0.4,d_i=0.04$. Bifurcations of Hopf-type are given by an unbroken curve, and those of fold-type are given by a dashed line. The points at $\eta_{e}=5.2$, $\eta_{i}=-8.384$ and $\eta_{i}=-9$ correspond to the parameter values of the Hopf bifurcation in Figure~\ref{fig:eigenvalues_l2} and of the simulation in Figures~\ref{fig:hopf-l2-tetrahedtal},\ref{fig:hopf-l2-mixed}, respectively.}
    \label{fig:bif_diag_l2}
\end{figure}

\begin{figure}%[H]
    \centering
    \includegraphics[scale=1,trim={0.8cm 0 0 0},clip]{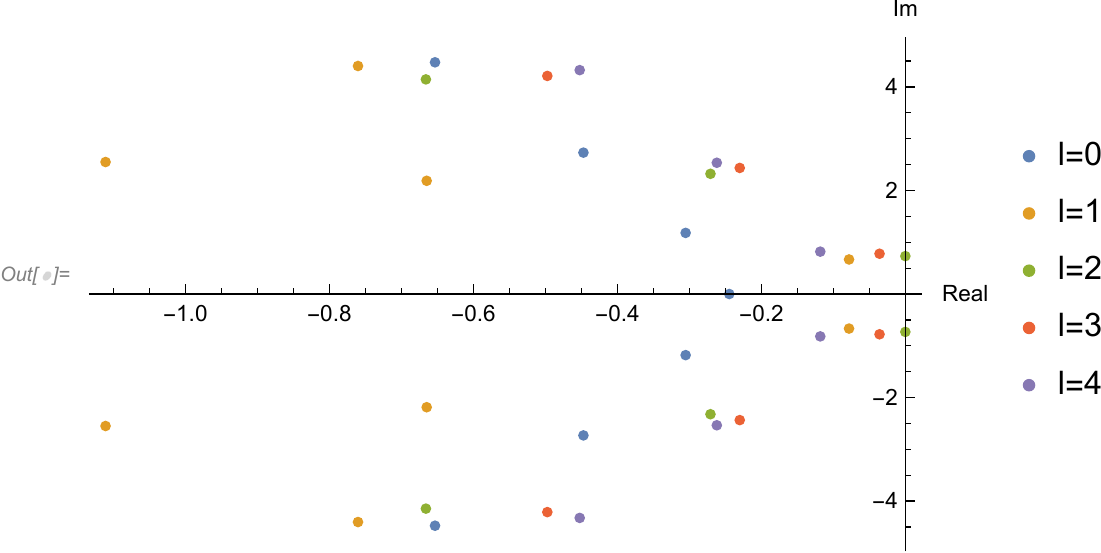}
    \caption{Eigenvalues at the $l=2$ Hopf bifurcation with parameters $\eta_{e}=5.2$, $\eta_{i}=-8.384$, $d_e=0.4$, $d_i=0.04$. Critical eigenvalues and corresponding eigenvectors are $\lambda = \pm 0.732 i$, $v_e=-0.224-0.310i$, $v_i=-0.750-0.540i$. Normal form coefficients are $g_{2,1}=-0.224-0.310i$, $g_{2,2}=-0.750-0.540i$, $g_{2,3}=0$.}
    \label{fig:eigenvalues_l2}
\end{figure}

Similar to the previous section, we consider the case where the diffusion in the excitatory cells is larger than in the inhibitory cells, but smaller magnitudes as before, $d_e = 0.4$, $d_i=0.04$. In the corresponding bifurcation diagram in Figure~\ref{fig:bif_diag_l2}, we now have an $l=2$ Hopf bifurcation curve that bounds the area where the trivial equilibrium is stable. Along this curve, we fix $\eta_{e}=5.2$, $\eta_{i}=-8.384$ and find the eigenvalues and eigenvectors $\lambda = \pm \omega i=\pm 0.723 i$, $v_e=-0.224-0.310i$, $v_i=-0.750+0.540i$, respectively, see Figure~\ref{fig:eigenvalues_l2}. The normal form coefficients are $g_{2,1}=-0.224-0.310i$, $g_{2,2}=-0.750-0.540i$, $g_{2,3}=0$.

According to \cite{iooss_hopf_1989}, $g_{2,3}=0$ is a degenerate case, so we change $\delta=0.1$ and $\gamma=10.332$. This leaves $S'(0)$ and hence the linearisation the same and therefore the eigenvalues too, but $S''(0)$ is now non-zero. For these parameters, the normal form coefficients are $g_{2,1}=-6.295-1.302i, g_{2,2}=-1.927-0.698i, g_{2,3}=-0.043-0.018i$. According to \cite{iooss_hopf_1989}, one maximal branch should be stable, which has a tetrahedral symmetry. In terms of the normal form variables, this branch can be described as $|z_1|=\sqrt{2}|z_{-2}|$ and $z_{-1}=z_{0}=z_2=0$. 

We simulate the full nonlinear neural field system with parameters beyond Hopf bifurcation as $\eta_{e}=5.2$, $\eta_{i}=-9$.  

\begin{figure}%[H]
    \centering
    \hspace{-1cm}
    \begin{subfigure}{0.3\textwidth}
    \includegraphics[width=\linewidth]{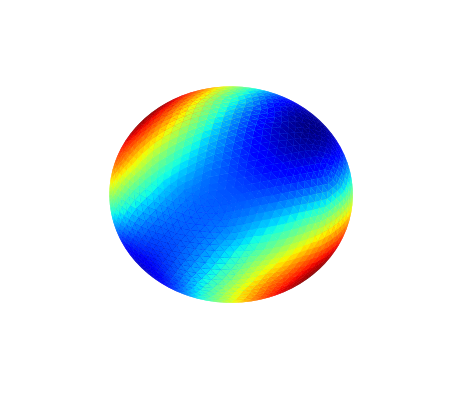}
    \end{subfigure}\hspace{-1cm}
    \begin{subfigure}{0.3\textwidth}
    \includegraphics[width=\linewidth]{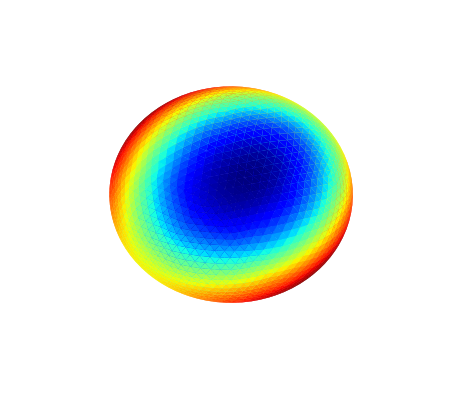}
    \end{subfigure}\hspace{-1cm}
    \begin{subfigure}{0.3\textwidth}
    \includegraphics[width=\linewidth]{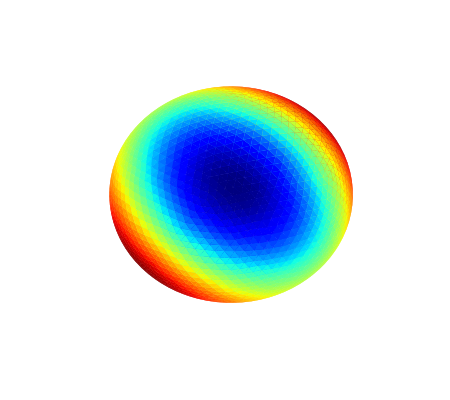}
    \end{subfigure}\hspace{-1cm}
    \begin{subfigure}{0.3\textwidth}
    \includegraphics[width=\linewidth]{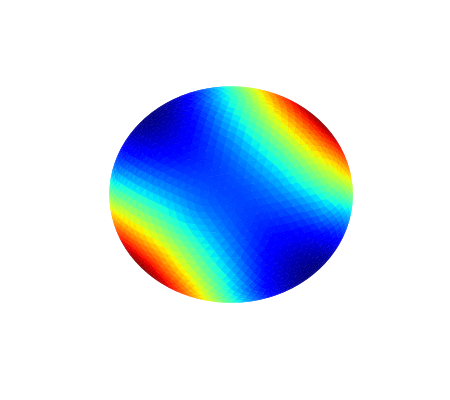}
    \end{subfigure}\\
    \vspace*{-1cm}
    \hspace{-1.7cm}
    \begin{subfigure}{0.3\textwidth}
    \includegraphics[width=\linewidth]{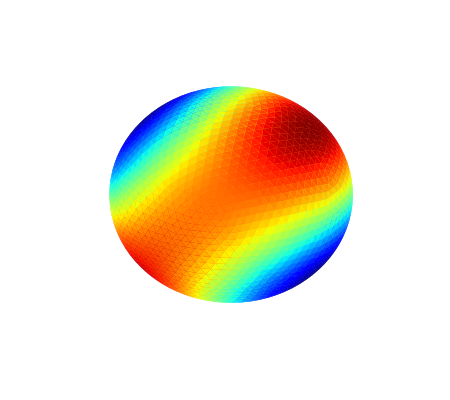}
    \end{subfigure}\hspace{-1cm}
    \begin{subfigure}{0.3\textwidth}
    \includegraphics[width=\linewidth]{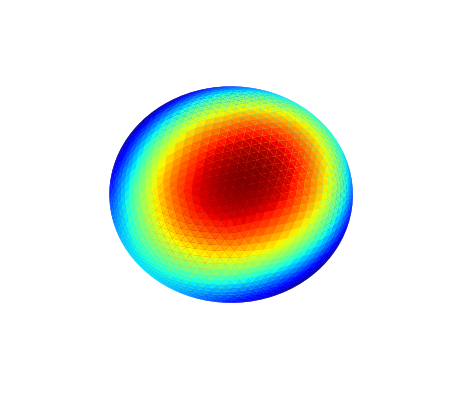}
    \end{subfigure}\hspace{-1cm}
    \begin{subfigure}{0.3\textwidth}
    \includegraphics[width=\linewidth]{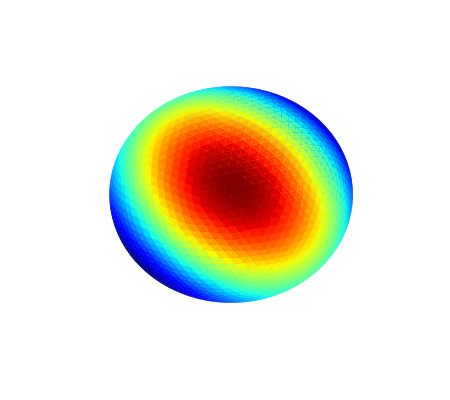}
    \end{subfigure}\hspace{-1cm}
    \begin{subfigure}{0.3\textwidth}
    \includegraphics[width=\linewidth]{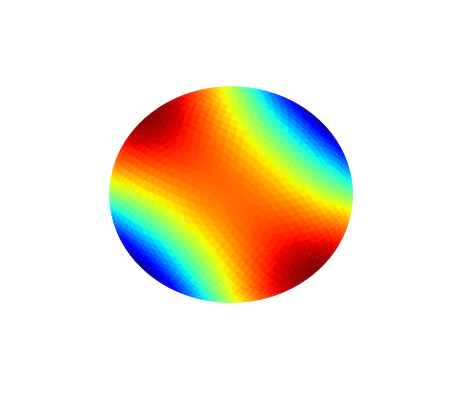}
    \end{subfigure}\hspace{-1cm}
    \vspace{-.5cm}
    \caption{Simulation of the solution beyond an $l=2$ Hopf bifurcation, with initial functions as in \eqref{eq:initial_solution_l2_a}. Snapshots of the resulting inhibitory solution with tetrahedral symmetry at times $nT/8$, $n=1,\cdots,8$, with time period $T$.}
    \label{fig:hopf-l2-tetrahedtal}
\end{figure}

\paragraph{Waves with tetrahedral symmetry} When the initial functions are 
     \begin{equation}\label{eq:initial_solution_l2_a}
         \varphi_e(t,r) = \varphi_i(t,r) = 0.1\, \re(e^{i\omega t} Y_2^{-2}(\br)) + 0.1 \sqrt{2} \, \re(e^{i\omega t} Y_2^{1}(\br))
     \end{equation}
    then we find the branch with tetrahedral symmetry. The solution pattern during 1/8 of a period can be observed in Figure~\ref{fig:hopf-l2-tetrahedtal}. This is not a standing wave, nor does it rotate in a single direction. 
    
    Next, we slightly perturb the initial condition, so it falls outside the invariant subspace
      \begin{equation}\label{eq:initial_solution_l2_small_perturb}
      \begin{split}
         \varphi_e(t,r) = \varphi_i(t,r) =& 0.1\, \re(e^{i\omega t} Y_2^{-2}(\br)) + 0.1\sqrt{2} \, \re(e^{i\omega t} Y_2^{1}(\br))\\
         &+ 10^{-4} \re(e^{i\omega t} Y_2^{-1}(\br)) + 10^{-4} \re(e^{i\omega t} Y_2^{0}(\br)) + 10^{-4} \re(e^{i\omega t} Y_2^{2}(\br)). 
      \end{split}
     \end{equation}    
     We find that the numerical solution converges to the same solution as Figure~\ref{fig:hopf-l2-tetrahedtal}. This matches with the theoretical prediction that this branch should be (locally) stable.
    
 %    \begin{figure}[H]
 %     \centering
 %     \hspace{-.5cm}
 %     \includegraphics[scale=0.4]{figures/Tetrahedral/Figure_1.png}
 %     \hspace{-1cm}
 %     \includegraphics[scale=0.4]{figures/Tetrahedral/Figure_2.png}
 %     \hspace{-1cm}
 %     \includegraphics[scale=0.4]{figures/Tetrahedral/Figure_3.png}
 %     \hspace{-1cm}
 %     \includegraphics[scale=0.4]{figures/Tetrahedral/Figure_4.png}\\
 %     \hspace{-.5cm}
 %     \includegraphics[scale=0.4]{figures/Tetrahedral/Figure_5.png}
 %     \hspace{-1cm}
 %     \includegraphics[scale=0.4]{figures/Tetrahedral/Figure_6.png}
 %     \hspace{-1cm}
 %     \includegraphics[scale=0.4]{figures/Tetrahedral/Figure_7.png}
 %     \hspace{-1cm}
 %     \includegraphics[scale=0.4]{figures/Tetrahedral/Figure_8.png}
 %     \caption{Simulation of the solution beyond an $l=2$ Hopf bifurcation, with initial functions as in \eqref{eq:initial_solution_l2_a}. Snapshots of the resulting inhibitory solution with tetrahedral symmetry at times $nT/8$, $n=1,\cdots,8$, with time period $T$.}
 %     \label{fig:hopf-l2-tetrahedtal}
 % \end{figure}
     
\begin{figure}%[H]
    \centering
    \hspace{-1cm}
    \begin{subfigure}{0.3\textwidth}
    \includegraphics[width=\linewidth]{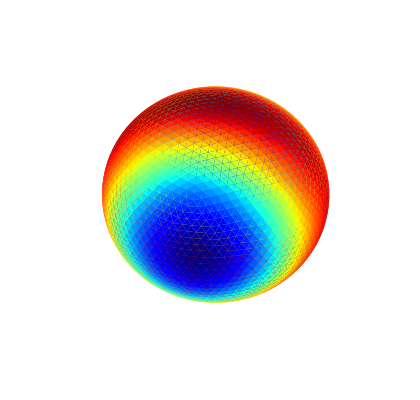}
    \end{subfigure}\hspace{-1cm}
    \begin{subfigure}{0.3\textwidth}
    \includegraphics[width=\linewidth]{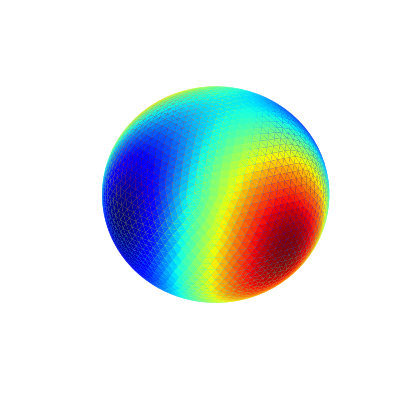}
    \end{subfigure}\hspace{-1cm}
    \begin{subfigure}{0.3\textwidth}
    \includegraphics[width=\linewidth]{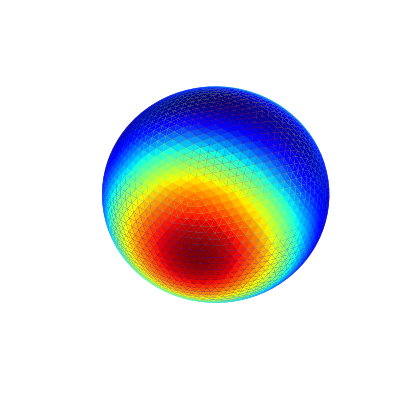}
    \end{subfigure}\hspace{-1cm}
    \begin{subfigure}{0.3\textwidth}
    \includegraphics[width=\linewidth]{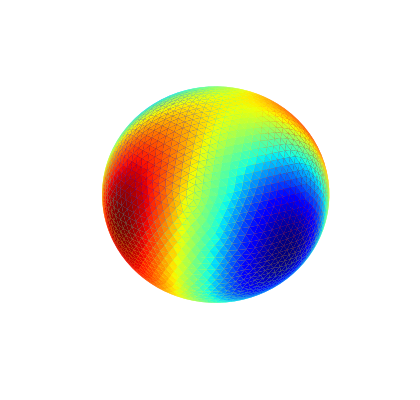}
    \end{subfigure}\hspace{-1cm}
    \vspace{-.5cm}
    \caption{Simulation of the solution beyond an $l=2$ Hopf bifurcation, with initial functions as in \eqref{eq:initial_solution_l2_perturbed}. Snapshots of the resulting inhibitory solution, called mixed waves, at times $nT/4$, $n=1,\cdots,4$, with time period $T$. }
     \label{fig:hopf-l2-mixed}
\end{figure}

\paragraph{Mixed waves} Next we examine an initial condition which is further away from \eqref{eq:initial_solution_l2_a}, hence not in a fixed-point subspace corresponding to any of the maximal branches as identified by \cite{iooss_hopf_1989}
\begin{equation}\label{eq:initial_solution_l2_perturbed}
\begin{split}
    \varphi_e(t,r) = \varphi_i(t,r) &= e^{-1}\, \re(e^{i\omega t} Y_2^{-2}(\br)) + e^{-\tfrac{1}{2}} \, \re(e^{i\omega t} Y_2^{-1}(\br)) \\
    &+ \re(e^{i\omega t} Y_2^{0}(\br)) + e^{\tfrac{1}{2}} \re(e^{i\omega t} Y_2^{1}(\br)) + e^2 \re(e^{i\omega t} Y_2^{2}(\br)).
\end{split}
\end{equation}  
The numerical solution is a wave that travels around the sphere as in Figure~\ref{fig:hopf-l2-mixed}. This solution does not correspond to any of the maximal branches mentioned in \cite{iooss_hopf_1989}. We conjecture that for these parameter values, the neural field is bi-stable with one maximal and one submaximal branch.

\subsection{Hopf bifurcation $l=3$}
\begin{figure}%[H]
    \centering
    \includegraphics[scale=.8,trim={0.6cm 0 0 0},clip]{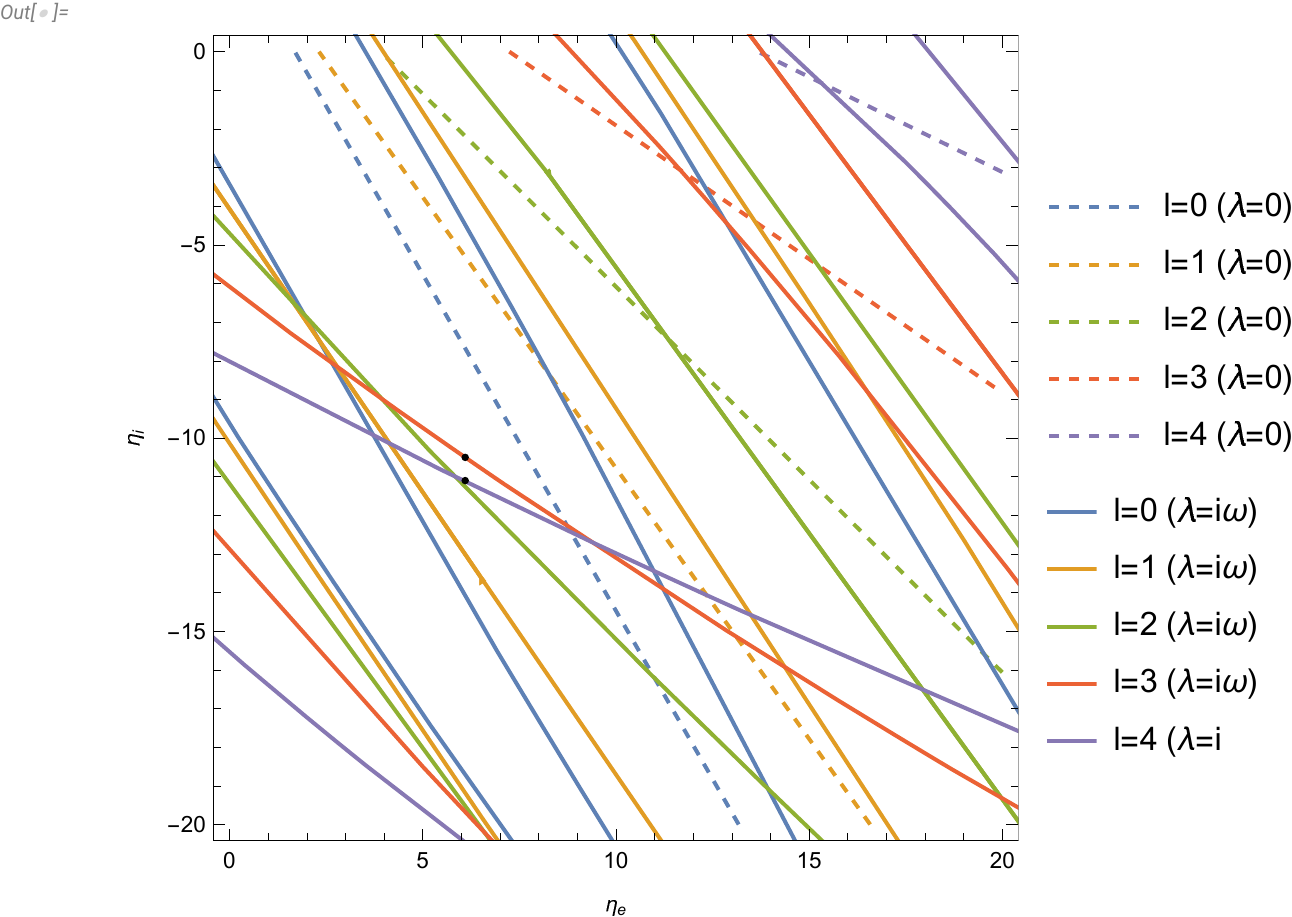}
    \caption{Bifurcation diagram in $\eta_{e}$ and $\eta_{i}$, up to $l=4$, with $d_e=0.1,d_i=0.01$. Bifurcations of Hopf-type are given by an unbroken curve, and those of fold-type are given by a dashed line. The points at $\eta_{e}=6.1$, $\eta_{i}=-10.500$ and $\eta_{i}=-11.1$ correspond to the parameter values of the Hopf bifurcation in Figure~\ref{fig:eigenvalues_l3} and of the simulation in Figures~\ref{fig:hopf-l3-Y31Y3m1},\ref{fig:Hopf-l3-rotating-perturbed}, respectively.}
    \label{fig:bif_diag_l3}
\end{figure}

\begin{figure}%[H]
    \centering
    \includegraphics[scale=1,trim={0.8cm 0 0 0},clip]{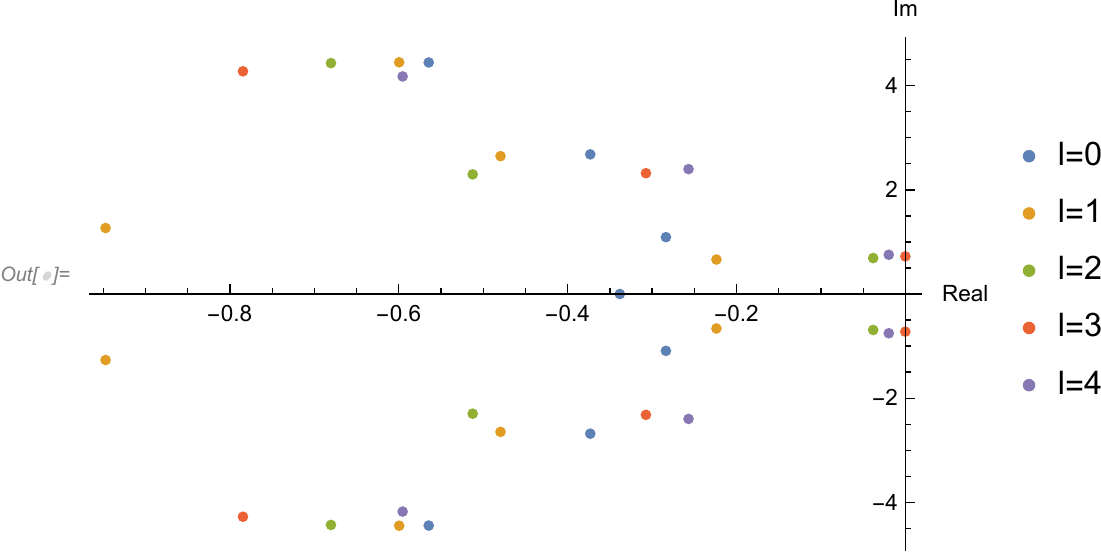}
    \caption{Eigenvalues at $l=3$ Hopf bifurcation with parameters $\eta_{e}=6.1$, $\eta_{i}=-10.500$, $d_e=0.1$, $d_i=0.01$. Critical eigenvalues and corresponding eigenvectors $\lambda = \pm \omega i=\pm 0.723 i$, $v_e=-0.496-0.049i$, $v_i=-0.856+0.135i$, respectively. Normal form coefficients are $g_{3,1}=-1.131-0.344i$, $g_{3,2}=-0.566-0.172i$, $g_{3,3}=-0.026-0.0078i$, $g_{3,4}=0.043+0.013i$.}
    \label{fig:eigenvalues_l3}
\end{figure}

If we reduce the magnitude of the diffusion even further to $d_e = 0.1$, $d_i=0.01$, the curve of $l=3$ Hopf bifurcations now bounds the area where the trivial equilibrium is stable, see Figure~\ref{fig:bif_diag_l3}. We find an $l=3$ Hopf bifurcation at $\eta_{e}=6.1$, $\eta_{i}=-10.500$, with eigenvalues and eigenvectors $\lambda = \pm 0.723 i$, $v_e=-0.496-0.049i$, $v_i=-0.856+0.135i$, and normal form coefficients $g_{3,1}=-1.131-0.344i$, $g_{3,2}=-0.566-0.172i$, $g_{3,3}=-0.026-0.0078i$, $g_{3,4}=0.043+0.013i$, see Figure~\ref{fig:eigenvalues_l3}. 

This corresponds to case (v) described in Section 5.2 in  \cite{sigrist_hopf_2010}: All branches of Figure~\ref{fig:branches} bifurcate supercritically, but all have at least one unstable direction. Continuing this Hopf bifurcation line until the point where the different Hopf curves cross and recalculating the normal form coefficients, we find that the stability of the branches does not change. Because all the branches bifurcate supercritically and all other eigenvalues have negative real part, the trivial equilibrium is not exponentially stable. On the side of the bifurcation where the equilibrium becomes unstable, we expect a stable (quasi-)periodic corresponding to a submaximal branch.

%$\eta_{e}=6.1$,  $\eta_{i}=-10.8$. 

Next, we give some numerical examples of periodic solutions corresponding to maximal branches that exist in the neural field \eqref{eq:model_sphere}. These are guaranteed to exist by the equivariant Hopf theorem, but as previously established unstable. However, we can find them by starting the simulation of the full nonlinear system with initial functions in the corresponding fixed-point subspace, \cite[Section 5.3.2]{sigrist_hopf_2010}.

\begin{figure}%[H]
     \centering
     \includegraphics[scale=0.52]{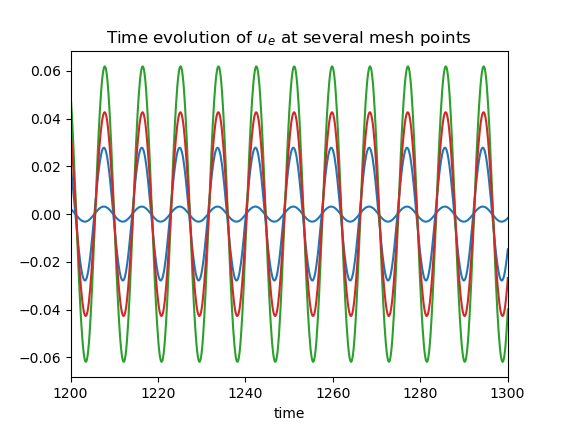}
     \includegraphics[scale=0.52]{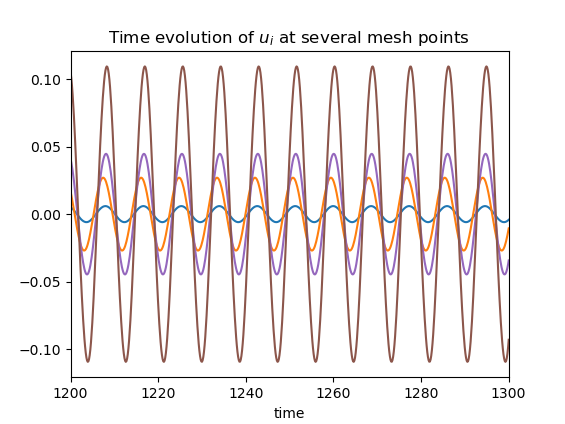}
     \includegraphics[scale=0.52]{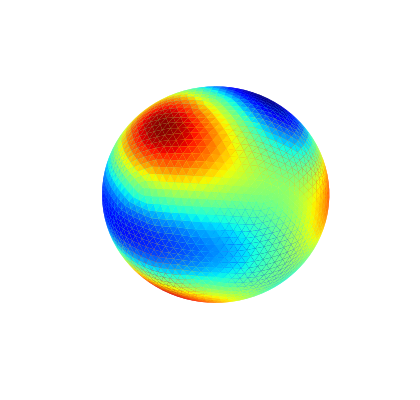}
     \caption{Time simulation beyond an $l=3$ Hopf bifurcation with initial functions in \eqref{eq:initial_solution_l3} (top). Rotating wave solution pattern that rotates around the $z$-axis (bottom).}
     \label{fig:hopf-l3-Y31Y3m1}
 \end{figure}

\paragraph{Rotating waves} When the initial functions are 
     \begin{equation}\label{eq:initial_solution_l3}
         \varphi_e(t,r) = \varphi_i(t,r) =s_0\re \left(e^{i\omega t} Y_3^{-1}(\br)\right), 
     \end{equation}
then a rotating wave solution emerges, see Figure~\ref{fig:hopf-l3-Y31Y3m1}, where the pattern rotates around the $z$-axis. This corresponds to the solution in Figure~\ref{fig:branches} (b).
 
\begin{figure}%[H]
     \centering
     \hspace{-.5cm}
     \includegraphics[scale=0.4]{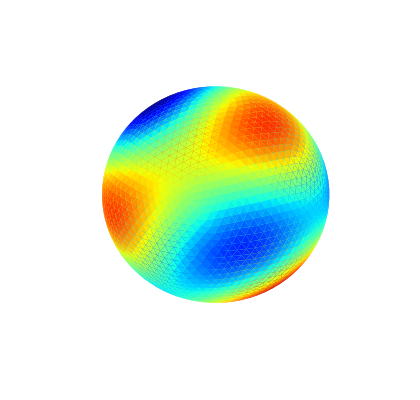}
     \hspace{-1cm}
     \includegraphics[scale=0.4]{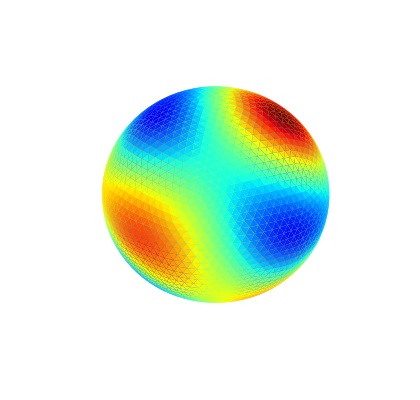}
     \hspace{-1cm}
     \includegraphics[scale=0.4]{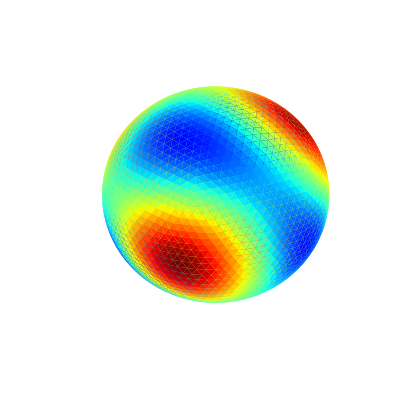}
     \hspace{-1cm}
     \includegraphics[scale=0.4]{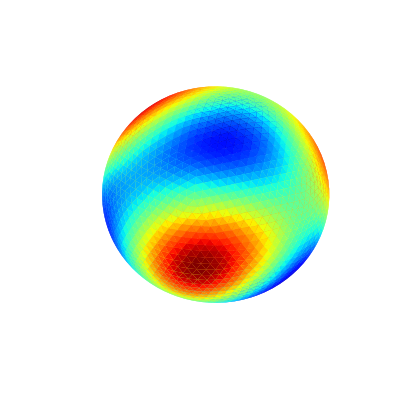}\\ \vspace*{-1.cm}
     \hspace{-.5cm}
     \includegraphics[scale=0.4]{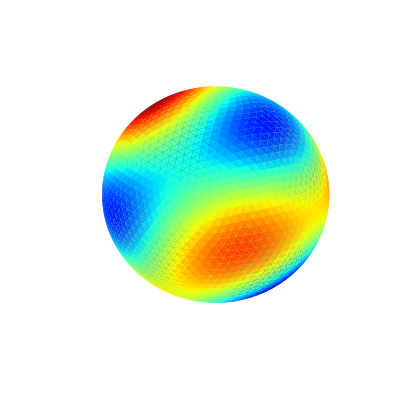}
     \hspace{-1cm}
     \includegraphics[scale=0.4]{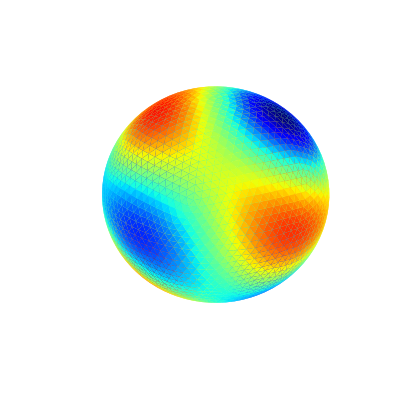}
     \hspace{-1cm}
     \includegraphics[scale=0.4]{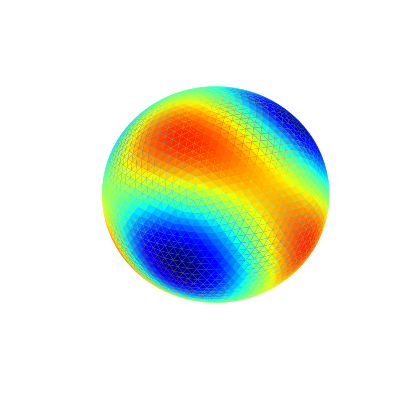}
     \hspace{-1cm}
     \includegraphics[scale=0.4]{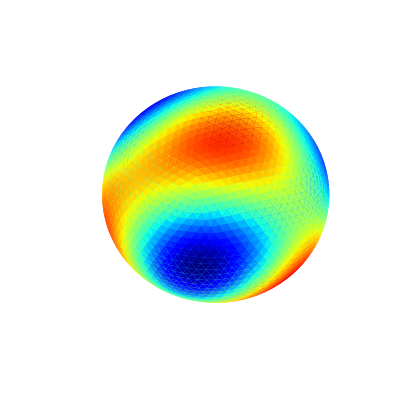}
     \vspace{-.5cm}
     \caption{Simulation of the solution beyond an $l=3$ Hopf bifurcation, with initial functions in \eqref{eq:initial_solution_l3-perturb}. Snapshots of the resulting inhibitory solution at times $nT/8$, $n=1,\cdots,8$, with time period $T$.}
     \label{fig:Hopf-l3-rotating-perturbed}
 \end{figure}

    \paragraph{Mixed waves} By perturbing slightly the initial functions as 
    \begin{equation}\label{eq:initial_solution_l3-perturb}
    \begin{split}
        \varphi_e(t,r) &=0.1 \re \left(e^{i\omega t} Y_3^{-1}(\br)\right)+ 10^{-3}\, \re \left(e^{i\omega t} Y_3^{-2}(\br)\right)\\
        \varphi_i(t,r) &=0.1 \re \left(e^{i\omega t} Y_3^{-1}(\br)\right) + 5\cdot 10^{-3}\, \re \left(e^{i\omega t} Y_3^{-3}(\br)\right),
    \end{split}
    \end{equation}
    the pattern in Figure~\ref{fig:Hopf-l3-rotating-perturbed} emerges. Although, it is a rotating wave, it does not rotate in a single direction. We conjecture that this solution corresponds to a submaximal branch, but it is still unclear which one. 
    
\subsection{Stationary spiral patterns}
Next to standing and rotating waves that arise from Hopf bifurcations with spherical symmetry, also stationary spiral patterns can arise from bifurcations with zero eigenvalue. In this section, we study the formation of spiral patterns in our model \eqref{eq:model_sphere}. In \cite{sigrist_symmetric_2011}, the authors give some examples of bifurcations with a double zero eigenvalue, where $\mathcal{E}_l(0)=\mathcal{E}_{l+1}(0)=0$.

\begin{figure}%[H]
    \centering
    \includegraphics[scale=.9,trim={0.8cm 0 0 0},clip]{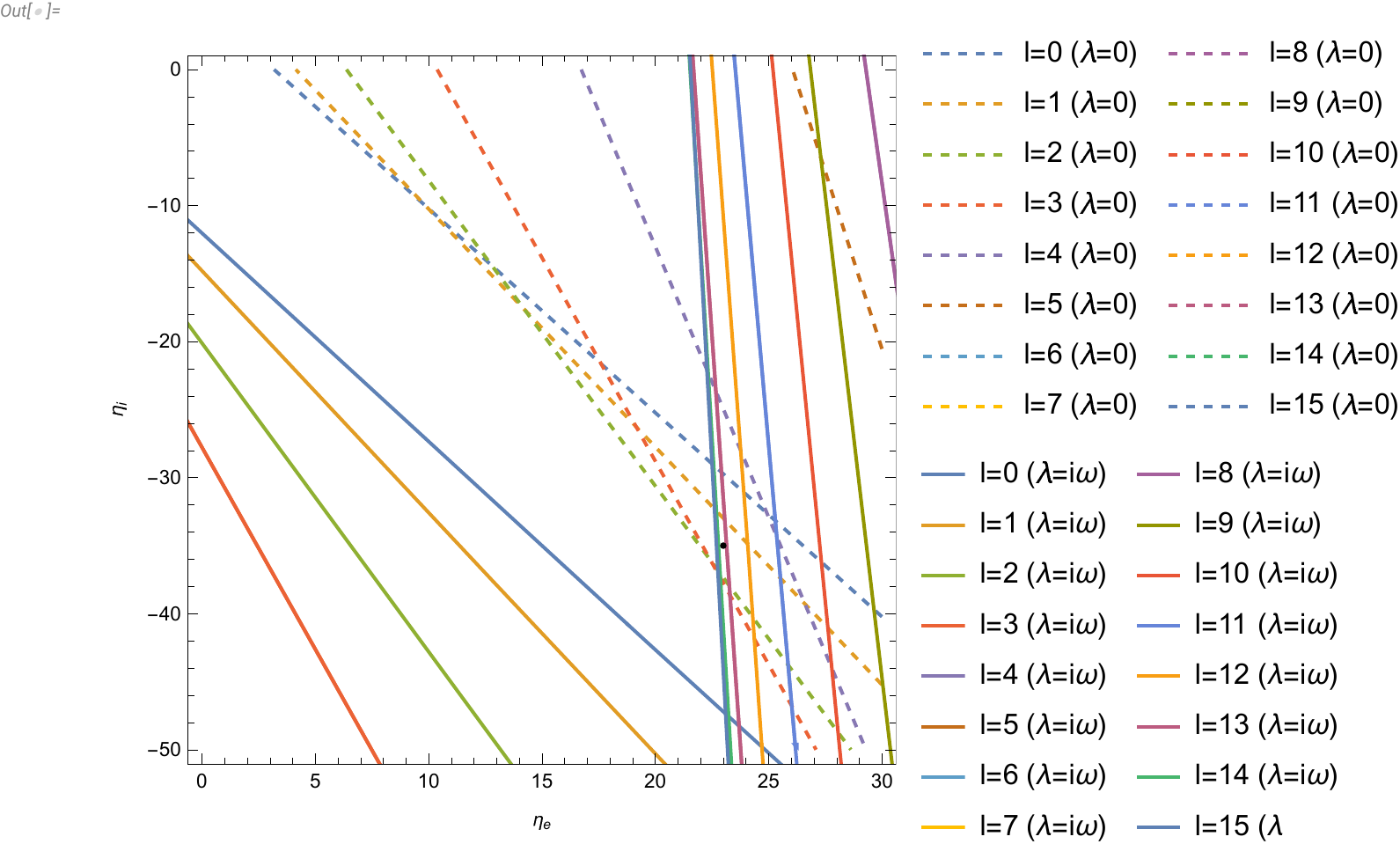}
    \caption{Bifurcation diagram in $\eta_{e}$ and $\eta_{i}$, up to $l=15$, with the parameters \eqref{eq:spiral_param}. Bifurcations of Hopf-type are given by an unbroken curve, and those of fold-type are given by a dashed line. The point at $\eta_{e}=23$, $\eta_{i}=-35$ corresponds to the parameter values of the simulation in Figures~\ref{fig:hopf-l3-Y31Y3m1}, \ref{fig:Hopf-l3-rotating-perturbed}, respectively.}
    \label{fig:bif_diag_spiral}
\end{figure}

First, we computed the bifurcation diagram in Figure~\ref{fig:bif_diag_spiral} with parameters
\begin{equation}\label{eq:spiral_param}
    \alpha_e=4,\ \alpha_i=6,\ \sigma_{e}=1/3,\  \sigma_{i}=1/3,\ \tau_0=3,\ c=0.1,\  \gamma=8,\ \delta = 0,\ d_e =0.001,\ d_i=0.5,
\end{equation}
for which we have both $l=2$ and $l=3$ fold-type bifurcation curves that bound the area where the trivial equilibrium is stable. At their intersection point, at $\eta_e=23$ and $\eta_i=-35$, we find the eigenvalue $\lambda=0$, corresponding to these bifurcations and all the other eigenvalues have negative real part. 

\begin{figure}%[H]
     \centering
     \includegraphics[scale=0.42]{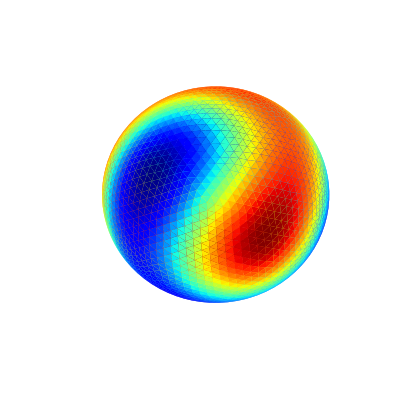}
     \includegraphics[scale=0.42]{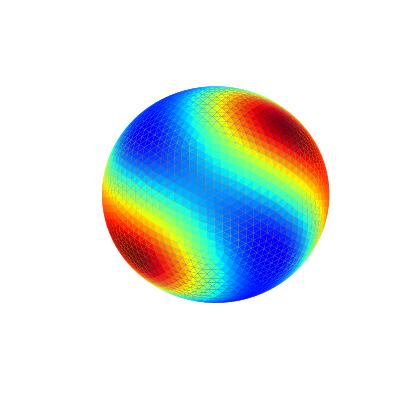}
     \vspace{-.5cm}
     \caption{Two views of a stationary spiral pattern, $u_e$ solution component.}
     \label{fig:Spiral-1}
 \end{figure}

By setting $\eta_e=23$ and $\eta_i=-35$, which lie beyond a fold-type bifurcation, we simulate the nonlinear system with initial functions
\begin{align}\label{eq:initial_solution_spiral1}
    \varphi_e(t,r) &= \sum_{m=-2}^2c_m^e\re \left(Y_2^{m}(\br)\right)+\sum_{m=-3}^3 d_m^e\re \left(Y_3^{m}(\br)\right), \\
    \varphi_i(t,r) &= \sum_{m=-2}^2c_m^i\re \left(Y_2^{m}(\br)\right)+\sum_{m=-3}^3 d_m^i\re \left(Y_3^{m}(\br)\right),
\end{align}
where $c_m^e,d_m^e, c_m^i,d_m^i\in\mathbb{R}$. The numerical solution converges to the stationary spiral pattern in Figure~\ref{fig:Spiral-1}, which we identified with the solution to the system (6.3)–(6.5) plotted in Figure~6 (last row) in \cite{sigrist_symmetric_2011}.  

\begin{figure}%[H]
     \centering
     \includegraphics[scale=0.42]{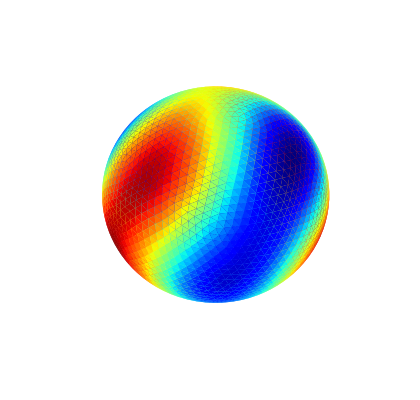}
     \includegraphics[scale=0.42]{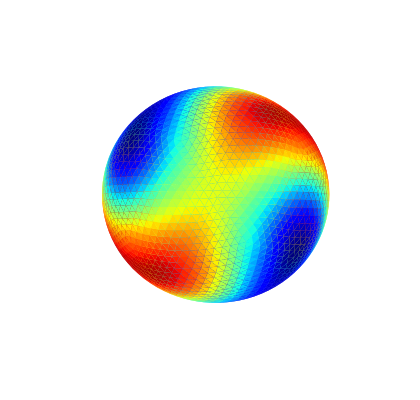}
     \vspace{-.5cm}
     \caption{Two views of a stationary spiral pattern, $u_i$ solution component.}
     \label{fig:Spiral-2}
 \end{figure}

Setting the initial function in \eqref{eq:initial_solution_spiral1} $\varphi_e=\varphi_i$, the solution converges to the spiral pattern in Figure~\ref{fig:Spiral-2}, again identified with the solution of the same problem in \cite{sigrist_symmetric_2011}, plotted in Figure~6 (images $\widetilde{{\mathbf D}_4}$ (i) and (ii)).  

\section{Discussion}
This paper presents a case study of local bifurcations leading to oscillatory behaviour in a distributed dynamical system with local and non-local interactions, caused by intrinsic dynamics, diffusion, and propagation delays. While the general strategy to study such phenomena is well-known and is based on the center manifold reduction and the computation of the normal forms, its actual implementation in the presence of spherical symmetry proved to be rather involved and non-trivial. Hopefully, our analysis will serve as a successful example, helping to study other such classes of dynamical systems appearing in applications.

We have analysed the behaviour of the neural field model on the sphere for different values of the diffusion constants between excitatory and inhibitory neurons. In particular, we examined equivariant Hopf bifurcations in this spherically symmetric setting. We have mostly looked at the case where $d_i<d_e$, where the Hopf bifurcations of order $l=1,2,3$ become the primary bifurcation. However, in some biological models, the diffusion is mostly between inhibitory cells \cite{schwab_synchronization_2014}. Studying the bifurcation diagram for this case, we find that bifurcations with a zero eigenvector, which requires a different bifurcation analysis. According to \cite{sigrist_symmetric_2011}, patterns like spiral waves can become stable as we have also found numerically, see Figures~\ref{fig:Spiral-1}, \ref{fig:Spiral-2}. 

In section~\ref{sec:branch} we considered the different maximal branches that can arise from Hopf bifurcations with spherical symmetry for $l=2,3$. In section~\ref{sec:numerics} we could indeed observe these branches in the simulations of the model \eqref{eq:model_sphere} when they were initialised in certain subspaces. However, when starting from arbitrary initial conditions, we find some other stable solutions. In fact, the theory predicts that for the computed normal form coefficients for $l=3$, that none of the maximal branches are stable. We conjecture that the cycles we observe in the numerics correspond to submaximal branches of the equivariant Hopf bifurcation. It is still an open question to characterise when these submaximal branches are stable. 

\section*{Acknowledgements}
The authors would like to thank Sid Visser for his help with the numerical discretization used in this paper and its Python implementation. 

\FloatBarrier

\appendix
\begingroup
\allowdisplaybreaks
\section{The $O(3)\times S^1$ normal forms for $l=2,3$}
\label{sec:app_normal_form}
\subsection{For $l=2$}\label{sec:app_normal_form_l2}
When $l=2$, we find the normal form 
\[\dot{\mathbf{z}} = \lambda \mathbf{z} + g_{2,1}\mathbf{z}|\mathbf{z}|^2 + g_{2,2}\hat{\mathbf{z}} P^2(\mathbf{z}) + g_{2,3} \mathbf{C}(\mathbf{z}) + \mathcal{O}(|\mathbf{z}|^5)\]
where $P^2(\mathbf{z})=z_0^2-2z_{-1}z_1+2z_{-2}z_{2}$ and  $\mathbf{C}(\mathbf{z})=(\mathbf{C}_{-2}(\mathbf{z}),\mathbf{C}_{-1}(\mathbf{z}),\mathbf{C}_{0}(\mathbf{z}),\mathbf{C}_{1}(\mathbf{z}),\mathbf{C}_{2}(\mathbf{z}))$,
\begin{align*}
\mathbf{C}_{-2}(\mathbf{z})=& -\sqrt{6} z_{-2} \left(\tfrac{2}{3} |z_{0}|^2 + |z_1|^2 + |z_2|^2 \right) + z_{-1}^2 \overline{z}_0 + z_{-1}z_{0}\overline{z}_1 + \tfrac{1}{\sqrt{6}} z_0^2 \overline{z}_2 \\
\mathbf{C}_{-1}(\mathbf{z})=& -\sqrt{\frac{3}{2}} z_{-1}\left(|z_{-1}|^2+\frac{1}{3}|z_0|^2+|z_1|^2+2 |z_{2}|^2\right) + \sqrt{\tfrac{2}{3}}z_0^2\overline{z}_1 + z_0z_1\overline{z}_2\\
&+z_{-2}z_1\overline{z}_0 + 2z_{-2}z_0\overline{z}_{-1} \\
\mathbf{C}_{0}(\mathbf{z})=& -\tfrac{1}{\sqrt{6}} z_0(4|z_{-2}|^2+|z_{-1}|^2+3|z_0|^2+|z_1|^2+4|z_2|^2) + 
\end{align*}
and $\mathbf{C}_m(\mathbf{z}) = \mathbf{C}_{-m}(\tilde{\mathbf{z}})$, where $\tilde{\mathbf{z}} = (z_2,z_1,z_0,z_{-1},z_{-2})$. 

\subsection{For $l=3$}\label{sec:app_normal_form_l3}
When $l=3$, we find the normal form as 
\[\dot{\mathbf{z}} = \lambda \mathbf{z} + g_{3,1}\mathbf{z}|\mathbf{z}|^2 + g_{3,2}\hat{\mathbf{z}} P^3(\mathbf{z}) + g_{3,3} \mathbf{Q}(\mathbf{z})+ g_{3,4} \mathbf{R}(\mathbf{z}) + \mathcal{O}(|\mathbf{z}|^5),\]
where $P^2(\mathbf{z})=z_0^2-2z_{-1}z_1+2z_{-2}z_{2}-2z_{-3}z_{3}$ and where $\mathbf{Q}(\mathbf{z})$ and $\mathbf{R}(\mathbf{z})$ 
\begin{align*}
\mathbf{Q}(\mathbf{z}) =&(\mathbf{Q}_{-3}(\mathbf{z}),\mathbf{Q}_{-2}(\mathbf{z}),\mathbf{Q}_{-1}(\mathbf{z}),\mathbf{Q}_{0}(\mathbf{z}),\mathbf{Q}_{1}(\mathbf{z}),\mathbf{Q}_{2}(\mathbf{z}), \mathbf{Q}_{3}(\mathbf{z}))\\
\mathbf{R}(\mathbf{z}) =&(\mathbf{R}_{-3}(\mathbf{z}),\mathbf{R}_{-2}(\mathbf{z}),\mathbf{R}_{-1}(\mathbf{z}),\mathbf{R}_{0}(\mathbf{z}),\mathbf{R}_{1}(\mathbf{z}),\mathbf{R}_{2}(\mathbf{z}), \mathbf{R}_{3}(\mathbf{z}))
\end{align*}
and $\mathbf{Q}_m(\mathbf{z}) = \mathbf{Q}_{-m}(\tilde{\mathbf{z}}), \mathbf{R}_m(\mathbf{z}) = \mathbf{R}_{-m}(\tilde{\mathbf{z}})$ where $\tilde{\mathbf{z}} = (z_3,z_2,z_1,z_0,z_{-1},z_{-2},z_{-3})$. Moreover,
\begin{align*}
\mathbf{Q}_{-3}(\mathbf{z}) =& 5 z_{-3}(5|z_{-3}|^2+5|z_{-2}|^2 -|z_{-1}|^2 -4|z_0|^2 -5|z_1|^2-5|z_2|^2-8|z_3|^2) + 10z_0^2\overline{z}_3\\
&-15z_1z_{-1}\overline{z}_3+15z_2z_{-2}\overline{z}_3 + 2\sqrt{15} z_{-1}^2\overline{z}_1 +5\sqrt{15}z_{-2}^2\overline{z}_{-1} + 5\sqrt{2}z_0z_{-1}\overline{z}_2\\
&+5\sqrt{2} z_0z_{-2}\overline{z}_{1} + 15\sqrt{2}z_{-2}z_{-1}\overline{z}_0\\
\mathbf{Q}_{-2}(\mathbf{z}) =& 5 z_{-2}(5|z_{-3}|^2 +3|z_{-1}|^2 -3|z_1|^2 -8|z_2|^2 -5|z_3|^2) + 4\sqrt{30} z_{-1}z_0\overline{z}_1+25z_1z_{-1}\overline{z}_2\\
&+15 z_3z_{-3}\overline{z}_2 + 10\sqrt{15}z_{-1}z_{-3}\overline{z}_{-2} + 3\sqrt{30}z_{-1}^2\overline{z}_0+5\sqrt{2}z_1z_{-3}\overline{z}_0\\
&+5\sqrt{2}z_0z_1\overline{z}_3 + 15\sqrt{2}z_0z_{-3}\overline{z}_{-1}\\
\mathbf{Q}_{-1}(\mathbf{z}) =& z_{-1}(-5|z_{-3}|^2+15|z_{-2}|^2 -3|z_{-1}|^2 +12|z_0|^2 -16|z_1|^2-15|z_2|^2-25|z_3|^2) +24z_0^2\overline{z}_1\\
&+25z_2z_{-2}\overline{z}_1 -15z_3z_{-3}\overline{z}_1+4\sqrt{15}z_1z_{-3}\overline{z}_{-1} +2\sqrt{15}z_1^2\overline{z}_3 +5\sqrt{15}z_{-2}^2\overline{z}_{-3}\\
&+5\sqrt{2}z_2z_{-3}\overline{z}_0+5\sqrt{2}z_0z_2\overline{z}_3+15\sqrt{2}z_0z_{-3}\overline{z}_{-2} +6\sqrt{30}z_{-2}z_0\overline{z}_{-1}\\
&+4\sqrt{30}z_{-2}z_1\overline{z}_0 +4\sqrt{30} z_1z_0\overline{z}_2\\
\mathbf{Q}_{0}(\mathbf{z}) =& z_0(-20|z_{-3}|^2 +12|z_{-1}|^2 -12|z_0|^2 +12|z_1|^2 -20|z_3|^2) +48z_1z_{-1}\overline{z}_0+20z_3z_{-3}\overline{z}_0\\
&+ 15\sqrt{2}z_1z_2\overline{z}_3+15\sqrt{2}z_{-2}z_{-1}\overline{z}_{-3} +5\sqrt{2}z_3z_{-2}\overline{z}_1+5\sqrt{2}z_2z_{-3}\overline{z}_{-1}+5\sqrt{2}z_1z_{-3}\overline{z}_{-2} \\
&+ 5\sqrt{2}z_3z_{-1}\overline{z}_2 +4\sqrt{30}z_{-2}z_1\overline{z}_{-1}+4\sqrt{30}z_2z_{-1}\overline{z}_1 +3\sqrt{30}z_1^2\overline{z}_2 +3\sqrt{30}z_{-1}^2\overline{z}_{-2}
\end{align*}
\begin{align*}
\mathbf{R}_{-3}(\mathbf{z}) =& 3z_{-3}(3|z_{-3}|^2+3|z_{-2}|^2 +|z_{-1}|^2 -|z_1|^2-2|z_2|^2-3|z_3|^2)+3z_{-2}z_2\overline{z}_3+3\sqrt{2}z_0z_{-2}\overline{z}_1\\
&+3\sqrt{2}z_{-2}z_{-1}\overline{z}_0 +\sqrt{15}z_{-2}^2\overline{z}_{-1} +\sqrt{15}z_1z_{-2}\overline{z}_2\\
\mathbf{R}_{-2}(\mathbf{z}) =& z_{-2}(9|z_{-3}|^2+4|z_{-2}|^2 +7|z_{-1}|^2 -2|z_1|^2-4|z_2|^2-6|z_3|^2) +3z_{-3}z_3\overline{z}_2+3\sqrt{2}z_0z_{-3}\overline{z}_1\\
&+3\sqrt{2}z_1z_{-3}\overline{z}_0 +5z_{-1}z_1\overline{z}_2 + \sqrt{30}z_{-1}^2\overline{z}_0 +\sqrt{30}z_{-1}z_0\overline{z}_1 +\sqrt{15}z_{-3}z_2\overline{z}_1\\
&+\sqrt{15}z_2z_{-1}\overline{z}_3+2\sqrt{15}z_{-1}z_{-3}\overline{z}_{-2}\\
\mathbf{R}_{-1}(\mathbf{z}) =& z_{-1}(3|z_{-3}|^2+7|z_{-2}|^2 +|z_{-1}|^2 +6|z_0|^2 -|z_1|^2-2|z_2|^2-3|z_3|^2) +6z_0^2\overline{z}_1\\
&+ 3\sqrt{2}z_0z_{-3}\overline{z}_{-2}+3\sqrt{2}z_0z_2\overline{z}_3+2\sqrt{30}z_{-2}z_0\overline{z}_{-1}+\sqrt{30}z_{-2}z_1\overline{z}_{0}+\sqrt{30}z_0z_1\overline{z}_2\\
&+\sqrt{15}z_{-2}^2\overline{z}_{-3} +\sqrt{15}z_{-2}z_3\overline{z}_2 + 5z_{-2}z_2\overline{z}_1\\
\mathbf{R}_{0}(\mathbf{z}) =& 6z_0(|z_{-1}|^2 +|z_1|^2) +3\sqrt{2}z_{-3}z_1\overline{z}_{-2}+3\sqrt{2}z_{-1}z_3\overline{z}_2 +3\sqrt{2}z_1z_2\overline{z}_3 +3\sqrt{2}z_{-2}z_{-1}\overline{z}_{-3}\\
&+12z_1z_{-1}\overline{z}_0+\sqrt{30}z_{-2}z_1\overline{z}_{-1}+\sqrt{30}z_{-1}z_2\overline{z}_1+\sqrt{30} z_1^2\overline{z}_2 +\sqrt{30}z_{-1}^2\overline{z}_{-2}
\end{align*}

\section{Computing normal form coefficients for $l=1,2,3$}
\label{sec:computing_nfc_app}
\subsection{For $l=1$}
Suppose that we have a Hopf bifurcation, with a pair of simple eigenvalues $\lambda= i \omega, \overline{\lambda} = -i\omega$ for which $\mathcal{E}_1(i \omega)=0$ and $\mathcal{E}_l(i \hat{\omega})\neq 0$ for all $l\neq 1$ or $\omega \neq \hat{\omega}\in \mathbb{R}$. In this case the eigenvalue $\lambda$ has geometric multiplicity $3$, with corresponding eigenvectors $\psi_m(\theta)(\br) = e^{i\omega \theta}Y^m_1(\br) v$ for $m=-1,0,1$, with the normalisation $\|v\|_2=1$.

The homological equation \eqref{eq:homological_equation}, gives us a polynomial for which it is sufficient if we solve for the coefficients of the monomials $z_{-1}^2 \overline{z_{-1}}$ and $z_0^2 \overline{z_1}$ to compute all the normal form coefficients.
\begin{align*}
(i\omega - A^{\odot *})jh_{200\,100} =& \ell D^3G(\psi_{-1},\psi_{-1},\overline{\psi_{-1}})+\ell D^2G(\overline{\psi_{-1}},h_{200\,000})+2\ell D^2G(\psi_{-1},h_{100\,100}) -2g_{1,1} \psi_{-1}\\
(i\omega - A^{\odot *})jh_{020\,001} =& \ell D^3G(\psi_0,\psi_0,\overline{\psi_1})+\ell D^2G(\overline{\psi_1},h_{020\,000})+2\ell D^2G(\psi_0,h_{010\,001}) -2g_{1,2} \psi_{-1}\\
\end{align*}

To solve these equations, we first need to find $4$ different second order center manifold coefficients
\begin{align*}
(2i\omega - A^{\odot *})jh_{200\,000}(\theta) =& \ell D^2G(\psi_{-1},\psi_{-1}),\\    
- A^{\odot *}jh_{100\,100}(\theta) =& \ell D^2G(\psi_{-1},\overline{\psi_{-1}}),\\    
(2i\omega - A^{\odot *})jh_{020\,000}(\theta) =& \ell D^2G(\psi_0,\psi_0),\\ 
- A^{\odot *}h_{010\,001}(\theta) =& \ell D^2G(\psi_0,\overline{\psi_1}).\\    
\end{align*}

Using the techniques from \cite{spek_neural_2020}, we rewrite these equations in terms of the resolvent $\Delta^{-1}$
\begin{align*}
h_{200\,000}(\theta) =& e^{2i\omega\theta}\Delta^{-1}(2i \omega)D^2G(\psi_{-1},\psi_{-1}),\\    
h_{100\,100}(\theta) =& \Delta^{-1}(0)D^2G(\psi_{-1},\overline{\psi_{-1}}),\\    
h_{020\,000}(\theta) =& e^{2i\omega\theta}\Delta^{-1}(2i \omega)D^2G(\psi_0,\psi_0),\\ 
h_{010\,001}(\theta) =& \Delta^{-1}(0)D^2G(\psi_0,\overline{\psi_1}).\\    
\end{align*}

For the computation of $h_{200\,000}$ we use a decomposition into spherical harmonics, like in Section~\ref{sec:computing_nfc}
\begin{align*}
D^2G(\psi_{-1},\psi_{-1})(\br)=&S''(0) \int_\Omega g(\br\cdot \br',2 i \omega)(v*v) Y_1^{-1}(\br') Y_1^{-1}(\br') d\br' ,\\
\int_\Omega D^2G(\psi_{-1},\psi_{-1})(\br) \overline{Y_l^m(\br)}d\br =& S''(0) G_l(2 i \omega) (v*v)\int_\Omega Y_1^{-1}(\br') Y_1^{-1}(\br') \overline{Y_l^m(\br)}d\br\\
=& \begin{cases}
\sqrt{\frac{3}{10 \pi}} S''(0) G_2(2 i \omega)(v*v) & m=-2,n=2\\
0 & \text{otherwise}.
\end{cases}\end{align*}

The other coefficients of the center manifold we obtain similarly  
\begin{align*}
h_{200\,000}(\theta) =& \sqrt{\frac{3}{10 \pi}} e^{2i\omega \theta} S''(0) Q_2(2 i \omega)(v*v) Y_2^{-2},\\
h_{100\,100}(\theta) =& \frac{1}{2\sqrt{\pi}}S''(0)Q_0(0)(v*\overline{v}) Y_0^0-\frac{1}{2\sqrt{5\pi}}S''(0)Q_2(0)(v*\overline{v})Y_2^0,\\
h_{020\,000}(\theta) =& \frac{1}{2\sqrt{\pi}}e^{2i\omega \theta}S''(0)Q_0(2 i \omega)(v*v)Y_0^0+\frac{1}{\sqrt{5\pi}}e^{2i\omega \theta}S''(0)Q_2(2 i \omega)(v*v)Y_2^0,\\
h_{010\,001}(\theta) =& -\frac{1}{2}\sqrt{\frac{3}{5\pi}} S''(0)Q_2(0)(v*\overline{v})Y_2^{-1}.
\end{align*}

The normal form coefficients $g_{1,1}$ and $g_{1,2}$ can also be rewritten in terms of an integral in the complex plane and the resolvent $\Delta^{-1}$ using the techniques from \cite{spek_neural_2020} 
\begin{align*}
2g_{1,1} \psi_{-1}(\theta) =& \frac{1}{2\pi i} \oint_{\partial C_\lambda} e^{z\theta}\Delta^{-1}(z) (D^3G(\psi_{-1},\psi_{-1},\overline{\psi_{-1}})+D^2G(\overline{\psi_{-1}},h_{200\,000})+2D^2G(\psi_{-1},h_{100\,100})) \, dz,\\
-2g_{1,2} \psi_{-1}(\theta) =& \frac{1}{2\pi i} \oint_{\partial C_\lambda} e^{z\theta}\Delta^{-1}(z) (D^3G(\psi_0,\psi_0,\overline{\psi_1})+D^2G(\overline{\psi_1},h_{020\,000})+2D^2G(\psi_0,h_{010\,001})) \, dz.\\
\end{align*}

Decomposing each term into spherical harmonics, we find that
\begin{align*}
\Delta^{-1}(z)D^3G(\psi_{-1},\psi_{-1},\overline{\psi_{-1}})=& \frac{3}{10\pi} S'''(0) E_1^{-1}(z) G_1(i\omega)(v*v*\overline{v}) Y_1^{-1} \\
&- \frac{3}{10\sqrt{14}\pi} S'''(0) E_3^{-1}(z) G_3(i\omega)(v*v*\overline{v}) Y_3^0\\
\Delta^{-1}(z)D^2G(\overline{\psi_{-1}},h_{200\,000})=& \frac{3}{10\pi}(S''(0))^2 E_1^{-1}(z) G_1(i\omega) (\overline{v}*(Q_2(2 i \omega)(v*v)))Y_1^{-1}\\
&- \frac{3}{10\sqrt{14}\pi}(S''(0))^2 E_3^{-1}(z) G_3(i\omega) (\overline{v}*(Q_2(2 i \omega)(v*v)))Y_3^0\\
\Delta^{-1}(z)D^2G(\psi_{-1},h_{100\,100})=& \frac{1}{4\pi}(S''(0))^2 E_1^{-1}(z) G_1(i\omega) (v*(Q_0(0)(v*\overline{v})))Y_1^{-1}\\
&+\frac{1}{20\pi}(S''(0))^2 E_1^{-1}(z) G_1(i\omega) (v*(Q_2(0)(v*\overline{v})))Y_1^{-1}\\
&-\frac{3}{10\sqrt{14}\pi}(S''(0))^2 E_3^{-1}(z) G_3(i\omega) (v*(Q_2(0)(v*\overline{v})))Y_3^0\\
\Delta^{-1}(z)D^3G(\psi_0,\psi_0,\overline{\psi_1})=& -\frac{3}{20\pi} S'''(0) E_1^{-1}(z) G_1(i\omega)(v*v*\overline{v}) Y_1^{-1}\\
&-\frac{3}{5\sqrt{14}\pi} S'''(0) E_3^{-1}(z) G_3(i\omega)(v*v*\overline{v}) Y_3^0\\
\Delta^{-1}(z)D^2G(\overline{\psi_1},h_{020\,000})=& -\frac{1}{4\pi}(S''(0))^2 E_1^{-1}(z) G_1(i\omega) (\overline{v}*(Q_0(2 i \omega)(v*v)))Y_1^{-1}\\
&+\frac{1}{10\pi}(S''(0))^2 E_1^{-1}(z) G_1(i\omega) (\overline{v}*(Q_2(2 i \omega)(v*v)))Y_1^{-1}\\
&-\frac{3}{5\sqrt{14}\pi}(S''(0))^2 E_3^{-1}(z) G_3(i\omega) (\overline{v}*(Q_2(2 i \omega)(v*v)))Y_3^0\\
\Delta^{-1}(z)D^2G(\psi_{-1},h_{010\,001}))=& -\frac{3}{20\pi}(S''(0))^2 E_1^{-1}(z) G_1(i\omega) (v*(Q_2(0)(v*\overline{v})))Y_1^{-1}\\
&-\frac{3}{5\sqrt{14}\pi}(S''(0))^2 E_3^{-1}(z) G_3(i\omega) (v*(Q_2(0)(v*\overline{v})))Y_3^0.
\end{align*}

By the assumptions at the start of this section, $\lambda=i\omega$ is a simple eigenvalue. Hence, $\mathcal{E}_3(z)=\det E_3(z)$ is analytic at $z=i\omega$, but $\mathcal{E}_1(z)$ is not. Hence, all the terms with $E_3^{-1}(z)$ vanishes after applying the integral. For the terms with $E_1^{-1}(z)$, we use the adjugate matrix trick of Section~\ref{sec:computing_nfc}
\begin{align*}
g_{1,1}=& \frac{1}{20 \pi} \overline{v}^T \, \frac{adj(E_1(i\omega))}{\mathcal{E}'_1(i\omega)}G_1(i\omega)(3S'''(0)(v*v*\overline{v})\\
&+3(S''(0))^2(Q_2(2 i \omega)(v*v))*\overline{v}+(S''(0))^2((5Q_0(0)+Q_2(0))(v*\overline{v}))*v)\\[7pt]
g_{1,2}=& \frac{1}{40 \pi} \overline{v}^T \, \frac{adj(E_1(i\omega))}{\mathcal{E}'_1(i\omega)}G_1(i\omega)(3S'''(0)(v*v*\overline{v})\\
&+(S''(0))^2 ((5Q_0(2 i\omega)-2Q_2(2 i\omega))(v*v))*\overline{v} -6(S''(0))^2(Q_2(0)(v*\overline{v}))*v).
\end{align*}

\subsection{For $l=2$}
Suppose we have a Hopf bifurcation, with a pair of simple eigenvalues $\lambda= i \omega, \overline{\lambda} = -i\omega$ for which  $\mathcal{E}_2(i \omega)=0$ and $\mathcal{E}_l(i \hat{\omega})\neq 0$ for all $l\neq 2$  or $\omega \neq \hat{\omega}\in \mathbb{R}$. In this case the eigenvalue $\lambda$ has geometric multiplicity $5$, with corresponding eigenvectors $\psi_m(\theta)(\br) = e^{i\omega \theta}Y^m_2(\br) v$ for $m=-2,\cdots,2$, with the normalisation $\|v\|_2=1$. 

The homological equation \eqref{eq:homological_equation}, gives us a polynomial for which it is sufficient to solve for the coefficients of the monomials $z_{-2}^2 \overline{z_{-2}}, z_{-1} z_1 \overline{z_2}$ and $z_{-1} z_0 \overline{z_1}$ to compute all the normal form coefficients
\begin{align*}
(i\omega - A^{\odot *})jh_{20000\,10000}=& \ell D^3G(\psi_{-2},\psi_{-2},\overline{\psi_{-2}}) + \ell D^2G(\overline{\psi_{-2}}, h_{20000\,00000})\\
&+ 2 \ell D^2G(\psi_{-2},h_{10000\,10000}) - 2 g_{2,1}j \psi_{-2}\\
(i\omega - A^{\odot *})jh_{01010\,00001}=& \ell D^3G(\psi_{-1},\psi_1,\overline{\psi_{2}}) + \ell D^2G(\overline{\psi_{2}},h_{01010\,00000})\\
&+ \ell D^2G(\psi_1,h_{01000\,00001}) + \ell D^2G(\psi_{-1},h_{00010\,00001}) + 2 g_{2,2}j \psi_{-2}\\
(i\omega - A^{\odot *})jh_{01100\,00010}=& \ell D^3G(\psi_{-1},\psi_0,\overline{\psi_1}) + \ell D^2G(\overline{\psi_{1}},h_{01100\,00000})\\
&+ \ell D^2G(\psi_{0},h_{01000\,00010}) + \ell D^2G(\psi_{-1},h_{00100\,00010}) -  g_{2,3}j \psi_{-2}
\end{align*}

To solve these equations, we first need to find $8$ different second order center manifold coefficients
\begin{align*}
(2i\omega - A^{\odot *})jh_{20000\,00000}=& \ell D^2G(\psi_{-2},\psi_{-2})\\
- A^{\odot *}jh_{10000\,10000}=& \ell D^2G(\psi_{-2},\overline{\psi_{-2}})\\
(2i\omega - A^{\odot *})jh_{01010\,00000}=& \ell D^2G(\psi_{-1},\psi_{1})\\
- A^{\odot *}jh_{01000\,00001}=& \ell D^2G(\psi_{-1},\overline{\psi_{2}})\\
- A^{\odot *}jh_{00010\,00001}=& \ell D^2G(\psi_{1},\overline{\psi_{2}})\\
(2i\omega - A^{\odot *})jh_{01100\,00000}=& \ell D^2G(\psi_{-1},\psi_{0})\\
- A^{\odot *}jh_{01000\,00010}=& \ell D^2G(\psi_{-1},\overline{\psi_{1}})\\
- A^{\odot *}jh_{00100\,00010}=& \ell D^2G(\psi_{0},\overline{\psi_{1}})
\end{align*}

Following the techniques from \cite{spek_neural_2020},
\begin{align*}
h_{20000\,00000}(\theta) =& e^{2i\omega\theta}\Delta^{-1}(2i \omega)D^2G(\psi_{-2},\psi_{-2})= \sqrt{\frac{5}{14\pi}} e^{2i\omega\theta} S''(0) Q_4(2i\omega)(v*v) Y_4^{-4}\\
h_{10000\,10000}(\theta) =& \Delta^{-1}(0)D^2G(\psi_{-2},\overline{\psi_{-2}})\\
=& \frac{1}{2\sqrt{\pi}}S''(0) Q_0(0)(v*\overline{v}) Y_0^{0}-\frac{1}{7}\sqrt{\frac{5}{\pi}}S''(0) Q_2(0)(v*\overline{v}) Y_2^{0}+\frac{1}{14\sqrt{\pi}}S''(0) Q_4(0)(v*\overline{v}) Y_4^{0}\\
h_{01010\,00000}(\theta) =& e^{2i\omega\theta}\Delta^{-1}(2i \omega)D^2G(\psi_{-1},\psi_{1}),\\
=& -\frac{1}{2\sqrt{\pi}} e^{2i\omega\theta}S''(0) Q_0(2i \omega)(v*v) Y_0^{0}-\frac{1}{14}\sqrt{\frac{5}{\pi}} e^{2i\omega\theta}S''(0) Q_2(2i \omega)(v*v) Y_2^{0}\\
&+\frac{2}{7\sqrt{\pi}} e^{2i\omega\theta}S''(0) Q_4(2i \omega)(v*v) Y_4^{0}\\
h_{01000\,00001}(\theta) =& \Delta^{-1}(0)D^2G(\psi_{-1},\overline{\psi_{2}})= \frac{1}{2}\sqrt{\frac{5}{7\pi}}S''(0) Q_4(0)(v*\overline{v}) Y_4^{-3}\\
h_{00010\,00001}(\theta) =& \Delta^{-1}(0)D^2G(\psi_{1},\overline{\psi_{2}})= -\frac{1}{7}\sqrt{\frac{15}{2\pi}}S''(0) Q_2(0)(v*\overline{v}) Y_2^{-1}+\frac{1}{14}\sqrt{\frac{5}{\pi}}S''(0) Q_4(0)(v*\overline{v}) Y_4^{-1}\\
h_{01100\,00000}(\theta) =& e^{2i\omega\theta}\Delta^{-1}(2i \omega)D^2G(\psi_{-1},\psi_{0})= \frac{1}{14}\sqrt{\frac{5}{\pi}}e^{2i\omega\theta}S''(0) Q_2(2i \omega)(v*v) Y_2^{-1}\\
&+\frac{1}{7}\sqrt{\frac{15}{2\pi}}e^{2i\omega\theta}S''(0) Q_4(2i \omega)(v*v) Y_4^{-1}\\
h_{01000\,00010}(\theta) =& \Delta^{-1}(0)D^2G(\psi_{-1},\overline{\psi_{1}})= -\frac{1}{7}\sqrt{\frac{15}{2\pi}}S''(0) Q_2(0)(v*\overline{v}) Y_2^{-2}-\frac{1}{7}\sqrt{\frac{10}{\pi}}S''(0) Q_4(0)(v*\overline{v}) Y_4^{-2}\\
h_{00100\,00010}(\theta) =& \Delta^{-1}(0)D^2G(\psi_{0},\overline{\psi_{1}})= -\frac{1}{14}\sqrt{\frac{5}{\pi}}S''(0) Q_2(0)(v*\overline{v}) Y_2^{-1}-\frac{1}{7}\sqrt{\frac{15}{2\pi}}S''(0) Q_4(0)(v*\overline{v}) Y_4^{-1}
\end{align*}

For the normal form coefficients, we get
\begin{align*}
2 g_{2,1}j \psi_{-2}(\theta) =& \frac{1}{2\pi i} \oint_{\partial C_\lambda} e^{z\theta}\Delta^{-1}(z) ( D^3G(\psi_{-2},\psi_{-2},\overline{\psi_{-2}})+ D^2G(\overline{\psi_{-2}}, h_{20000\,00000}) \\
&+ 2 D^2G(\psi_{-2},h_{10000\,10000}) )\, dz,\\
-2 g_{2,2}j \psi_{-2}(\theta) =& \frac{1}{2\pi i} \oint_{\partial C_\lambda} e^{z\theta}\Delta^{-1}(z) ( D^3G(\psi_{-1},\psi_1,\overline{\psi_{2}}) + D^2G(\overline{\psi_{2}},h_{01010\,00000}) \\
&+  D^2G(\psi_1,h_{01000\,00001}) +  D^2G(\psi_{-1},h_{00010\,00001}) ) \, dz,\\
g_{2,3}j \psi_{-2}(\theta) =& \frac{1}{2\pi i} \oint_{\partial C_\lambda} e^{z\theta}\Delta^{-1}(z) ( D^3G(\psi_{-1},\psi_0,\overline{\psi_1}) + D^2G(\overline{\psi_{1}},h_{01100\,00000})\\
&+  D^2G(\psi_{0},h_{01000\,00010}) +  D^2G(\psi_{-1},h_{00100\,00010}) )\, dz.
\end{align*}

Using the same procedure as above, where we use that $Q_l(z)$ is analytic at $z=\lambda$ if and only if $l=2$
\begin{align*}
g_{2,1} =& \frac{1}{196 \pi} \overline{v}^T \, \frac{adj(E_2(i\omega))}{\mathcal{E}'_2(i\omega)}G_2(i\omega)(35S'''(0)(v*v*\overline{v})+35(S''(0))^2(Q_4(2i\omega)(v*v))*\overline{v}\\
&+(S''(0))^2((49Q_0(0)+20Q_2(0)+Q_4(0))(v*\overline{v}))*v)\\
g_{2,2} =& \frac{1}{392 \pi} \overline{v}^T \, \frac{adj(E_2(i\omega))}{\mathcal{E}'_2(i\omega)}G_2(i\omega)(35S'''(0)(v*v*\overline{v})\\
&+(S''(0))^2((49Q_0(2i\omega)-10Q_2(2i\omega)-4Q_4(2i\omega))(v*v))*\overline{v}\\
&+(S''(0))^2((30Q_2(0)+40Q_4(0))(v*\overline{v}))*v)\\
g_{2,3} =& \frac{5}{98 \pi}\sqrt{\frac{3}{2}} \overline{v}^T \, \frac{adj(E_2(i\omega))}{\mathcal{E}'_2(i\omega)}G_2(i\omega)((S''(0))^2((-Q_2(2i\omega)+Q_4(2i\omega))(v*v))*\overline{v}\\
&+(S''(0))^2((Q_2(0)-Q_4(0))(v*\overline{v}))*v).
\end{align*}
One interesting note is that in the decomposition of  $D^3G(\psi_{-1},\psi_0,\overline{\psi_1})$ the coefficient corresponding to $Y_2^m$ is zero for all $m$. This makes $\Delta^{-1}(z) D^3G(\psi_{-1},\psi_0,\overline{\psi_1})$ analytic at $z=\lambda$, so this term vanishes. Therefore, when $S''(0)=0$, we have a degeneracy with $g_{2,3}=0$. 

\subsection{For $l=3$}
Suppose we have a Hopf bifurcation, with a pair of simple eigenvalues $\lambda= i \omega, \overline{\lambda} = -i\omega$ for which $\mathcal{E}_3(i \omega)=0$ and $\mathcal{E}_l(i \hat{\omega})\neq 0$ for all $l\neq 3$ or $\omega \neq \hat{\omega}\in \mathbb{R}$. In this case the eigenvalue has geometric multiplicity $7$, with corresponding eigenvectors $\psi_m(\theta)(\br) = e^{i\omega \theta}Y^m_3(\br) v$ for $m=-3,\cdots,3$, with the normalisation $\|v\|_2=1$. 

The homological equation \eqref{eq:homological_equation}, gives us a polynomial for which it is sufficient is we solve for the coefficients of the monomials $z_{-2} z_0 \overline{z_0}, z_0^2 \overline{z_2}, z_0 z_{-1} \overline{z_2}$ and $z_{-1} z_2 \overline{z_3}$ to compute all the normal form coefficients
\begin{align*}
(i\omega - A^{\odot *})jh_{0101000\,0001000}=& \ell D^3G(\psi_{-2},\psi_0,\overline{\psi_0}) + \ell D^2G(\overline{\psi_{0}},h_{0101000\,0000000})\\
&+ \ell D^2G(\psi_{0},h_{0100000\,0001000}) + \ell D^2G(\psi_{-2},h_{0001000\,0001000}) -  g_{3,1}j \psi_{-2}\\
(i\omega - A^{\odot *})jh_{0002000\,0000010}=& \ell D^3G(\psi_0,\psi_0,\overline{\psi_2}) + \ell D^2G(\overline{\psi_{2}},h_{0002000\,0000000})\\
&+ 2\ell D^2G(\psi_{0},h_{0001000\,0000010}) - 2 g_{3,2}j \psi_{-2}\\
(i\omega - A^{\odot *})jh_{0011000\,0000010}=& \ell D^3G(\psi_{-1},\psi_0,\overline{\psi_2}) + \ell D^2G(\overline{\psi_{2}},h_{0011000\,0000000})\\
&+ \ell D^2G(\psi_{0},h_{0010000\,0000010}) + \ell D^2G(\psi_{-1},h_{0001000\,0000010}) - 5\sqrt{2} g_{3,3}j \psi_{-3}\\
(i\omega - A^{\odot *})jh_{0010010\,0000001}=& \ell D^3G(\psi_{-1},\psi_2,\overline{\psi_3}) + \ell D^2G(\overline{\psi_{3}},h_{0010010\,0000000})\\
&+ \ell D^2G(\psi_{2},h_{0010000\,0000001}) + \ell D^2G(\psi_{-1},h_{0000010\,0000001}) - \sqrt{15} g_{3,4}j \psi_{-2}.
\end{align*}

To solve these equation, we first need to find $15$ different second order center manifold coefficients
\begin{align*}
(2i\omega - A^{\odot *})jh_{0101000\,0000000}=& \ell D^2G(\psi_{-2},\psi_{0})\\
- A^{\odot *}jh_{0100000\,0001000}=& \ell D^2G(\psi_{-2},\overline{\psi_{0}})\\
- A^{\odot *}jh_{0001000\,0001000}=& \ell D^2G(\psi_{0},\overline{\psi_{0}})\\
(2i\omega - A^{\odot *})jh_{0002000\,0000000}=& \ell D^2G(\psi_{0},\psi_{0})\\
- A^{\odot *}jh_{0001000\,0000010}=& \ell D^2G(\psi_{0},\overline{\psi_{2}})\\
(2i\omega - A^{\odot *})jh_{0011000\,0000000}=& \ell D^2G(\psi_{-1},\psi_{0})\\
- A^{\odot *}jh_{0010000\,0000010}=& \ell D^2G(\psi_{-1},\overline{\psi_{2}})\\
- A^{\odot *}jh_{0001000\,0000010}=& \ell D^2G(\psi_{0},\overline{\psi_{2}})\\
(2i\omega - A^{\odot *})jh_{0010010\,0000000}=& \ell D^2G(\psi_{-1},\psi_{2})\\
- A^{\odot *}jh_{0010000\,0000001}=& \ell D^2G(\psi_{-1},\overline{\psi_{3}})\\
- A^{\odot *}jh_{0000010\,0000001}=& \ell D^2G(\psi_{2},\overline{\psi_{3}}).
\end{align*}

We get that
\begin{align*}
h_{0101000\,0000000}(\theta)=& -\frac{1}{3\sqrt{\pi}}e^{2i\omega\theta}S''(0) Q_2(2i\omega)(v*\overline{v})Y_2^{-2} - \frac{1}{22}\sqrt{\frac{3}{\pi}}e^{2i\omega\theta}S''(0) Q_4(2i\omega)(v*\overline{v})Y_4^{-2}\\
&+ \frac{10}{33}\sqrt{\frac{14}{13\pi}}e^{2i\omega\theta}S''(0) Q_6(2i\omega)(v*\overline{v})Y_6^{-2}\\
h_{0100000\,0001000}(\theta)=&  -\frac{1}{3\sqrt{\pi}}S''(0) Q_2(0)(v*\overline{v})Y_2^{-2}  - \frac{1}{22}\sqrt{\frac{3}{\pi}} S''(0) Q_4(0)(v*\overline{v})Y_4^{-2}\\
&+ \frac{10}{33}\sqrt{\frac{14}{13\pi}}S''(0) Q_6(0)(v*\overline{v})Y_6^{-2}\\
h_{0001000\,0001000}(\theta)=& \frac{1}{2\sqrt{\pi}}S''(0) Q_0(0)(v*\overline{v})Y_0^{0} + \frac{2}{3\sqrt{5\pi}}S''(0) Q_2(0)(v*\overline{v})Y_2^{0} + \frac{6}{11\sqrt{\pi}}S''(0) Q_4(0)(v*\overline{v})Y_4^{0}\\
&+\frac{50}{33\sqrt{13\pi}}S''(0) Q_6(0)(v*\overline{v})Y_6^{0}\\
h_{0002000\,0000000}(\theta)=& \frac{1}{2\sqrt{\pi}}e^{2i\omega\theta}S''(0) Q_0(2i\omega)(v*\overline{v})Y_0^{0} + \frac{2}{3\sqrt{5\pi}}e^{2i\omega\theta}S''(0) Q_2(2i\omega)(v*\overline{v})Y_2^{0}\\
&+ \frac{3}{11\sqrt{\pi}}e^{2i\omega\theta}S''(0) Q_4(2i\omega)(v*\overline{v})Y_4^{0} + \frac{50}{33\sqrt{13\pi}}e^{2i\omega\theta}S''(0) Q_6(2i\omega)(v*\overline{v})Y_6^{0}\\
h_{0001000\,0000010}(\theta)=& -\frac{1}{3\sqrt{\pi}}S''(0) Q_2(0)(v*\overline{v})Y_2^{-2} + \frac{1}{22}\sqrt{\frac{3}{\pi}}S''(0) Q_4(0)(v*\overline{v})Y_4^{-2}\\
&+ \frac{10}{33}\sqrt{\frac{14}{13\pi}}S''(0) Q_6(0)(v*\overline{v})Y_6^{-2}\\
h_{0011000\,0000000}(\theta)=& \frac{1}{3\sqrt{10\pi}}e^{2i\omega\theta}S''(0) Q_2(2i\omega)(v*\overline{v})Y_2^{-1} - \frac{1}{22}\sqrt{\frac{15}{\pi}}e^{2i\omega\theta}S''(0) Q_4(2i\omega)(v*\overline{v})Y_4^{-1}\\
&+ \frac{25}{33}\sqrt{\frac{7}{26\pi}}e^{2i\omega\theta}S''(0) Q_6(2i\omega)(v*\overline{v})Y_6^{-1}\\
h_{0010000\,0000010}(\theta)=& \frac{1}{11}\sqrt{\frac{7}{2\pi}}S''(0) Q_4(0)(v*\overline{v})Y_4^{-3} + \frac{5}{11}\sqrt{\frac{21}{26\pi}}S''(0) Q_6(0)(v*\overline{v})Y_6^{-3}\\
h_{0001000\,0000010}(\theta)=& -\frac{1}{3\sqrt{\pi}}S''(0) Q_2(0)(v*\overline{v})Y_2^{-2} + \frac{1}{22}\sqrt{\frac{3}{\pi}}S''(0) Q_4(0)(v*\overline{v})Y_4^{-2}\\
&+ \frac{10}{33}\sqrt{\frac{14}{13\pi}}S''(0) Q_6(0)(v*\overline{v})Y_6^{-2}\\
h_{0010010\,0000000}(\theta)=&  -\frac{1}{2\sqrt{3\pi}}e^{2i\omega\theta}S''(0) Q_2(2i\omega)(v*\overline{v})Y_2^{1} - \frac{2}{11}\sqrt{\frac{2}{\pi}}e^{2i\omega\theta}S''(0) Q_4(2i\omega)(v*\overline{v})Y_4^{1}\\
&+ \frac{5}{22}\sqrt{\frac{35}{39\pi}}e^{2i\omega\theta}S''(0) Q_6(2i\omega)(v*\overline{v})Y_6^{1}\\
h_{0010000\,0000001}(\theta)=& \frac{1}{11}\sqrt{\frac{21}{2\pi}}S''(0) Q_4(0)(v*\overline{v})Y_4^{-4} - \frac{5}{11}\sqrt{\frac{35}{78\pi}}S''(0) Q_6(0)(v*\overline{v})Y_6^{-4}\\
h_{0000010\,0000001}(\theta)=& -\frac{1}{6}\sqrt{\frac{5}{\pi}}S''(0) Q_2(0)(v*\overline{v})Y_2^{-1} + \frac{1}{11}\sqrt{\frac{15}{2\pi}}S''(0) Q_4(0)(v*\overline{v})Y_4^{-1}\\
&- \frac{5}{66}\sqrt{\frac{7}{13\pi}}S''(0) Q_6(0)(v*\overline{v})Y_6^{-1}.
\end{align*}

Following the techniques from \cite{spek_neural_2020}.
\begin{align*}
g_{3,1}j \psi_{-2}(\theta) =& \frac{1}{2\pi i} \oint_{\partial C_\lambda} e^{z\theta}\Delta^{-1}(z)  (D^3G(\psi_{-2},\psi_0,\overline{\psi_0})+D^2G(\overline{\psi_{0}},h_{0101000\,0000000})\\
&+  D^2G(\psi_{0},h_{0100000\,0001000}) +  D^2G(\psi_{-2},h_{0001000\,0001000}))\, dz,\\
2 g_{3,2}j \psi_{-2}(\theta) =& \frac{1}{2\pi i} \oint_{\partial C_\lambda} e^{z\theta}\Delta^{-1}(z)  (D^3G(\psi_0,\psi_0,\overline{\psi_2})+D^2G(\overline{\psi_{2}},h_{0002000\,0000000}) \\
&+ 2D^2G(\psi_{0},h_{0001000\,0000010})) \, dz,\\
5\sqrt{2} g_{3,3}j \psi_{-3}(\theta) =& \frac{1}{2\pi i} \oint_{\partial C_\lambda} e^{z\theta}\Delta^{-1}(z)  (D^3G(\psi_{-1},\psi_0,\overline{\psi_2})+D^2G(\overline{\psi_{2}},h_{0011000\,0000000}) \\
&+ D^2G(\psi_{0},h_{0010000\,0000010}) + D^2G(\psi_{-1},h_{0001000\,0000010}))\, dz,\\
\sqrt{15} g_{3,4}j \psi_{-2}(\theta) =& \frac{1}{2\pi i} \oint_{\partial C_\lambda} e^{z\theta}\Delta^{-1}(z) (D^3G(\psi_{-1},\psi_2,\overline{\psi_3})+D^2G(\overline{\psi_{3}},h_{0010010\,0000000}) \\
&+ D^2G(\psi_{2},h_{0010000\,0000001}) + D^2G(\psi_{-1},h_{0000010\,0000001}))\, dz.
\end{align*}

Using the same procedure as above, where we use that $Q_l(z)$ is analytic at $z=\lambda$ if and only if $l=3$
\begin{align*}
g_{3,1}=& \frac{1}{56628 \pi} \overline{v}^T \, \frac{adj(E_3(i\omega))}{\mathcal{E}'_3(i\omega)}G_3(i\omega)(12243S'''(0)(v*v*\overline{v})\\
&+(S''(0))^2((6292Q_2(2i\omega)+351Q_4(2i\omega)+5600Q_6(2i\omega))(v*v))*\overline{v}\\
&+(S''(0))^2((14157Q_0(0)+6292Q_2(0)-4563Q_4(0)+8600Q_6(0))(v*\overline{v}))*v)\\
g_{3,2} =& \frac{1}{113256 \pi} \overline{v}^T \, \frac{adj(E_3(i\omega))}{\mathcal{E}'_3(i\omega)}G_3(i\omega)(12243S'''(0)(v*v*\overline{v})\\
&+(S''(0))^2((14157Q_0(2i\omega)-4914Q_4(2i\omega)+3000Q_6(2i\omega))(v*v))*\overline{v}\\
&+(S''(0))^2((12584Q_2(0)+702Q_4(0)+11200Q_6(0))(v*\overline{v}))*v)\\
g_{3,3}=& \frac{1}{283140 \pi} \overline{v}^T \, \frac{adj(E_3(i\omega))}{\mathcal{E}'_3(i\omega)}G_3(i\omega)(693S'''(0)(v*v*\overline{v})\\
&+(S''(0))^2((1573Q_2(2i\omega)-1755Q_4(2i\omega)+875Q_6(2i\omega))(v*v))*\overline{v}\\
&+(S''(0))^2((3146Q_2(0)-3510Q_4(0)+1750Q_6(0))(v*\overline{v}))*v)\\
g_{3,4} =& -\frac{1}{56628 \pi} \overline{v}^T \, \frac{adj(E_3(i\omega))}{\mathcal{E}'_3(i\omega)}G_3(i\omega)(462S'''(0)(v*v*\overline{v})\\
&+(S''(0))^2((1573Q_2(2i\omega)-936Q_4(2i\omega)-175Q_6(2i\omega))(v*v))*\overline{v}\\
&+(S''(0))^2((1573Q_2(0)-2574Q_4(0)+1925Q_6(0))(v*\overline{v}))*v).
\end{align*}
\endgroup

\bibliographystyle{plainnat} 
\bibliography{NeuralField}

\end{document}